\documentclass[11pt]{amsart}

\usepackage{amsfonts,amsmath,amssymb,amsthm}
\allowdisplaybreaks[4]

\usepackage[colorlinks,linkcolor=red,citecolor=green]{hyperref}
\usepackage{makecell}
\usepackage{graphics}
\usepackage{tikz}
\newcommand{\circled}[1]{
  \tikz[baseline=(char.base)]{
    \node[shape=circle,draw,inner sep=1pt] (char) {#1};
  }
}

\newtheorem{theorem}{Theorem}
\newtheorem{lemma}[theorem]{Lemma}
\newtheorem{remark}[theorem]{Remark}
\newtheorem{proposition}[theorem]{Proposition}

\newtheorem{corollary}[theorem]{Corollary}

\numberwithin{theorem}{section}
\numberwithin{equation}{section}

\begin{document}

\setlength{\baselineskip}{1.1\baselineskip}

\title[Power convexity of solutions to complex Monge-Amp\`{e}re equation]{Power convexity of solutions to complex Monge-Amp\`{e}re equation in $\mathbb{C}^2$}

\author{Wei Zhang}
\address{School of Mathematics and Statistics\\
Lanzhou University\\
Lanzhou, 730000, Gansu Province, China.}
\email{zhangw@lzu.edu.cn}

\author{Qi Zhou}
\address{School of Mathematics and Statistics\\
Lanzhou University\\
Lanzhou, 730000, Gansu Province, China.}
\email{zhouqi2025@lzu.edu.cn, zhouqimath20@lzu.edu.cn}

\date{\today}

\maketitle

\begin{abstract}
The convexity of solutions to boundary value problems for fully nonlinear elliptic partial differential equations (such as real or complex $k$-Hessian equations) is a challenging topic. In this paper, we establish the power convexity of solutions to the Dirichlet problem for the complex Monge-Amp\`{e}re equation on a bounded, smooth, strictly convex domain in $\mathbb{C}^2$. Our approach is based on the constant rank theorem and the deformation process. A key obstacle in establishing the constant rank theorem lies in the construction of suitable auxiliary functions for deriving the associated differential inequalities. To address this difficulty, we refine the auxiliary function introduced by Bian and Guan \cite{Bian-Guan2009}.
\end{abstract}



\tableofcontents

\section{Introduction}

\subsection{Main result} In 1971, Makar-Limanov \cite{Makar-Limanov1971} studied the Dirichlet problem for the Laplace equation
\begin{align}\label{Equ-Laplace}
\begin{cases}
\Delta u=-1 & \mathrm{in}\ \Omega, \\
\hspace{3mm} u=0 & \mathrm{on}\ \partial\Omega,
\end{cases}
\end{align}
where $\Omega\subset\mathbb{R}^2$ is a bounded, strictly convex domain with smooth boundary. He proved that the function $\sqrt{u}$ is strictly concave in $\Omega$. Recall that the unique solution $u$ of the Dirichlet problem \eqref{Equ-Laplace} is commonly referred to as the torsion function of $\Omega$. That is to say, Makar-Limanov \cite{Makar-Limanov1971} established the $1/2$-concavity of the torsion function of a planar bounded convex domain. Since then, many authors have investigated the convexity of solutions to various kinds of elliptic partial differential equations (PDEs); see Section \ref{Section 1.2} below.

In this work, we will explore the convexity of solutions of the Dirichlet problem for the complex Monge-Amp\`{e}re equation. The complex Monge-Amp\`{e}re equation is a fully nonlinear elliptic PDE of the form $$\det(u_{i\bar{j}})=f, $$ where $u$ is a plurisubharmonic function, and $(u_{i\bar{j}})$ is its complex Hessian in local holomorphic coordinates. This equation plays a fundamental role in complex geometry and fully nonlinear PDE theory. 

Let us first recall the solvability of the Dirichlet problem for complex Monge-Amp\`{e}re equation on bounded domain in $\mathbb{C}^n$. Assume that $\Omega\subset\mathbb{C}^n$ is a bounded strongly pseudoconvex domain with smooth boundary, and that $0<f\in C^{\infty}(\bar\Omega)$, $\varphi\in C^{\infty}(\partial\Omega)$. Then Caffarelli, Kohn, Nirenberg and Spruck \cite{CKNS1985} proved that the Dirichlet problem
\begin{equation}\label{Equ-CMA}
\begin{cases}
\det(u_{i\bar{j}})=f & \mathrm{in}\ \Omega, \\
\hspace{11mm} u=\varphi & \mathrm{on}\ \partial\Omega,
\end{cases}
\end{equation}
admits a unique strictly plurisubharmonic solution $u\in C^{\infty}(\bar\Omega)$. 

Our goal in this paper is to investigate the power convexity of the unique solution $u$ to problem \eqref{Equ-CMA} when $f=1$, $\varphi=0$ and $n=2$. Precisely, our main result is as follows.
\begin{theorem}\label{Theorem-Main}
Let $\Omega\subset\mathbb{C}^2$ be a bounded, strictly convex domain with smooth boundary. Suppose that $u\in C^\infty(\bar{\Omega})$ is the unique solution of Dirichlet problem 
\begin{equation}\label{Equ-CMA1}
\begin{cases}
\det(u_{i\bar{j}})=1 & \mathrm{in}\ \Omega, \\
\hspace{11mm} u=0 & \mathrm{on}\ \partial\Omega.
\end{cases}
\end{equation}
Then the function $v=-\sqrt{-u/2}$ is strictly convex in $\Omega$.
\end{theorem}

\begin{remark}
To the best of our knowledge, this is the first known result addressing the power convexity of solutions to complex Monge-Amp\`{e}re equations. For $n\geq 3$, whether the unique solution of Dirichlet problem \eqref{Equ-CMA1} exhibits the power convexity property remains an open question.
\end{remark}

\subsection{Previous works}\label{Section 1.2}

Convexity has long been a fundamental and intriguing topic in the study of elliptic PDEs. Over the decades, many authors have made significant contributions to this field. Below, we provide a brief overview of the historical developments on the convexity of solutions to elliptic PDEs.

As we mentioned earlier, the $1/2$-concavity of the two-dimensional torsion function was first studied by Makar-Limanov \cite{Makar-Limanov1971}. His key idea was to introduce the auxiliary function $$P=v^4\det(\nabla^2v), $$ where $v=\sqrt{u}$ and $u$ satisfies the Dirichlet problem \eqref{Equ-Laplace}. He then could show that $P$ is superharmonic in $\Omega$. In 1976, Brascamp and Lieb \cite{Brascamp-Lieb1976} proved the log-concavity of the first Dirichlet eigenfunction of the Laplacian of a bounded convex domain by parabolic method. Acker, Payne and Philippin \cite{Acker-PP1981} later provided an alternative proof of Brascamp and Lieb's result along the route of Makar-Limanov \cite{Makar-Limanov1971} for the two-dimensional case. 

In the early 1980s,	Korevaar \cite{Korevaar1983a, Korevaar1983b} made a breakthrough in studying the convexity of solutions to linear and quasilinear elliptic equations. Specifically, for a bounded convex domain $\Omega\subset\mathbb{R}^n$ and a solution $u$, he introduced the convexity test function $$\mathcal{C}(x,y,\lambda)=u((1-\lambda)x+\lambda y)-(1-\lambda)u(x) -\lambda u(y)$$ for $(x, y, \lambda)\in\bar{\Omega} \times\bar{\Omega}\times[0, 1]$. Clearly, the nonnegativity of $\mathcal{C}(x, y, \lambda)$ is equivalent to the concavity of $u$. This method is commonly known as Korevaar's concavity maximum principle. Using this, Korevaar reproved the result of Brascamp and Lieb \cite{Brascamp-Lieb1976} (see also Caffarelli and Spruck \cite{Caffarelli-Spruck1982}). Later, Kennington \cite{Kennington1985} and Kawohl \cite{Kawohl1986} presented an improved version that could be used to establish the $1/2$-concavity of the torsion function in high dimensions. For problems involving the $p$-Laplacian or the Finsler Laplacian, due to the lack of regularity in solutions, Korevaar's concavity maximum principle cannot be applied directly. However, by suitable approximation arguments, power convexity results have been derived in \cite{Borrelli-MS2024-AdvCalc, Mosconi-RS2024-DCDS, Sakaguchi1987}, among others. In 2021, Langford and Scheuer \cite{Langford-Scheuer2021} extended Korevaar's concavity maximum principle to spherical settings, and proved the convexity of solutions to certain fully nonlinear degenerate elliptic equations under appropriate conditions.

In 1997, Alvarez, Lasry and Lions \cite{Alvarez-LL1997} developed a new method, the convex envelope method, to establish the convexity of viscosity solution of general fully nonlinear degenerate elliptic equation $$F(\nabla^2v, \nabla v, v, x)=0$$ in a bounded convex domain $\Omega$ of $\mathbb{R}^n$. In particular, they gave alternative proof of previous results, including the $1/2$-concavity of the torsion function and the log-concavity of the first Dirichlet eigenfunction of the Laplacian. This method was also extend to handle the power convexity of certain elliptic boundary value problems directly (still working on the solution before transformation) as well as the power convexity of solutions of semilinear elliptic systems in \cite{Bianchini-Salani2015, Ishige-NS2016}. Recently, one more application of convex envelope method was ascribed by Crasta and Fragal\`{a} \cite{Crasta-Fragala-2015}, who proved that the unique viscosity solution $u$ to Dirichlet problem 
\begin{align*}
\begin{cases}
\Delta_\infty u=-1 & \mathrm{in}\ \Omega, \\
\hspace{0.64cm}u=0 & \mathrm{on}\ \partial\Omega
\end{cases}
\end{align*}
is $3/4$-concave and in $C^1(\Omega)$. Here $\Delta_{\infty}$ denotes the usual infinity Laplacian.

Another powerful tool to produce convex solution is the so-called microscopic convexity principle (also called constant rank theorem). This approach was first discovered by Caffarelli and Friedman \cite{Caffarelli-Friedman1985} for semilinear elliptic equations in two-dimensional case. Combining the deformation process, they established the $1/2$-concavity of the torsion function and the log-concavity of the first Dirichlet eigenfunction of the Laplacian in dimension two. Later, Korevaar and Lewis \cite{Korevaar-Lewis1987} proved  the analogous results in high dimensions. Similar idea also appeared at the same time in Singer, Wong, Yau and Yau \cite{Singer-WYY1985} to deal with the log-concavity of the first Dirichlet eigenfunction of the Laplacian.

For fully nonlinear elliptic equations, the constant rank theorem remains highly effective in dealing with problems related to convexity. One of the earliest results in this direction was obtained by Ma and Xu \cite{Ma-Xu2008}, who studied the Dirichlet problem $$
\begin{cases}
\sigma_2(\nabla^2u)=1 & \mathrm{in}\ \Omega, \\
\hspace{12mm}u=0 & \mathrm{on}\ \partial\Omega,
\end{cases}$$
where $\Omega\subset\mathbb{R}^3$ is a bounded, convex domain with smooth boundary, and demonstrated that $-\sqrt{-u}$ is strictly convex in $\Omega$. In the same way, Liu, Ma and Xu \cite{Liu-Ma-Xu2010} investigated the $\log$-concavity of the eigenfunction of the $2$-Hessian operator and derived a Brunn-Minkowski inequality for the $2$-Hessian eigenvalue. See Salani \cite{Salani2012} for new proofs of the power convexity results in \cite{Liu-Ma-Xu2010, Ma-Xu2008} as well as discussions on Brunn-Minkowski inequalities of the $2$-torsion functional. Much more recently, several lower dimensional convexity results of solutions to fully nonlinear elliptic boundary value problems appeared via constant rank theorems, for example \cite{Chen-Jia-Xiong2025, Huang2019, Zhang-Zhou2023}. For historical developments on constant rank theorems of convex solutions of  fully nonlinear elliptic PDEs and their applications in geometry, one can see \cite{Bhattacharya-Shankar2024, Bian-Guan2009, Bryan-Ivaki-Scheuer2023, Caffarelli-Guan-Ma2007,Ogden-Yuan2025, Szekelyhidi-Weinkove2016, Szekelyhidi-Weinkove2021}, \cite{Chen-Xu2022,Colesanti-Focardi-Guan-Salani-arXiv2025, Guan-Li-Zhang2009, Guan-Ma2003, Guan-Ma-Zhou2006}, and the references therein.

Similar to the methods of concavity maximum principle, convex envelope and constant rank theorem, Makar-Limanov's approach has also been extended to high dimensions to deal with $1/2$-concavity of the torsion function and the log-concavity of the first Dirichlet eigenfunction of the Laplacian, see \cite{Jia-Ma-Shi2023-PA, Jia-Ma-Shi2023-SCM, Ma-Shi-Ye2012}. There are also many other miscellaneous interesting results on convexity of solutions of elliptic PDEs. For example, it is worth mentioning that Kulczycki \cite{Kulczycki2017} drew inspiration from the constant rank theorem and showed the concavity of solution to the Dirichlet problem of $1/2$-Laplacian in dimension two. In Hu \cite{Hu2024}, it was shown that the complex Hessian of the Dirichlet problem for the homogeneous complex Monge-Amp\`{e}re equation on a ring-shaped domain in $\mathbb{C}^n$ has constant rank $n-1$ in the weak sense. More recently, Hu and Sheng \cite{Hu-Sheng-arXiv2026} proved that the potential function of a complete K\"{a}hler-Einstein metric on a bounded strictly convex domain in $\mathbb{C}^n$ is strictly convex.

\subsection{Strategy of the proof}\label{Section 1.3}

Let us outline the fundamental framework for proving Theorem \ref{Theorem-Main}. We will follow the approach of Caffarelli and Friedman \cite{Caffarelli-Friedman1985}, and Korevaar and Lewis \cite{Korevaar-Lewis1987}, who established the power convexity of the torsion function and the first Dirichlet eigenfunction of the Laplacian by microscopic convexity principle. It is well-known that two key steps are required to implement this approach. The first is to obtain convexity estimates near the boundary, and the second is to establish the constant rank theorem. Recall the Dirichlet problem \eqref{Equ-CMA1} for the complex Monge-Amp\`{e}re equation that we considered in this paper. Letting $v=-\sqrt{-u/2}$, we know that $v$ satisfies the Dirichlet problem 
\begin{equation*}
\begin{cases}
F(\nabla^2 v, \nabla v, v)=0 & \mathrm{in} \ \Omega, \\
\hspace{22mm}v=0 & \mathrm{on} \ \partial\Omega,
\end{cases}
\end{equation*}
where
\begin{align*}
F(\nabla^2v, \nabla v, v)
=&~v^2\left[(v_{11}+v_{33})(v_{22}+v_{44})-(v_{12}+v_{34})(v_{21}+v_{43})-(v_{14}-v_{32})(v_{41}-v_{23})\right]\\
&+v\left[(v_2^2+v_4^2)(v_{11}+v_{33})+(v_1^2+v_3^2)(v_{22}+v_{44})\right.\\
&\hspace{8mm}\left.-(v_1v_2+v_3v_4)(v_{12}+v_{34}+v_{21}+v_{43})\right.\\
&\hspace{8mm}\left.-(v_1v_4-v_2v_3)(v_{14}-v_{32}+v_{41}-v_{23})\right]-1.
\end{align*}
For more details on the derivation of this equation, we refer to the next section (Section \ref{Section 2.1}). Since the domain $\Omega$ is strictly convex, we can easily obtain the convexity estimates for $v$ near the boundary $\partial\Omega$ (see Proposition \ref{Proposition-BCE} for precise clarification). It suffices to establish the constant rank theorem, namely that $\nabla^2v$ has constant rank in $\Omega$, for the transformed equation 
\begin{equation}\label{Equ-F-v}
F(\nabla^2 v, \nabla v, v)=0. 
\end{equation}

Let us explain why establishing the constant rank theorem for equation \eqref{Equ-F-v} is particularly challenging. Note that equation \eqref{Equ-F-v} is non-homogeneous with respect to the variable $\nabla^2v$. In general, it is difficult to verify the crucial inverse convexity structural condition \eqref{Inverse convexity structural condition} outlined in Bian and Guan \cite{Bian-Guan2009} for non-homogeneous equations. Although Bryan, Ivaki and Scheuer \cite{Bryan-Ivaki-Scheuer2023} presents a constant rank theorem for non-homogeneous equations, their conditions do not hold for $k$-Hessian type equations after power transformation. So we have to prove the constant rank theorem directly. This is why, up to now only convexity results in low dimensions \cite{Chen-Jia-Xiong2025,Huang2019,Liu-Ma-Xu2010,Ma-Xu2008,Zhang-Zhou2023} can be established for fully nonlinear elliptic PDEs.

In the complex setting, Hermitian matrix can be diagonalized by unitary transformations. Our goal, however, is to study the real convexity of solutions, for which the relevant object is the real Hessian $\nabla^2v$. A major difficulty is that the unitary transformations naturally associated with the complex structure generally do not diagonalize $\nabla^2v$, making direct computations difficult. To address this, we seek a coordinate system in which $\nabla^2v$ is diagonal while both the equation and the auxiliary quantities are invariant under the coordinate transformation. By exploiting the Autonne-Takagi factorization, we construct a suitable complex linear transformation satisfying these requirements. The details are presented in Section \ref{Section-auxiliary function}.

Define $$l=\min_{x\in\Omega}\mathrm{rank}(\nabla^2v(x)).$$ At a point $x_0\in\Omega$ where $\nabla^2v(x_0)$ attains the minimal rank $l$ ($0\leq l\leq 4$), we can choose suitable coordinates such that $\nabla^2v(x_0)$ is diagonal. If $l=0$ or $l=1$, this would contradict with equation \eqref{Equ-F-v}, since the left-hand side of \eqref{Equ-F-v} is negative while the right-hand side is $0$. If $l=4$, it is evident that $\nabla^2v(x)$ is of full rank for any $x\in\Omega$. Therefore, we only need to consider two cases, i.e.,

$\bullet$ $l=2$ (Theorem \ref{Theorem-rank2});

$\bullet$ $l=3$ (Theorem \ref{Theorem-rank3}).

To establish the constant rank theorem, the key step is to choose an appropriate auxiliary function, derive suitable differential inequalities, and apply the strong maximum principle.

In the celebrated work of Bian and Guan \cite{Bian-Guan2009} on fully nonlinear elliptic equations $$F(\nabla^2v,\nabla v,v,x)=0\quad\mathrm{in}\ \Omega\subset\mathbb{R}^n,$$ they introduced the auxiliary function $$\phi(x)=\sigma_{l+1}(\nabla^2v(x))+q(\nabla^2v(x)),$$ where $$q(\nabla^2v(x))=\begin{cases}
\frac{\sigma_{l+2}(\nabla^2v(x))}{\sigma_{l+1}(\nabla^2v(x))} & \mathrm{if}\ \sigma_{l+1}(\nabla^2v(x))>0,\\
0 & \mathrm{if}\ \sigma_{l+1}(\nabla^2v(x))=0.
\end{cases}$$ See Section \ref{Section 3.1} for more discussions on the elementary symmetric functions. For $l=2$, this function depends on $\sigma_3(\nabla^2v)$, which is not invariant under the admissible complex linear transformations (i.e., those constructed to diagonalize $\nabla^2v$ locally while preserving the form of the equation) because eigenvalues of $\nabla^2v$ are not preserved. To overcome this difficulty, we introduce a $2\times 2$ complex matrix $K$ defined in \eqref{Auxiliary matrix K}, satisfying $$\mathrm{rank}(\nabla^2v)=\mathrm{rank}(K)+2,$$ as stated in Lemma \ref{Lemma-rank identity}. This identity shows that the constant rank problem for $\nabla^2v$ is reduced to that for $K$. Under complex linear coordinate changes, $K$ transforms by similarity. Hence its eigenvalues are preserved, and consequently $\mathrm{tr}(K)$ and $\det(K)$ are invariant. This allows us to construct invariant auxiliary functions from $\mathrm{tr}(K)$ and $\det(K)$, which play a central role in our proof. A detailed discussion is given in Section \ref{Section-auxiliary function}.

For the case $l=2$, which corresponds to $\min_{x\in\Omega}\mathrm{rank}(K(x))=0$, we define the auxiliary function on $K$ in the sense of Bian and Guan as $$\phi(x)=\mathrm{tr}(K(x))+q(K(x)),$$ where $$q(K(x))=\begin{cases}
\frac{\det (K(x))}{\mathrm{tr}(K(x))} & \mathrm{if}\ \mathrm{tr}(K(x))>0,\\
0 & \mathrm{if}\ \mathrm{tr}(K(x))=0.
\end{cases}$$ As pointed out by Bian and Guan \cite{Bian-Guan2009}, using only the standard auxiliary function, in our framework simply $\mathrm{tr}(K)$, to prove the maximum principle is difficult for general fully nonlinear elliptic equations. Choosing $\mathrm{tr}(K)$ as the auxiliary function would lead to computational complexity that is unacceptable, even for the specific equation \eqref{Equ-F-v} considered here. Thus, the inclusion of the quotient term $q(K)$ is necessary.

The case $l=3$, corresponding to $\min_{x\in\Omega}\mathrm{rank}(K(x))=1$, is significantly challenging. To prove the differential inequality for the auxiliary function $\det(K)$ in this case, we must verify the semi-positivity of a six-variable quadratic form with most of its coefficients being non-zero. This part is highly involved and technical. Our strategy is performing elementary column and row operations cleverly on the corresponding $6\times6$ symmetric matrix. Fortunately, we can diagonalize this $6\times6$ matrix into two $3\times3$ blocks. By carefully computing and arranging the leading principal minors, and utilizing a basic lemma (Lemma \ref{Lemma-linear algebra}) from linear algebra, we finally complete the proof of Theorem \ref{Theorem-rank3}.

The rest of this paper is organized as follows. In Section \ref{Section 2}, we present the necessary preliminaries for the subsequent discussions, including the real and complex coordinate settings, the notations adopted throughout the paper, the choice of suitable coordinates, the definition of the auxiliary function, and so on. In Sections \ref{Section 3} and \ref{Section 4}, we establish the differential inequalities corresponding to the cases where the minimal rank is $2$ and $3$, respectively, which leads to the required constant rank theorem (Corollary \ref{Corollary-CRT}). This constitutes the main part of the paper. Finally, Section \ref{Section 5} is devoted to the proof of the main result in this paper, Theorem \ref{Theorem-Main}.

\section{Notation and Preliminaries}\label{Section 2}

In this section, the necessary preliminaries and notation for later use are provided. We first derive the real form of the two-dimensional complex Monge-Amp\`{e}re equation under a power transformation. We then discuss the choice of adapted coordinates and the construction of affine invariant auxiliary functions. Finally, we recall a linear algebra lemma that will be used repeatedly.

\subsection{Real form of the complex Monge-Amp\`{e}re equation}\label{Section 2.1}

This subsection provides a brief discussion on the real form of the equation for $v$, obtained from the complex Monge-Amp\`{e}re equation via a power transformation, on a domain $\Omega\subset\mathbb{C}^2$.

We begin with the complex Monge-Amp\`{e}re equation
\begin{equation}\label{Equ-complexMA}
\det(u_{i\bar{j}})=1.
\end{equation}
To avoid confusion with the real Hessian notation used later, we shall write all derivatives in terms of the complex coordinates $z_i$ and $\bar{z}_i$. In particular, we set $$u_i=u_{z_i}=\frac{\partial u}{\partial z_i},\quad u_{\bar{i}}=u_{\bar{z}_i}=\frac{\partial u}{\partial\bar{z}^i}.$$ With $v=-\sqrt{-
u/2}$, a direct computation gives $$u_{z_i}=-4vv_{z_i},\quad u_{z_i\bar{z}_j}=-4(vv_{z_i\bar{z}_j}+v_{z_i}v_{\bar{z}_j}),$$ so that the equation \eqref{Equ-complexMA} becomes
\begin{equation}\label{Equ-complexMA-v}
\det(vv_{z_i\bar{z}_j}+v_{z_i}v_{\bar{z}_j})=\frac{1}{16}.
\end{equation}

To work with real coordinates, we express the complex derivatives in terms of $z_i=x_i+\sqrt{-1}y_i$. For a real-valued function $v$ defined on a domain $\Omega\subset\mathbb{C}^2$, then the complex derivatives $v_{z_i}$ and $v_{\bar{z}_i}$ are given by
\begin{align*}
&v_{z_i}=\frac{\partial v}{\partial z_i}=\frac{1}{2}\left(\frac{\partial v}{\partial x_i}-\sqrt{-1}\frac{\partial v}{\partial y_i}\right)=\frac{1}{2}(v_{x_i}-\sqrt{-1}v_{y_i})
\intertext{and}
&v_{\bar{z}_i}=\frac{\partial v}{\partial\bar{z}_i}=\frac{1}{2}\left( \frac{\partial v}{\partial x_i}+\sqrt{-1}\frac{\partial v}{\partial y_i}\right)=\frac{1}{2}(v_{x_i}+\sqrt{-1}v_{y_i}).
\end{align*}
Furthermore, we have
\begin{equation}\label{Equality-v_{ibarj}}
v_{z_i\bar{z}_j}=\frac{\partial^2v}{\partial z_i\partial\bar{z}_j}=\frac{1}{4}(v_{x_ix_j}+v_{y_iy_j})+\frac{1}{4}\sqrt{-1}(v_{x_iy_j}-v_{y_ix_j})
\end{equation}
and $$v_{z_iz_j}=\frac{\partial^2v}{\partial z_i\partial z_j}=\frac{1}{4}(v_{x_ix_j}-v_{y_iy_j})-\frac{1}{4}\sqrt{-1}(v_{x_iy_j}+v_{y_ix_j}).$$ 

Let $A=(v_{z_i\bar{z}_j})$, $B=(v_{z_iz_j})$, and write the real Hessian of $v$ in block form as $$\nabla^2v=\begin{pmatrix}
U & V\\
V^T & W\end{pmatrix},$$ where $U=(v_{x_ix_j})$, $V=(v_{x_iy_j})$, and $W=(v_{y_iy_j})$. Then $$A=\frac{1}{4}(U+W)+\frac{\sqrt{-1}}{4}(V-V^T),$$ and $$B=\frac{1}{4}(U-W)-\frac{\sqrt{-1}}{4}(V+V^T).$$ Since $u$ is strictly plurisubharmonic and $v<0$, it follows from $u_{z_i\bar{z}_j}=-4(vv_{z_i\bar{z}_j}+v_{z_i}v_{\bar{z}_j})$ that $A=(v_{z_i\bar{z}_j})$ is positive definite. In fact, for any nonzero vector $\xi\in\mathbb{C}^2$, $$\sum_{i,j}v_{z_i\bar{z}_j}\xi_i\bar{\xi}_j=-\frac{1}{4v}\sum_{i,j}u_{z_i\bar{z}_j}\xi_i\bar{\xi}_j-\frac{1}{v}\left|\sum_iv_{z_i}\xi_i\right|^2>0.$$ In particular, $A$ is invertible.

For the sake of convenience, we adopt the following notation for real partial derivatives, which will be used throughout the rest of the paper
\begin{align*}
&v_{x_1}=v_1, \ v_{x_2}=v_2, \ v_{y_1}=v_3, \ v_{y_2}=v_4, \\
&v_{x_1x_1}=v_{11}, \ v_{y_1y_1}=v_{33}, \ v_{x_1y_1}=v_{13}, \ \ldots
\end{align*}
Substituting \eqref{Equality-v_{ibarj}} into \eqref{Equ-complexMA-v}, we obtain the following real equation
\begin{equation}\label{Equ-real equation(v)}
\begin{aligned}
1=&~v^2\left[(v_{11}+v_{33})(v_{22}+v_{44})-(v_{12}+v_{34})(v_{21}+v_{43})-(v_{14}-v_{32})(v_{41}-v_{23})\right]\\
&+v\left[(v_2^2+v_4^2)(v_{11}+v_{33})+(v_1^2+v_3^2)(v_{22}+v_{44})\right.\\
&\hspace{8mm}\left.-(v_1v_2+v_3v_4)(v_{12}+v_{34}+v_{21}+v_{43})\right.\\
&\hspace{8mm}\left.-(v_1v_4-v_2v_3)(v_{14}-v_{32}+v_{41}-v_{23})\right].
\end{aligned}
\end{equation}
Hence $v$ is the unique solution of the Dirichlet problem
\begin{equation}\label{Equ-Dirichlet problem-real}
\begin{cases}
F(\nabla^2 v, \nabla v, v)=0 & \mathrm{in} \ \Omega, \\
\hspace{22mm}v=0 & \mathrm{on} \ \partial\Omega,
\end{cases}
\end{equation}
where
\begin{align*}
F(\nabla^2v, \nabla v, v)
=&~v^2\left[(v_{11}+v_{33})(v_{22}+v_{44})-(v_{12}+v_{34})(v_{21}+v_{43})-(v_{14}-v_{32})(v_{41}-v_{23})\right]\\
 &+v\left[(v_2^2+v_4^2)(v_{11}+v_{33})+(v_1^2+v_3^2)(v_{22}+v_{44})\right.\\
 &\hspace{8mm}\left.-(v_1v_2+v_3v_4)(v_{12}+v_{34}+v_{21}+v_{43})\right.\\
 &\hspace{8mm}\left.-(v_1v_4-v_2v_3)(v_{14}-v_{32}+v_{41}-v_{23})\right]-1.
\end{align*}

For simplicity, we denote
\begin{alignat*}{3}
F^{ij}&=\frac{\partial F}{\partial v_{ij}}, \quad &F^{v_k}&=\frac{\partial F}{\partial v_k}, \quad &F^v&=\frac{\partial F}{\partial v}, \\
F^{ij, kl}&=\frac{\partial^2F}{\partial v_{ij}\partial v_{kl}}, \quad &F^{ij, v_k}&=\frac{\partial^2F}{\partial v_{ij}\partial v_k}, \quad &F^{ij, v}&=\frac{\partial^2F}{\partial v_{ij}\partial v}, \\
F^{v_k, v_l}&=\frac{\partial^2F}{\partial v_k\partial v_l}, \quad &F^{v_k, v}&=\frac{\partial^2F}{\partial v_k\partial v}, \quad &F^{v, v}&=\frac{\partial^2F}{\partial v^2}, \hspace{0.4cm}\mathrm{etc}.
\end{alignat*}
These notations will be used throughout Sections \ref{Section 3} and \ref{Section 4}.

\subsection{Choice of an auxiliary function}\label{Section-auxiliary function}

The primary objective of this subsection is to construct appropriate affine invariant auxiliary functions to be used later in the proof of the constant rank theorems for the real Hessian $\nabla^2v$. The construction relies on two fundamental properties, namely the existence of an adapted coordinate system in which $\nabla^2v$ can be locally diagonalized, and the fact that the transformed complex Monge-Amp\`{e}re equation is invariant under such a coordinate choice.

We define the auxiliary matrix $K$ by 
\begin{equation}\label{Auxiliary matrix K}
K=I-B\overline{A^{-1}}\overline{B}A^{-1},
\end{equation}
where $I$ denotes the identity matrix of size $2\times 2$, $A=(v_{z_i\bar{z}_j})_{1\leq i,j\leq 2}$ and $B=(v_{z_iz_j})_{1\leq i,j\leq 2}$. The rank relation between the real Hessian $\nabla^2v$ and the auxiliary matrix $K$ is now considered. 
\begin{lemma}\label{Lemma-rank identity}
The real Hessian $\nabla^2v$ and the auxiliary matrix $K$ satisfy
\begin{equation}\label{rank D^2v=rank K+2}
\mathrm{rank}(\nabla^2v)=\mathrm{rank}(K)+2.
\end{equation}
\end{lemma}
\begin{proof}
The standard real-complex correspondence gives the following matrix identity
\begin{equation}\label{Matrix identity-K and D^2v}
4\begin{pmatrix}
A & B\\[6pt]
\overline{B} & \overline{A}
\end{pmatrix}=\begin{pmatrix}
I & -\sqrt{-1}I\\[6pt]
I & \sqrt{-1}I
\end{pmatrix}\begin{pmatrix}
U & V\\[6pt]
V^T & W
\end{pmatrix}\begin{pmatrix}
I & I\\[6pt]
\sqrt{-1}I & -\sqrt{-1}I
\end{pmatrix}=T\nabla^2vT^H,
\end{equation}
where $T^H$ is the conjugate transpose of the complex matrix $T$. Since the transformation matrix $T$ is invertible, it preserves the matrix rank, yielding $$\mathrm{rank}(\nabla^2v)=\mathrm{rank}\begin{pmatrix}
 A & B\\[6pt]
 \overline{B} & \overline{A}
\end{pmatrix}.$$ On the other hand, a standard block decomposition yields $$\begin{pmatrix}
I & -B\overline{A^{-1}}\\[6pt]
0 & I
\end{pmatrix}\begin{pmatrix}
A & B\\[6pt]
\overline{B} & \overline{A}
\end{pmatrix}\begin{pmatrix}
I & 0\\[6pt]
-\overline{A^{-1}}\overline{B} & I
\end{pmatrix}=\begin{pmatrix}
A-B\overline{A^{-1}}\overline{B} & 0\\[6pt]
0 & \overline{A}
\end{pmatrix},$$ which implies the rank identity $$\mathrm{rank}\begin{pmatrix}
A & B\\[6pt]
\overline{B} & \overline{A}
\end{pmatrix}=\mathrm{rank}(A-B\overline{A^{-1}}\overline{B})+\mathrm{rank}(A).$$ Since $A$ is an invertible $2\times $2 matrix and $K$ is defined by \eqref{Auxiliary matrix K}, we obtain $$\mathrm{rank}(\nabla^2v)=\mathrm{rank}(K)+2.$$
\end{proof}

To simplify local computations, it is necessary to choose coordinates in which the real Hessian $\nabla^2v$ is diagonal at a given point $x_0\in\Omega$. According to Lemma A.4 in Hu \cite{Hu2025}, which is based on the Autonne-Takagi factorization (see also Lemma A.3 therein), there exists an invertible matrix $P$ such that $$PAP^H=I,\quad PBP^T=\Lambda,$$ where $\Lambda$ is a nonnegative diagonal matrix. We now recall the explicit expressions of the complex Hessian blocks $A$ and $B$ in terms of the real second derivatives of $v$, namely 
\begin{equation}\label{Matrix-A}
A=\begin{pmatrix}
\frac{1}{4}(v_{11}+v_{33}) & \frac{1}{4}(v_{12}+v_{34})+\frac{\sqrt{-1}}{4}(v_{14}-v_{32})\\[9pt]
\frac{1}{4}(v_{21}+v_{43})+\frac{\sqrt{-1}}{4}(v_{23}-v_{41}) & \frac{1}{4}(v_{22}+v_{44})
\end{pmatrix}
\end{equation}
and
\begin{equation}
B=\begin{pmatrix}\label{Matrix-B}
\frac{1}{4}(v_{11}-v_{33})-\frac{\sqrt{-1}}{4}(v_{13}+v_{31}) & \frac{1}{4}(v_{12}-v_{34})-\frac{\sqrt{-1}}{4}(v_{14}+v_{32})\\[9pt]
\frac{1}{4}(v_{21}-v_{43})-\frac{\sqrt{-1}}{4}(v_{23}+v_{41}) & \frac{1}{4}(v_{22}-v_{44})-\frac{\sqrt{-1}}{4}(v_{24}+v_{42})
\end{pmatrix}.
\end{equation}
The diagonalizability of $A$ requires $$v_{12}+v_{34}=0,\quad v_{14}-v_{32}=0,$$ while the diagonalizability of $B$ yields $$v_{12}-v_{34}=0,\quad v_{14}+v_{32}=0.$$ Combining these relations leads to $$v_{12}=v_{34}=v_{14}=v_{23}=0.$$ Moreover, since $PBP^T=\Lambda$ is a real diagonal matrix, we further obtain $$v_{13}=v_{24}=0.$$ This confirms the existence of an adapted coordinate system at $x_0$ in which the real Hessian $\nabla^2v$ is diagonal.

We then normalize the coordinate transformation such that $|\det(P)|=1$, which preserves the affine invariance of equation \eqref{Equ-complexMA-v}. Under this normalization, $$PAP^H=\sqrt{\det A}I,\quad PBP^T=\sqrt{\det A}\Lambda,$$ and $$\det(P(vv_{z_i\bar{z}_j}+v_{z_i}v_{\bar{z}_j})P^H)=\det(vv_{z_i\bar{z}_j}+v_{z_i}v_{\bar{z}_j})=\frac{1}{16},$$ which shows that the equation is invariant under complex linear transformations with $|\det(P)|=1$.

Applying the above transformation, the matrices $A$ and $B$ transform as $$A\mapsto PAP^H,\quad B\mapsto PBP^T.$$ Substituting these into the definition of $K$, we obtain
\begin{align*}
\widetilde{K}
=&~I-(PBP^T)\overline{(PAP^H)^{-1}}\overline{(PBP^T)}(PAP^H)^{-1}\\
=&~I-PBP^T\overline{(P^H)^{-1}A^{-1}P^{-1}}\overline{PBP^T}(P^H)^{-1}A^{-1}P^{-1}\\
=&~PP^{-1}-PB\overline{A^{-1}}\overline{B}A^{-1}P^{-1}\\
=&~P(I-B\overline{A^{-1}}\overline{B}A^{-1})P^{-1}\\
=&~PKP^{-1}.
\end{align*}
Hence $K$ is transformed via a similarity transformation. Consequently, its eigenvalues are preserved, and therefore both $\mathrm{tr}(K)$ and $\det(K)$ are invariant under this transformation.

As discussed in Section \ref{Section 1.3}, it suffices to prove the constant rank theorem for $\nabla^2v$ when its minimal rank is $2$ or $3$. In view of \eqref{rank D^2v=rank K+2}, the constant rank problem for $\nabla^2v$ is reduced to that for the auxiliary matrix $K$, namely proving the constant rank theorem in the cases $\min_{x\in\Omega}\mathrm{rank}(K(x))=0$ and $\min_{x\in\Omega}\mathrm{rank}(K(x))=1$, respectively. We therefore select suitable auxiliary functions for there two cases.

$\bullet$ $\min_{x\in\Omega}\mathrm{rank}(K(x))=0$. In this case, we employ the quotient type auxiliary function introduced by Bian and Guan \cite{Bian-Guan2009}, 
\begin{equation}\label{Auxiliary function-rankK=0}
\phi(x)=\mathrm{tr}(K(x))+q(K(x)), 
\end{equation}
where $$q(K(x))=\begin{cases}
\frac{\det(K(x))}{\mathrm{tr}(K(x))} & \mathrm{if}\ \mathrm{tr}(K(x))>0,\\
0 & \mathrm{if}\ \mathrm{tr}(K(x))=0.
\end{cases}$$ The affine invariance of $\mathrm{tr}(K(x))$ and $\det(K(x))$ ensures that $\phi$ is also invariant. Thus, $\phi$ may be computed in the above diagonalized coordinate system without loss of generality.

$\bullet$ $\min_{x\in\Omega}\mathrm{rank}(K(x))=1$. We take the standard auxiliary function
\begin{equation}\label{Auxiliary function-rankK=1}
\phi=\det(K(x)),
\end{equation}
as in Caffarelli and Friedman \cite{Caffarelli-Friedman1985}, Korevaar and Lewis \cite{Korevaar-Lewis1987}, and other related works. In fact, this is consistent with the construction of Bian and Guan \cite{Bian-Guan2009}.

We shall apply these two auxiliary functions in Section \ref{Section 3} and \ref{Section 4} to obtain the differential inequalities stated in Theorem \ref{Theorem-rank2} and \ref{Theorem-rank3}, which allow the strong maximum principle to establish the constant rank theorem.

\subsection{A linear algebra lemma}

We end this section with a useful result from linear algebra. Although it is standard and can be found in many textbooks, it will play a key role in the proof of the constant rank theorem for both the minimal rank $2$ and minimal rank $3$ cases. For completeness, we include a proof here.

\begin{lemma}\label{Lemma-linear algebra}
Let $M$ be an $n\times n$ real symmetric matrix. Assume that $M$ has $n-1$ nonzero leading principal minors, denoted by $P_1$, $P_2$, $\cdots$, $P_{n-1}$. Then there exists an $n\times n$ invertible matrix $C$ such that  $$C^TMC=\mathrm{diag}\left(P_1, \frac{P_2}{P_1}, \cdots, \frac{P_n}{P_{n-1}}\right), $$ where $P_n=\det(M)$.
\end{lemma}
\begin{proof}
We proceed by induction on $n$. The result is evident when $n=1$. Assume that the statement holds for $n-1$, and let $$M=\begin{pmatrix}
M_{n-1} & \alpha\\[6pt]
\alpha^T & a_{nn}
\end{pmatrix}, $$ where $M_{n-1}$ is the $(n-1)\times (n-1)$ leading principal submatrix of $M$. Since $\det(M_{n-1})=P_{n-1}\neq 0$, we can define an $n\times n$ invertible matrix $$C_1=\begin{pmatrix}
I_{n-1} & -M_{n-1}^{-1}\alpha\\[6pt]
O & 1
\end{pmatrix}. $$ Then we have $$C_1^TMC_1=\begin{pmatrix}
M_{n-1} & O\\[6pt]
O & a_{nn}-\alpha^TM_{n-1}^{-1}\alpha
\end{pmatrix}. $$
Note that $\det(C_1)=1$. It yields that $$\det(M)=(a_{nn}-\alpha^T M_{n-1}^{-1}\alpha) \det(M_{n-1}), $$ that is, $$a_{nn}-\alpha^TM_{n-1}^{-1} \alpha=\frac{P_n}{P_{n-1}}. $$ By the inductive hypothesis, there exists an $(n-1)\times (n-1)$ invertible matrix $C_2$ such that
$$C_2^TM_{n-1}C_2=\mathrm{diag}\left(P_1, \frac{P_2}{P_1}, \cdots, \frac{P_{n-1}}{P_{n-2}}\right). $$ If we let $C_3=\begin{pmatrix}
C_2 & O\\
O & 1
\end{pmatrix}$ and $C=C_1C_3$, then we have $$C^TMC=\mathrm{diag}\left(P_1, \frac{P_2}{P_1}, \cdots, \frac{P_n}{P_{n-1}}\right). $$
\end{proof}

Up to this point, we have completed the setup of notation, the choice of coordinates, and the construction of the key auxiliary functions used in the sequel.

\section{Constant rank theorem: the case of minimal rank $2$}\label{Section 3}

This section is devoted to the constant rank theorem for the real Hessian $\nabla^2v$ in the case where it attains its minimal rank $2$, equivalently, where the auxiliary matrix $K$ has rank $0$. More precisely, our proof is inspired by Bian and Guan \cite{Bian-Guan2009}, where we employ the quotient type auxiliary function defined in \eqref{Auxiliary function-rankK=0}. The first step is to perform a careful computation of the first and second order derivatives of the quotient component of the test function. Based on these expansions, we then establish a strong maximum principle for $\phi$, which yields the constant rank theorem in this setting.

Throughout this and next section, we adopt the convention that all summation indices $i$, $j$, $k$, $l$, $p$, $q$, $r$, $s$ and $\alpha$, $\beta$, $\gamma$, $\delta$ run from $1$ to $4$, unless stated otherwise.

\subsection{Estimates for the auxiliary function}\label{Section 3.1}

In this subsection, we focus on the quotient term in the auxiliary function \eqref{Auxiliary function-rankK=0}, whose calculation is rather delicate.

We recall the definition of $k$-th elementary symmetric functions. For $1\leq k\leq n$ and $\lambda=(\lambda_1, \lambda_2, \cdots, \lambda_n)\in\mathbb{R}^n$, we define $$\sigma_k(\lambda)=\sum_{1\leq i_1<i_2<\cdots<i_k\leq n}\lambda_{i_1}\lambda_{i_2}\cdots\lambda_{i_k}. $$ We also make the convention that $\sigma_0(\lambda)=1$ and $\sigma_{-1}(\lambda)=\sigma_{n+1}(\lambda)=0$. The definition of $\sigma_k$ for vectors in $\mathbb{R}^n$ can be extended to $n\times n$ symmetric matrices in a natural way by letting $\sigma_k(W)=\sigma_k(\lambda(W))$, where $\lambda(W)=(\lambda_1(W), \lambda_2(W), \cdots, \lambda_n(W))$ are the eigenvalues of the symmetric matrix $W$. We denote by $\sigma_{k-1}(W|i)$ the $(k-1)$-th elementary symmetric function with $\lambda_i=0$ and $\sigma_{k-2}(W|ij)$ $(i\neq j)$ the $(k-2)$-th elementary symmetric function with $\lambda_i=\lambda_j=0$.
	
The first and second derivatives of the elementary symmetric functions are listed in the following lemma.
\begin{lemma}\label{Lemma-sigma_k}
Suppose that $W=(W_{ij})_{n\times n}$ is diagonal. For $k=1, 2, \cdots, n$, then we have $$\frac{\partial\sigma_k(W)}{\partial W_{ij}}=
\begin{cases}
\sigma_{k-1}(W|i) & \mathrm{if}\ i=j\\
0 & \mathrm{if}\ i\neq j,
\end{cases}$$ and
$$\frac{\partial^2\sigma_k(W)}{\partial W_{ij}\partial W_{pq}}=
\begin{cases}
\sigma_{k-2}(W|ip) & \mathrm{if}\ i=j, p=q, i\neq p, \\
-\sigma_{k-2}(W|ip) & \mathrm{if}\ i=q, j=p, i\neq p, \\
0, & \mathrm{otherwise}.
\end{cases}$$
\end{lemma}

We briefly recall the result of Bian and Guan \cite{Bian-Guan2009}, which has been repeatedly referred to in the preceding discussion. They introduced the quotient type auxiliary function
\begin{equation}\label{Auxiliary function-Bian and Guan-D^2v}
\phi(x)=\sigma_{l+1}(\nabla^2v(x))+q(\nabla^2v(x)),
\end{equation}
where $$q(\nabla^2v(x))=\begin{cases}
\frac{\sigma_{l+2}(\nabla^2v(x))}{\sigma_{l+1}(\nabla^2v(x))} & \mathrm{if}\ \sigma_{l+1}(\nabla^2v(x))>0,\\
0 & \mathrm{if}\ \sigma_{l+1}(\nabla^2v(x))=0.
\end{cases}$$
Using the auxiliary function \eqref{Auxiliary function-Bian and Guan-D^2v}, they established a constant rank theorem for convex solutions of the general fully nonlinear elliptic equation $$F(\nabla^2v,\nabla v,v,x)=0\quad \mathrm{in}\ \Omega\subset\mathbb{R}^n,$$ in the case where $\nabla^2v$ attains its minimal rank $l$, under the structural condition that 
\begin{equation}\label{Inverse convexity structural condition}
F(M^{-1},p,v,x)\ \mbox{is locally convex in $(M,v,x)$ for each $p$}.
\end{equation}

Below we review the notations from Bian and Guan \cite{Bian-Guan2009} and adapt it to the present setting. Suppose that $\nabla^2v(x)$ attains its minimal rank $2$ at $x_0\in\Omega$. We pick an open neighborhood $\mathcal{O}$ of $x_0$. For any $x\in\mathcal{O}$, let $\lambda_1(x)\geq\lambda_2(x)\geq\lambda_3(x)\geq\lambda_4(x)\geq 0$ be the eigenvalues of $\nabla^2v(x)$. There are exists a positive constant $c$ independent of $x$, such that $\lambda_1(x)\geq\lambda_2(x)\geq c$. Let $G=\{1,2\}$ and $B=\{3,4\}$ be the ``good" and ``bad" sets of indices, respectively. If there is no confusion, we also denote the ``good" and ``bad" eigenvalues of $\nabla^2v(x)$ by $G$ and $B$, respectively. Note that for any $\delta>0$, we may choose $\mathcal{O}$ small enough such that $\lambda_i(x)<\delta$ for all $i\in B$ and $x\in\mathcal{O}$. For two functions $h(x)$ and $k(x)$ defined on $\mathcal{O}$, by $h(x)=O(k(x))$ we mean that there exists a positive constant $C$ such that $|h(x)|\leq C|k(x)|$.

As mentioned in Section \ref{Section-auxiliary function}, see \eqref{rank D^2v=rank K+2}, we established the rank identity $$\mathrm{rank}(\nabla^2v)=\mathrm{rank}(K)+2.$$ Therefore, the case where $\nabla^2v$ has minimal rank $2$ is equivalent to the auxiliary matrix $K$ having minimal rank $0$. Section \ref{Section 1.3} and Section \ref{Section-auxiliary function} indicate that we now apply the framework of Bian and Guan to the auxiliary matrix $K$ instead of the real Hessian $\nabla^2v$. Accordingly, we consider the auxiliary function
\begin{equation}\label{Auxiliary function-Bian and Guan-K}
\phi(x)=\mathrm{tr}(K(x))+q(K(x)),
\end{equation}
where $$q(K(x))=\begin{cases}
\frac{\det(K(x))}{\mathrm{tr}(K(x))} & \mathrm{if}\ \mathrm{tr}(K(x))>0,\\
0 & \mathrm{if}\ \mathrm{tr}(K(x))=0.
\end{cases}$$

Next, we express $\phi(x)$ in terms of $v_{ij}$ for $1\leq i,j\leq 4$. According to the definition of the matrix $K$ in \eqref{Auxiliary matrix K} and \eqref{Matrix identity-K and D^2v}, we have
\begin{equation}\label{Equality-det K and det D^2v}
\det(K)=(\det(A))^{-2}\det(\nabla^2v),
\end{equation}
where $$\det(A)=\frac{1}{16}\big((v_{11}+v_{33})(v_{22}+v_{44})-(v_{12}+v_{34})^2-(v_{14}-v_{32})^2\big).$$ Rewriting $$K=(\det(A))^{-2}\big((\det(A))^2I-B\overline{A^\ast}\overline{B}A^\ast\big),$$ where $A^\ast$ denotes the adjugate matrix of $A$, it follows that $\mathrm{tr}(K)=(\det (A))^{-2}\mathrm{tr}\big((\det(A))^2I-B\overline{A^\ast}\overline{B}A^\ast\big)$. Denoting $$\psi=\mathrm{tr}\big((\det(A))^2I-B\overline{A^\ast}\overline{B}A^\ast\big),$$ it follows from \eqref{Equality-det K and det D^2v} that the auxiliary function in \eqref{Auxiliary function-Bian and Guan-K} can be rewritten as 
\begin{equation}\label{Auxiliary function-Bian and Guan-K-1}
\phi(x)=(\det(A))^{-2}\psi(x)+\frac{\sigma_4(\nabla^2v(x))}{\psi(x)}.
\end{equation}
For notational convenience, we shall still denote
\begin{equation}\label{q(D^2v(x))-new}
q(\nabla^2v(x))=\frac{\sigma_4(\nabla^2v(x))}{\psi(x)}.
\end{equation}

To compute $\phi(x)$, it remains to determine $\psi(x)$. From the explicit expressions of the matrices $A$ and $B$ in \eqref{Matrix-A} and \eqref{Matrix-B}, a straightforward calculation gives 
\begin{align*}
&~B\overline{A^\ast}\\
=&~\frac{1}{16}\scalebox{0.85}{$\begin{pmatrix}
\makecell{(v_{11}-v_{33})(v_{22}+v_{44})+(v_{23}^2+v_{34}^2-v_{12}^2-v_{14}^2)\\+2\sqrt{-1}(v_{12}v_{23}+v_{14}v_{34}-v_{13}(v_{22}+v_{44}))} & \makecell{2(v_{12}v_{33}-v_{11}v_{34}+v_{13}(v_{14}-v_{32}))\\+2\sqrt{-1}(v_{13}(v_{12}+v_{34})-v_{11}v_{32}-v_{33}v_{14})}\\[20pt]
\makecell{2(v_{21}v_{44}-v_{22}v_{43}+v_{24}(v_{23}-v_{41}))\\+2\sqrt{-1}(v_{24}(v_{21}+v_{43})-v_{23}v_{44}-v_{22}v_{14})} & \makecell{(v_{22}-v_{44})(v_{11}+v_{33})+(v_{14}^2+v_{34}^2-v_{12}^2-v_{23}^2)\\+2\sqrt{-1}(v_{12}v_{14}+v_{23}v_{34}-v_{24}(v_{11}+v_{33}))}
\end{pmatrix}$}.
\end{align*}
Since $\overline{B}A^\ast=\overline{B\overline{A^\ast}}$, we have $$\mathrm{tr}(B\overline{A^\ast}\overline{B}A^\ast)=\frac{1}{16^2}\big(|a|^2+|d|^2+2\Re(b\bar{c})\big),$$ where $a,b,c,d$ represents the $(1,1)$, $(1,2)$, $(2,1)$, and $(2,2)$ entries of the matrix $B\overline{A^\ast}$, respectively. Therefore, direct computation yields
\begin{equation}\label{Equality-psi}
\begin{aligned}
\psi
=&~\frac{2}{16^2}\big((v_{11}+v_{33})(v_{22}+v_{44})-(v_{12}+v_{34})^2-(v_{14}-v_{32})^2\big)^2\\
&-\frac{1}{16^2}\bigg[\big((v_{11}-v_{33})(v_{22}+v_{44})+(v_{23}^2+v_{34}^2-v_{12}^2-v_{14}^2)\big)^2\\
&\hspace{1.2cm}+\big((v_{22}-v_{44})(v_{11}+v_{33})+(v_{14}^2+v_{34}^2-v_{12}^2-v_{23}^2)\big)^2\\[6pt]
&\hspace{1.2cm}+4\big(v_{12}v_{23}+v_{14}v_{34}-v_{13}(v_{22}+v_{44})\big)^2+4\big(v_{12}v_{14}+v_{23}v_{34}-v_{24}(v_{11}+v_{33})\big)^2\\[6pt]
&\hspace{1.2cm}+8\big(v_{12}v_{33}-v_{11}v_{34}+v_{13}(v_{14}-v_{32})\big)\big(v_{21}v_{44}-v_{22}v_{43}+v_{24}(v_{23}-v_{41})\big)\\
&\hspace{1.2cm}+8\big(v_{13}(v_{12}+v_{34})-v_{11}v_{32}-v_{33}v_{14}\big)\big(v_{24}(v_{21}+v_{43})-v_{23}v_{44}-v_{22}v_{14}\big)\bigg].
\end{aligned}
\end{equation}

To avoid $\psi(x)=0$ in \eqref{q(D^2v(x))-new}, we follow the perturbation argument of Bian and Guan \cite{Bian-Guan2009}. For a sufficiently small $\varepsilon>0$, define $v_\varepsilon(x)=v(x)+\frac{\varepsilon}{2}|x|^2$. Then $\nabla^2v_\varepsilon=\nabla^2v+\varepsilon I_4$, where $I_4$ is the $4\times 4$ identity matrix. We define the perturbed quantities $A_\varepsilon$ and $B_\varepsilon$ by replacing $\nabla^2v$ with $\nabla^2v_\varepsilon$ in the definitions \eqref{Matrix-A} and \eqref{Matrix-B}. Accordingly, we set $$q_\varepsilon(\nabla^2v_\varepsilon)=\frac{\sigma_4(\nabla^2v_\varepsilon)}{\psi_\varepsilon},\quad \phi_\varepsilon=(\det(A_\varepsilon))^{-2}\psi_\varepsilon+q_\varepsilon(\nabla^2v_\varepsilon).$$ Here $$\psi_\varepsilon=\mathrm{tr}\big((\det(A_\varepsilon))^2I-B_\varepsilon\overline{A_\varepsilon^\ast}\overline{B_\varepsilon}A_\varepsilon^\ast\big).$$ We will work with $\phi_\varepsilon$ in order to derive uniform $C^2$ estimates independent of $\varepsilon$. One may also work directly with $\phi$ at points where $\psi(x)\neq 0$, with all constants independent of the chosen point. For simplicity of notation, we will write $v$, $A$, $B$, $K$, $\psi$, $\phi$ in place of $v_\varepsilon$, $A_\varepsilon$, $B_\varepsilon$, $K_\varepsilon$, $\psi_\varepsilon$, $\phi_\varepsilon$, respectively. All estimates below are uniform in $\varepsilon$.

Without loss of generality, we work at an arbitrary fixed point $x\in\mathcal{O}$ in the coordinate system introduced in Section \ref{Section-auxiliary function}, in which $\nabla^2v$ is diagonal and satisfies
\begin{equation}\label{Equality-v_{11}+v_{33}=v_{22}+v_{44}}
v_{11}+v_{33}=v_{22}+v_{44}.
\end{equation}
From the definition of $\psi$ in \eqref{Equality-psi}, the diagonal form of $\nabla^2v$ yields
\begin{equation}\label{Equality-psi-diagonal}
\psi=\frac{4}{16^2}(v_{22}+v_{44})^2v_{11}v_{33}+\frac{4}{16^2}(v_{11}+v_{33})^2v_{22}v_{44}.
\end{equation}
Choosing $\mathcal{O}$ sufficiently small, we may assume that $c\leq v_{22}\leq v_{11}\leq C$ for some positive constants $c$ and $C$ independent of $x$. Therefore,
\begin{equation}\label{Estimate-psi}
\psi\geq\frac{4}{16^2}(v_{11}^2v_{22}v_{44}+v_{11}v_{22}^2v_{33})\geq c(v_{33}+v_{44}).
\end{equation}
Since the perturbation ensures $v_{33},v_{44}\geq\varepsilon$, we further obtain
\begin{equation}\label{Estimate-Perturbation}
\psi(x)\geq C\varepsilon,\quad \mbox{for all}\ x\in\mathcal{O},
\end{equation}
where $C$ is independent of $\varepsilon$. According to \eqref{Equality-psi-diagonal}, one has $$\phi\geq(\det(A))^{-2}\psi=4v_{11}v_{33}(v_{11}+v_{33})^{-2}+4v_{22}v_{44}(v_{22}+v_{44})^{-2}.$$ Since $v_{11}$ and $v_{22}$ are uniformly bounded from above and below, it follows that $$\phi\geq c(v_{33}+v_{44}).$$ On the other hand, combining $\psi\leq C(v_{33}+v_{44})$, with \eqref{Estimate-psi}, i.e., $\psi\geq c(v_{33}+v_{44})$, we have $$\phi
=(\det A)^{-2}\psi+\frac{\sigma_4(\nabla^2v)}{\psi}\leq C'(v_{33}+v_{44})+C''\frac{v_{33}v_{44}}{v_{33}+v_{44}}\leq C(v_{33}+v_{44}).$$ Consequently,
\begin{equation}\label{Estimate-v_{33}+v_{44}=phi}
0\leq c\phi\leq v_{33}+v_{44}\leq C\phi.
\end{equation}
In particular,
\begin{equation}\label{Estimate-v_{33},v_{44}}
 v_{33}=O(\phi),\quad v_{44}=O(\phi),
\end{equation}
and from \eqref{Equality-v_{11}+v_{33}=v_{22}+v_{44}}, we conclude
\begin{equation}\label{Estimate-v_{22}}
v_{22}=v_{11}+O(\phi).
\end{equation}

We are now ready to compute the first-order derivative of $q$ with respect to $x_\alpha$. Here and throughout this section, $\alpha=1,2,3,4$ correspond to the variables $x_1, x_2, y_1, y_2$, respectively.
\begin{lemma}\label{Lemma-first derivative of q}
Let $q(\nabla^2v)$ be defined as in \eqref{q(D^2v(x))-new}. Assume that $\nabla^2v$ is diagonal. Then
$$\frac{\partial q}{\partial x_\alpha}=O\bigg(\phi+\sum_{p,q\in B}|\nabla v_{pq}|\bigg).$$
\end{lemma}
\begin{proof}
A direct computation yields $$\frac{\partial q}{\partial x_\alpha}=-\psi^{-2}\psi_\alpha\sigma_4(\nabla^2v)+\psi^{-1}\sum_{i,j}\frac{\partial\sigma_4}{\partial v_{ij}}v_{ij\alpha}.$$ Since $\nabla^2v$ is diagonal, it follows from Lemma \ref{Lemma-sigma_k} that
\begin{equation}\label{Equality-q_alpha}
\frac{\partial q}{\partial x_\alpha}=-\psi^{-2}\psi_\alpha v_{11}v_{22}v_{33}v_{44}+\psi^{-1}\sum_i\sigma_3(\lambda|i)v_{ii\alpha}.
\end{equation}
In view of \eqref{Estimate-psi}, we have $\psi^2\geq 4c^2v_{33}v_{44}$, which implies $$0\leq\psi^{-2}v_{11}v_{22}v_{33}v_{44}\leq\frac{v_{11}v_{22}}{4c^2}.$$ Differentiating $\psi$ with respect to $x_\alpha$ and using again the fact that $\nabla^2v$ is diagonal, we obtain
\begin{align*}
\psi_\alpha
=&~\frac{4}{16^2}(v_{11}+v_{33})^2(v_{22}+v_{44})(v_{22\alpha}+v_{44\alpha})+\frac{4}{16^2}(v_{11}+v_{33})(v_{22}+v_{44})^2(v_{11\alpha}+v_{33\alpha})\\
&-\frac{2}{16^2}(v_{11}-v_{33})(v_{22}+v_{44})^2(v_{11\alpha}-v_{33\alpha})-\frac{2}{16^2}(v_{11}-v_{33})^2(v_{22}+v_{44})(v_{22\alpha}+v_{44\alpha})\\
&-\frac{2}{16^2}(v_{22}-v_{44})(v_{11}+v_{33})^2(v_{22\alpha}-v_{44\alpha})-\frac{2}{16^2}(v_{22}-v_{44})^2(v_{11}+v_{33})(v_{11\alpha}+v_{33\alpha}).
\end{align*}
Since $v_{33}=O(\phi)$, $v_{44}=O(\phi)$, and $v_{22}=v_{11}+O(\phi)$, it follows that
\begin{equation}\label{Equality-psi_alpha}
\begin{aligned}
\psi_\alpha
=&~\frac{4}{16^2}v_{11}^3(v_{22\alpha}+v_{44\alpha})+\frac{4}{16^2}v_{11}^3(v_{11\alpha}+v_{33\alpha})-\frac{2}{16^2}v_{11}^3(v_{11\alpha}-v_{33\alpha})\\
&-\frac{2}{16^2}v_{11}^3(v_{22\alpha}+v_{44\alpha})-\frac{2}{16^2}v_{11}^3(v_{22\alpha}-v_{44\alpha})-\frac{2}{16^2}v_{11}^3(v_{11\alpha}+v_{33\alpha})+O(\phi)\\
=&~\frac{4}{16^2}v_{11}^3(v_{33\alpha}+v_{44\alpha})+O(\phi)=O\bigg(\phi+\sum_{p,q\in B}|\nabla v_{pq}|\bigg).
\end{aligned}
\end{equation}
Therefore, for the first term on the right-hand side of \eqref{Equality-q_alpha} satisfies $$-\psi^{-2}\psi_\alpha v_{11}v_{22}v_{33}v_{44}=O\bigg(\phi+\sum_{p,q\in B}|\nabla v_{pq}|\bigg).$$ Using \eqref{Estimate-psi}, the second term can be estimated as
\begin{align*}
\psi^{-1}\sum_i\sigma_3(\lambda|i)v_{ii\alpha}
=&~\psi^{-1}(v_{22}v_{33}v_{44}v_{11\alpha}+v_{11}v_{33}v_{44}v_{22\alpha}+v_{11}v_{22}v_{44}v_{33\alpha}+v_{11}v_{22}v_{33}v_{44\alpha})\\
=&~O\bigg(\phi+\sum_{p,q\in B}|\nabla v_{pq}|\bigg).
\end{align*}
Consequently, $$\frac{\partial q}{\partial x_\alpha}=O\bigg(\phi+\sum_{p,q\in B}|\nabla v_{pq}|\bigg).$$
\end{proof}

Before computing the second-order derivatives of $q$, we recall the following lemma, which provides an estimate for the third-order derivatives of convex functions. We refer the reader to Bian and Guan \cite{Bian-Guan2009} for the proof.
\begin{lemma}[Lemma 2.5 in \cite{Bian-Guan2009}]\label{Lemma-v_{ijalpha}}
Assume $v\in C^{3,1}(\Omega)$ is a convex function. Then there exists a positive constant $C$ depending only on $\mathrm{dist}(\mathcal{O},\partial\Omega)$ and $\|v\|_{C^{3,1}(\Omega)}$ such that
\begin{equation}
|v_{ij\alpha}|\leq C\bigg(\sqrt{v_{ii}(x)}+\sqrt{v_{jj}(x)}\bigg)
\end{equation}
for all $x\in\mathcal{O}$ and $1\leq i,j,\alpha\leq n$.
\end{lemma}

We now carry out the main computation of this subsection.

\begin{proposition}\label{Proposition-q_{ij}}
Let $q(\nabla^2v)$ be defined as in \eqref{q(D^2v(x))-new}. Assume that $\nabla^2v(x)$ is diagonal for $x\in\mathcal{O}$. Then, for any $\alpha,\beta\in\{1,2,3,4\}$, we have
\begin{align*}
\frac{\partial^2q}{\partial x_\alpha\partial x_\beta}
=&~4\cdot 16v_{11}^{-1}\frac{v_{44}^2}{\sigma_1^2(B)}v_{33\alpha\beta}-8\cdot 16v_{11}^{-2}\frac{v_{44}^2}{\sigma_1^2(B)}v_{13\alpha}v_{13\beta}-8\cdot 16v_{11}^{-2}\frac{v_{44}^2}{\sigma_1^2(B)}v_{23\alpha}v_{23\beta}\\
&+4\cdot16v_{11}^{-1}\frac{v_{33}^2}{\sigma_1^2(B)}v_{44\alpha\beta}-8\cdot 16v_{11}^{-2}\frac{v_{33}^2}{\sigma_1^2(B)}v_{14\alpha}v_{14\beta}-8\cdot16v_{11}^{-2}\frac{v_{33}^2}{\sigma_1^2(B)}v_{24\alpha}v_{24\beta}\\
&-4\cdot16v_{11}^{-1}\frac{1}{\sigma_1^3(B)}(V_{3\alpha}V_{3\beta}+V_{4\alpha}V_{4\beta})-8\cdot16v_{11}^{-1}\frac{1}{\sigma_1(B)}v_{34\alpha}v_{34\beta}\\
&+O\bigg(\phi+\sum_{p,q\in B}|\nabla v_{pq}|\bigg),
\end{align*}
where $$V_{3\alpha}=v_{33\alpha}\sigma_1(B)-v_{33}\bigg(\sum_{i\in B}v_{ii\alpha}\bigg),\quad\mathrm{and}\ V_{4\alpha}=v_{44\alpha}\sigma_1(B)-v_{44}\bigg(\sum_{i\in B}v_{ii\alpha}\bigg).$$
\end{proposition}
\begin{proof}
Differentiating $q$ twice with respect to $x_\alpha$ and $x_\beta$, we have
\begin{align*}
\frac{\partial^2q}{\partial x_\alpha\partial x_\beta}
=&~2\psi^{-3}\psi_\alpha\psi_\beta\sigma_4(\nabla^2v)-\psi^{-2}\psi_{\alpha\beta}\sigma_4(\nabla^2v)-\psi^{-2}\psi_\alpha\sum_{i,j}\frac{\partial\sigma_4}{\partial v_{ij}}v_{ij\beta}\\
&-\psi^{-2}\psi_\beta\sum_{i,j}\frac{\partial\sigma_4}{\partial v_{ij}}v_{ij\alpha}+\psi^{-1}\sum_{i,j,k,l}\frac{\partial^2\sigma_4}{\partial v_{ij}\partial v_{kl}}v_{ij\alpha}v_{kl\beta}+\psi^{-1}\sum_{i,j}\frac{\partial\sigma_4}{\partial v_{ij}}v_{ij\alpha\beta}.
\end{align*}
Since $\nabla^2v$ is diagonal at the point under consideration, Lemma \ref{Lemma-sigma_k} implies that
\begin{equation}\label{Equality-q_{alphabeta}}
\begin{aligned}
\frac{\partial^2q}{\partial x_\alpha\partial x_\beta}
=&~2\psi^{-3}\psi_\alpha\psi_\beta\sigma_4(\nabla^2v)-\psi^{-2}\psi_{\alpha\beta}\sigma_4(\nabla^2v)-\psi^{-2}\psi_\alpha\sum_i\sigma_3(\lambda|i)v_{ii\beta}\\
&-\psi^{-2}\psi_\beta\sum_i\sigma_3(\lambda|i)v_{ii\alpha}+\psi^{-1}\sum_i\sigma_3(\lambda|i)v_{ii\alpha\beta}\\
&+\psi^{-1}\sum_{i,j}\big(\sigma_2(\lambda|ij)v_{ii\alpha}v_{jj\beta}-\sigma_2(\lambda|ij)v_{ij\alpha}v_{ji\beta}\big)\\
=&~\circled{1}+\circled{2}+\circled{3}+\circled{4}+\circled{5}+\circled{6}.
\end{aligned}
\end{equation}

We now analyze the terms \circled{1}-\circled{6} individually.

$\bullet$ \circled{1}. By \eqref{Equality-psi-diagonal}, together with the estimates \eqref{Estimate-v_{33},v_{44}} and \eqref{Estimate-v_{22}}, we derive $$\psi=\frac{4}{16^2}v_{11}^3(v_{33}+v_{44})+O(\phi^2).$$ Denoting $\psi_0=\frac{4}{16^2}v_{11}^3(v_{33}+v_{44})$, we have $\psi=\psi_0+O(\phi^2)$. It then follows from \eqref{Estimate-v_{33}+v_{44}=phi}  that $\frac{\psi-\psi_0}{\psi_0}=O(\phi)$. Applying the Taylor expansion, we have
\begin{align*}
\psi^{-3}
=&~\psi_0^{-3}\bigg(1+\frac{\psi-\psi_0}{\psi_0}\bigg)^{-3}\\
=&~\psi_0^{-3}\bigg(1-3\frac{\psi-\psi_0}{\psi_0}+6\bigg(\frac{\psi-\psi_0}{\psi_0}\bigg)^2+O(\phi^3)\bigg)\\
=&~\psi_0^{-3}-3\psi_0^{-4}(\psi-\psi_0)+6\psi_0^{-5}(\psi-\psi_0)^2+O(1). 
\end{align*}
In view of \eqref{Equality-psi_alpha}, we further obtain 
\begin{align*}
\psi_\alpha\psi_\beta
=&~\frac{4^2}{16^4}v_{11}^6(v_{33\alpha}+v_{44\alpha})(v_{33\beta}+v_{44\beta})\\
&+O(\phi)\bigg(\frac{4}{16^2}v_{11}^3(v_{33\alpha}+v_{44\alpha}+v_{33\beta}+v_{44\beta})\bigg)+O(\phi^2).
\end{align*}
Hence,
\begin{align*}
2\psi^{-3}\psi_\alpha\psi_\beta\sigma_4(\nabla^2v)
=&~2\big(\psi_0^{-3}-3\psi_0^{-4}(\psi-\psi_0)+6\psi_0^{-5}(\psi-\psi_0)^2+O(1)\big)v_{11}v_{22}v_{33}v_{44}\\
&\cdot\bigg(\frac{4^2}{16^4}v_{11}^6(v_{33\alpha}+v_{44\alpha})(v_{33\beta}+v_{44\beta})\\
&\hspace{0.5cm}+O(\phi)\bigg(\frac{4}{16^2}v_{11}^3(v_{33\alpha}+v_{44\alpha}+v_{33\beta}+v_{44\beta})\bigg)+O(\phi^2)\bigg)\\
=&~\bigg(\frac{2v_{11}^2v_{33}v_{44}}{\frac{4^3}{16^6}v_{11}^9(v_{33}+v_{44})^3}-6\frac{(\psi-\psi_0)v_{11}^2v_{33}v_{44}}{\psi_0^4}\bigg)\\
&\cdot\bigg(\frac{4^2}{16^4}v_{11}^6(v_{33\alpha}+v_{44\alpha})(v_{33\beta}+v_{44\beta})\\
&\hspace{0.5cm}+O(\phi)\bigg(\frac{4}{16^2}v_{11}^3(v_{33\alpha}+v_{44\alpha}+v_{33\beta}+v_{44\beta})\bigg)+O(\phi^2)\bigg)\\
&+O(\phi).
\end{align*}
On the other hand, Lemma \ref{Lemma-v_{ijalpha}} and Cauchy-Schwarz inequality imply that 
\begin{equation}\label{Estimate-v_{pqalpha}^2=O(phi)}
v_{pq\alpha}^2=O(\phi),\quad p,q\in B,\ \alpha\in\{1,2,3,4\}.
\end{equation} 
Therefore, we may estimate $$2\psi^{-3}\psi_\alpha\psi_\beta\sigma_4(\nabla^2v)=\frac{1}{2}\cdot 16^2v_{11}^{-1}\frac{\sigma_2(B)}{\sigma_1^3(B)}(v_{33\alpha}+v_{44\alpha})(v_{33\beta}+v_{44\beta})+O\bigg(\phi+\sum_{p,q\in B}|\nabla v_{pq}|\bigg).$$

$\bullet$ \circled{2}. Similar to the analysis of \circled{1}, we have
\begin{equation}\label{Estimate-psi^{-2}}
\begin{aligned}
\psi^{-2}
=&~\psi_0^{-2}\bigg(1+\frac{\psi-\psi_0}{\psi_0}\bigg)^{-2}\\
=&~\psi_0^{-2}\bigg(1-2\frac{\psi-\psi_0}{\psi_0}+O(\phi^2)\bigg)\\
=&~\psi_0^{-2}-2\psi_0^{-3}(\psi-\psi_0)+O(1).
\end{aligned}
\end{equation}
Here $\psi_0=\frac{4}{16^2}v_{11}^3(v_{33}+v_{44})$. For the term \circled{2} we obtain
\begin{align*}
 -\psi^{-2}\psi_{\alpha\beta}\sigma_4(\nabla^2v)
 =&-\left(\psi_0^{-2}-2\psi_0^{-3}(\psi-\psi_0)+O(1)\right)v_{11}v_{22}v_{33}v_{44}\psi_{\alpha\beta}\\
 =&-\frac{v_{11}^2v_{33}v_{44}}{\frac{4^2}{16^4}v_{11}^6(v_{33}+v_{44})^2}\psi_{\alpha\beta}+O(\phi)\\
 =&-16^3v_{11}^{-4}\frac{\sigma_2(B)}{\sigma_1^2(B)}\psi_{\alpha\beta}+O(\phi).
\end{align*}
It remains to calculate $\psi_{\alpha\beta}$. The diagonal form of $\nabla^2v$ allows us to compute explicitly
\begin{align*}
\psi_{\alpha\beta}
=&~\frac{4}{16^2}\big[(v_{11\alpha}+v_{33\alpha})(v_{22}+v_{44})+(v_{11}+v_{33})(v_{22\alpha}+v_{44\alpha})\big]\\
&\hspace{0.8cm}\cdot\big[(v_{11\beta}+v_{33\beta})(v_{22}+v_{44})+(v_{11}+v_{33})(v_{22\beta}+v_{44\beta})\big]\\
&+\frac{4}{16^2}(v_{11}+v_{33})(v_{22}+v_{44})\\
&\hspace{0.8cm}\cdot\big[(v_{11\alpha\beta}+v_{33\alpha\beta})(v_{22}+v_{44})+(v_{11}+v_{33})(v_{22\alpha\beta}+v_{44\alpha\beta})\\
&\hspace{1.2cm}+(v_{11\alpha}+v_{33\alpha})(v_{22\beta}+v_{44\beta})+(v_{11\beta}+v_{33\beta})(v_{22\alpha}+v_{44\alpha})\\
&\hspace{1.2cm}-2(v_{12\alpha}+v_{34\alpha})(v_{12\beta}+v_{34\beta})-2(v_{14\alpha}-v_{32\alpha})(v_{14\beta}-v_{32\beta})\big]\\
&-\frac{1}{16^2}\big\{2\big[(v_{11\alpha}-v_{33\alpha})(v_{22}+v_{44})+(v_{11}-v_{33})(v_{22\alpha}+v_{44\alpha})\big]\\
&\hspace{1.3cm}\cdot\big[(v_{11\beta}-v_{33\beta})(v_{22}+v_{44})+(v_{11}-v_{33})(v_{22\beta}+v_{44\beta})\big]\\
&\hspace{1.2cm}+2(v_{11}-v_{33})(v_{22}+v_{44})\\
&\hspace{1.4cm}\cdot\big[(v_{11\alpha\beta}-v_{33\alpha\beta})(v_{22}+v_{44})+(v_{11}-v_{33})(v_{22\alpha\beta}+v_{44\alpha\beta})\\
&\hspace{1.7cm}+(v_{11\alpha}-v_{33\alpha})(v_{22\beta}+v_{44\beta})+(v_{11\beta}-v_{33\beta})(v_{22\alpha}+v_{44\alpha})\\
&\hspace{1.7cm}+2v_{23\alpha}v_{23\beta}+2v_{34\alpha}v_{34\beta}-2v_{12\alpha}v_{12\beta}-2v_{14\alpha}v_{14\beta}\big]\\
&\hspace{1.2cm}+2\big[(v_{22\alpha}-v_{44\alpha})(v_{11}+v_{33})+(v_{22}-v_{44})(v_{11\alpha}+v_{33\alpha})\big]\\
&\hspace{1.6cm}\cdot\big[(v_{22\beta}-v_{44\beta})(v_{11}+v_{33})+(v_{22}-v_{44})(v_{11\beta}+v_{33\beta})\big]\\
&\hspace{1.2cm}+2(v_{22}-v_{44})(v_{11}+v_{33})\\
&\hspace{1.6cm}\cdot\big[(v_{22\alpha\beta}-v_{44\alpha\beta})(v_{11}+v_{33})+(v_{22}-v_{44})(v_{11\alpha\beta}+v_{33\alpha\beta})\\
&\hspace{1.8cm}+(v_{22\alpha}-v_{44\alpha})(v_{11\beta}+v_{33\beta})+(v_{22\beta}-v_{44\beta})(v_{11\alpha}+v_{33\alpha})\\
&\hspace{1.8cm}+2v_{14\alpha}v_{14\beta}+2v_{34\alpha}v_{34\beta}-2v_{12\alpha}v_{12\beta}-2v_{23\alpha}v_{23\beta}\big]\\
&\hspace{1.2cm}+8(v_{22}+v_{44})^2v_{12\alpha}v_{13\beta}+8(v_{11}+v_{33})^2v_{24\alpha}v_{24\beta}\\
&\hspace{1.2cm}+8(v_{12\alpha}v_{33}-v_{11}v_{34\alpha})(v_{21\beta}v_{44}-v_{22}v_{43\beta})\\
&\hspace{1.2cm}+8(v_{12\beta}v_{33}-v_{11}v_{34\beta})(v_{21\alpha}v_{44}-v_{22}v_{43\alpha})\\
&\hspace{1.2cm}+8(v_{11}v_{32\alpha}+v_{33}v_{14\alpha})(v_{44}v_{23\beta}+v_{22}v_{14\beta})\\
&\hspace{1.2cm}+8(v_{11}v_{32\beta}+v_{33}v_{14\beta})(v_{44}v_{23\alpha}+v_{22}v_{14\alpha})\big\}.
\end{align*}
Substituting the estimates \eqref{Estimate-v_{33},v_{44}}, \eqref{Estimate-v_{22}}, and \eqref{Estimate-v_{pqalpha}^2=O(phi)} into the above formula, we deduce the following estimate for $\psi_{\alpha\beta}$
\begin{equation}\label{Equality-psi_{alphabeta}}
\begin{aligned}
\psi_{\alpha\beta}
=&~\frac{4}{16^2}v_{11}^3(v_{33\alpha\beta}+v_{44\alpha\beta})-\frac{8}{16^2}v_{11}^2v_{13\alpha}v_{13\beta}-\frac{8}{16^2}v_{11}^2v_{14\alpha}v_{14\beta}\\
&-\frac{8}{16^2}v_{11}^2v_{32\alpha}v_{32\beta}-\frac{8}{16^2}v_{11}^2v_{24\alpha}v_{24\beta}+O\bigg(\phi+\sum_{p,q\in B}|\nabla v_{pq}|\bigg).
\end{aligned}
\end{equation}
Thus,
\begin{align*}
-\psi^{-2}\psi_{\alpha\beta}\sigma_4(\nabla^2v)
=&-16^3v_{11}^{-4}\frac{\sigma_2(B)}{\sigma_1^2(B)}\\
&~\cdot\bigg(\frac{4}{16^2}v_{11}^3(v_{33\alpha\beta}+v_{44\alpha\beta})-\frac{8}{16^2}v_{11}^2v_{13\alpha}v_{13\beta}-\frac{8}{16^2}v_{11}^2v_{14\alpha}v_{14\beta}\\
&\hspace{0.6cm}-\frac{8}{16^2}v_{11}^2v_{32\alpha}v_{32\beta}-\frac{8}{16^2}v_{11}^2v_{24\alpha}v_{24\beta}\bigg)+O\bigg(\phi+\sum_{p,q\in B}|\nabla v_{pq}|\bigg).
\end{align*}

$\bullet$ \circled{3} and \circled{4}. For the term \circled{3}, combining \eqref{Equality-psi_alpha} and \eqref{Estimate-psi^{-2}} with the estimates \eqref{Estimate-v_{33},v_{44}}, \eqref{Estimate-v_{22}}, and \eqref{Estimate-v_{pqalpha}^2=O(phi)}, yields
\begin{align*}
-\psi^{-2}\psi_\alpha\sigma_3(\lambda|i)v_{ii\beta}
=&-\big(\psi_0^{-2}-2\psi_0^{-3}(\psi-\psi_0)+O(1)\big)\cdot\bigg(\frac{4}{16^2}v_{11}^3(v_{33\alpha}+v_{44\alpha})+O(\phi)\bigg)\\
&~\cdot\bigg(v_{11}v_{33}v_{44}v_{11\beta}+v_{11}v_{33}v_{44}v_{22\beta}+v_{11}^2v_{44}v_{33\beta}+v_{11}^2v_{33}v_{44\beta}+O(\phi^2)\bigg)\\
=&-\frac{1}{\frac{4^2}{16^4}v_{11}^6(v_{33}+v_{44})^2}\cdot\frac{4}{16^2}v_{11}^3(v_{33\alpha}+v_{44\alpha})\cdot\big(v_{11}^2v_{44}v_{33\beta}+v_{11}^2v_{33}v_{44\beta}\big)\\[6pt]
&~+O\bigg(\phi+\sum_{p,q\in B}|\nabla v_{pq}|\bigg)\\
=&-\frac{1}{4}\cdot16^2v_{11}^{-1}\frac{1}{\sigma_1^2(B)}\bigg(v_{44}v_{33\beta}(v_{33\alpha}+v_{44\alpha})+v_{33}v_{44\beta}(v_{33\alpha}+v_{44\alpha})\bigg)\\[6pt]
&~+O\bigg(\phi+\sum_{p,q\in B}|\nabla v_{pq}|\bigg).
\end{align*}
The term \circled{4} can be treated in the same way, and therefore
\begin{align*}
-\psi^{-2}\psi_\beta\sigma_3(\lambda|i)v_{ii\alpha}
=&-\frac{1}{4}\cdot16^2v_{11}^{-1}\frac{1}{\sigma_1^2(B)}\bigg(v_{44}v_{33\alpha}(v_{33\beta}+v_{44\beta})+v_{33}v_{44\alpha}(v_{33\beta}+v_{44\beta})\bigg)\\[6pt]
&+O\bigg(\phi+\sum_{p,q\in B}|\nabla v_{pq}|\bigg).
\end{align*}

$\bullet$ \circled{5}. Using the Taylor expansion of $\psi^{-1}$, we have
\begin{equation}\label{Estimate-psi^{-1}}
\begin{aligned}
\psi^{-1}
=&~\psi_0^{-1}\bigg(1+\frac{\psi-\psi_0}{\psi_0}\bigg)^{-1}\\
=&~\psi_0^{-1}\bigg(1-\frac{\psi-\psi_0}{\psi_0}+O(\phi^2)\bigg)\\
=&~\psi_0^{-1}-\psi_0^{-2}(\psi-\psi_0)+O(\phi).
\end{aligned}
\end{equation}
Proceeding as in the analysis of the previous terms, we arrive at
\begin{align*}
\psi^{-1}\sigma_3(\lambda|i)v_{ii\alpha\beta}
=&~\psi_0^{-1}\\
&\cdot\big(v_{11}v_{33}v_{44}v_{11\alpha\beta}+v_{11}v_{33}v_{44}v_{22\alpha\beta}+v_{11}^2v_{44}v_{33\alpha\beta}+v_{11}^2v_{33}v_{44\alpha\beta}+O(\phi^2)\big)\\
&+O(\phi)\\
=&~\frac{1}{4}\cdot 16^2v_{11}^{-1}\frac{1}{\sigma_1(B)}(v_{44}v_{33\alpha\beta}+v_{33}v_{44\alpha\beta})+O(\phi)
\end{align*}

$\bullet$ \circled{6}. Invoking \eqref{Estimate-psi^{-1}} together with the estimates \eqref{Estimate-v_{33},v_{44}}, \eqref{Estimate-v_{22}}, and \eqref{Estimate-v_{pqalpha}^2=O(phi)}, a straightforward computation gives
\begin{align*}
&~\psi^{-1}\sum_{i,j}\left(\sigma_2(\lambda|ij)v_{ii\alpha}v_{jj\beta}-\sigma_2(\lambda|ij)v_{ij\alpha}v_{ji\beta}\right)\\
=&~\frac{1}{4}\cdot16^2v_{11}^{-3}\frac{1}{\sigma_1(B)}\\
&\cdot\big[v_{11}^2(v_{33\alpha}v_{44\beta}+v_{33\beta}v_{44\alpha})+v_{11}v_{33}(v_{22\alpha}v_{44\beta}+v_{22\beta}v_{44\alpha})\\
&\hspace{0.4cm}+v_{11}v_{44}(v_{22\alpha}v_{33\beta}+v_{22\beta}v_{33\alpha})+v_{11}v_{33}(v_{11\alpha}v_{44\beta}+v_{11\beta}v_{44\alpha})\\
&\hspace{0.4cm}+v_{11}v_{44}(v_{11\alpha}v_{33\beta}+v_{11\beta}v_{33\alpha})+v_{33}v_{44}(v_{11\alpha}v_{22\beta}+v_{11\beta}v_{22\alpha})\\
&\hspace{0.4cm}-2v_{11}^2v_{34\alpha}v_{34\beta}-2v_{11}v_{33}v_{24\alpha}v_{24\beta}-2v_{11}v_{44}v_{23\alpha}v_{23\beta}-2v_{11}v_{33}v_{14\alpha}v_{14\beta}\\
&\hspace{0.4cm}-2v_{11}v_{44}v_{13\alpha}v_{13\beta}-2v_{33}v_{44}v_{12\alpha}v_{12\beta}\big]\\
&+O(\phi)\\
=&~\frac{1}{4}\cdot16^2v_{11}^{-3}\frac{1}{\sigma_1(B)}\\
&\cdot\big[v_{11}^2(v_{33\alpha}v_{44\beta}+v_{33\beta}v_{44\alpha}-2v_{34\alpha}v_{34\beta})\\
&\hspace{0.4cm}-2v_{11}v_{33}(v_{14\alpha}v_{14\beta}+v_{24\alpha}v_{24\beta})-2v_{11}v_{44}(v_{13\alpha}v_{13\beta}+v_{23\alpha}v_{23\beta})\big]\\
&+O\bigg(\phi+\sum_{p,q\in B}|\nabla v_{pq}|\bigg).
\end{align*}

Combining the estimates obtained above for \circled{1}-\circled{6}, \eqref{Equality-q_{alphabeta}} can be rewritten as
\begin{align*}
\frac{\partial^2q}{\partial x_\alpha\partial x_\beta}
=&~4\cdot 16v_{11}^{-1}\frac{v_{44}^2}{\sigma_1^2(B)}v_{33\alpha\beta}-8\cdot 16v_{11}^{-2}\frac{v_{44}^2}{\sigma_1^2(B)}v_{13\alpha}v_{13\beta}-8\cdot 16v_{11}^{-2}\frac{v_{44}^2}{\sigma_1^2(B)}v_{23\alpha}v_{23\beta}\\
&+4\cdot16v_{11}^{-1}\frac{v_{33}^2}{\sigma_1^2(B)}v_{44\alpha\beta}-8\cdot 16v_{11}^{-2}\frac{v_{33}^2}{\sigma_1^2(B)}v_{14\alpha}v_{14\beta}-8\cdot16v_{11}^{-2}\frac{v_{33}^2}{\sigma_1^2(B)}v_{24\alpha}v_{24\beta}\\
&-4\cdot16v_{11}^{-1}\frac{1}{\sigma_1^3(B)}(V_{3\alpha}V_{3\beta}+V_{4\alpha}V_{4\beta})-8\cdot16v_{11}^{-1}\frac{1}{\sigma_1(B)}v_{34\alpha}v_{34\beta}\\
&+O\bigg(\phi+\sum_{p,q\in B}|\nabla v_{pq}|\bigg),
\end{align*}
where $$V_{3\alpha}=v_{33\alpha}\sigma_1(B)-v_{33}\bigg(\sum_{i\in B}v_{ii\alpha}\bigg),\quad\mathrm{and}\ V_{4\alpha}=v_{44\alpha}\sigma_1(B)-v_{44}\bigg(\sum_{i\in B}v_{ii\alpha}\bigg).$$
The proof of the proposition is now complete.
\end{proof}

To establish the differential inequality in the next subsection, we shall need the following lemma from Bian and Guan \cite{Bian-Guan2009}. The proof can be found in their original work.
\begin{lemma}[Lemma 3.3 in \cite{Bian-Guan2009}]\label{Lemma-v_{ijk}}
Let $M$ be a positive constant such that $0<\lambda_i\leq M$ and $\frac{1}{M}\leq\gamma_i\leq M$ for $i=l+1, l+2, \cdots, n$. Assume that $v_{ij\alpha}=v_{ji\alpha}$ for $i, j=l+1, l+2, \cdots, n$ and $\alpha=1, 2, \cdots, n$. Then there exists a positive constant $C$ depending only on $n$ and $M$, such that for each $\alpha$, for any $D>0$ and any $\delta>0$, we have
\begin{align*}
\sum_{l+1\leq i, j\leq n}|v_{ij\alpha}|
\leq&~C\left(1+\frac{2D}{\delta}+D\right)\left(\sigma_1(\lambda) +\left|\sum_{l+1\leq i\leq n}\gamma_iv_{ii\alpha}\right|\right)\\
&+\frac{\delta}{2D}\frac{1}{\sigma_1(\lambda)}\sum_{\substack{l+1\leq i, j\leq n\\ i\neq j}}|v_{ij\alpha}|^2+\frac{C}{D}\frac{1}{\sigma_1^3(\lambda)}
\sum_{l+1\leq i\leq n}V_{i\alpha}^2,
\end{align*}
where $$\sigma_1(\lambda)=\sum_{l+1\leq i\leq n}\lambda_i \quad \mathrm{and} \ V_{i\alpha}=v_{ii\alpha}\sigma_1(\lambda)-\lambda_i\bigg(\sum_{l+1\leq j\leq n}v_{jj\alpha}\bigg). $$
\end{lemma}

\subsection{A strong maximum principle}

The purpose of this subsection is to establish a strong maximum principle for the auxiliary function $\phi$ defined in \eqref{Auxiliary function-Bian and Guan-K}, associated with equation \eqref{Equ-F-v}, which serves as a key ingredient in the proof of the constant rank theorem. More precisely, we derive a suitable differential inequality for $\sum_{i,j}F^{ij}\phi_{ij}$, where $F^{ij}$ denotes the linearized operator of \eqref{Equ-F-v}. This estimates will allow us to apply the strong maximum principle to $\phi$.

\begin{theorem}\label{Theorem-rank2}
Suppose $\Omega\subset\mathbb{R}^4$ is a domain and $v\in C^4(\Omega)$ is a convex solution of equation \eqref{Equ-real equation(v)}. Let $K(x)$ be the matrix defined by \eqref{Auxiliary matrix K}. For each $x\in\Omega$, set $$\phi(x)=\mathrm{tr}(K(x))+q(K(x)),$$ where $$q(K(x))=\begin{cases}
\frac{\det(K(x))}{\mathrm{tr}(K(x))} & \mathrm{if}\ \mathrm{tr}(K(x))>0\\
0 & \mathrm{if}\ \mathrm{tr}(K(x))=0.
\end{cases}$$ If the Hessian $\nabla^2v(x)$ attains the minimal rank $2$ at some point $x_0\in\Omega$, then there exists a neighborhood $\mathcal{O}$ of $x_0$ and positive constant $C$ independent of $\phi$, such that
\begin{equation}\label{Inequ-Theorem-minimal rank=2}
\sum_{i,j}F^{ij}\phi_{ij}\leq C(\phi+|\nabla\phi|)\quad \mathrm{in}\ \mathcal{O}.
\end{equation}
\end{theorem}
\begin{proof}
As in the previous subsection, we introduce the perturbation $v_\varepsilon(x)=v(x)+\frac{\varepsilon}{2}|x|^2$ to ensure that $\mathrm{tr}(K(x))$ remains strictly positive. We denote by $K_\varepsilon$ the matrix $K$ associated with $v_\varepsilon$. Then, the function $v_\varepsilon(x)$ satisfies equation
\begin{equation}\label{Equ-F-varepsilon}
F(\nabla^2v_\varepsilon(x),\nabla v_\varepsilon(x),v_\varepsilon(x))=R_\varepsilon(x),
\end{equation}
where $R_\varepsilon=F(\nabla^2v_\varepsilon,\nabla v_\varepsilon,v_\varepsilon)-F(\nabla^2v,\nabla v,v)$ and $F(\nabla^2v,\nabla v,v)$ is defined as in \eqref{Equ-Dirichlet problem-real}. Since $v\in C^4(\Omega)$, for any $x\in\mathcal{O}$, we have
\begin{equation}\label{Estimate-R}
|R_\varepsilon(x)|\leq C\varepsilon,\quad |\nabla R_\varepsilon(x)|\leq C\varepsilon\ \ \mathrm{and}\ \  |\nabla^2R_\varepsilon(x)|\leq C\varepsilon.
\end{equation}
Recalling \eqref{Estimate-Perturbation}, namely
\begin{equation}\label{Estimate-varepsilon}
\mathrm{tr}(K_\varepsilon(x))\geq C\varepsilon\quad \mbox{for all}\ x\in\mathcal{O},
\end{equation}
where $C>0$ is independent of $\varepsilon$. In view of \eqref{Auxiliary function-Bian and Guan-K-1}, we may set $$\phi_\varepsilon(x)=(\det(A_\varepsilon))^{-2}\psi_\varepsilon(x)+q_\varepsilon(\nabla^2v_\varepsilon(x)),$$ where $$q_\varepsilon(\nabla^2v_\varepsilon(x))=\frac{\sigma_4(\nabla^2v_\varepsilon(x))}{\psi_\varepsilon(x)}.$$ Here $A_\varepsilon$ and $\psi_\varepsilon$ are defined by \eqref{Matrix-A} and \eqref{Equality-psi}, respectively, with $v$ replaced by $v_\varepsilon$. Working on equation \eqref{Equ-F-varepsilon}, we will establish a differential inequality analogous to \eqref{Inequ-Theorem-minimal rank=2} for $v_\varepsilon(x)$. By letting $\varepsilon\rightarrow 0$, we obtain \eqref{Inequ-Theorem-minimal rank=2}.

We will continue to denote $v_\varepsilon$, $\nabla^2v_\varepsilon$, $R_\varepsilon$, $\phi_\varepsilon$, $A_\varepsilon$, $q_\varepsilon$, $\psi_\varepsilon$, and so on by $v$, $\nabla^2v_\varepsilon$, $R_\varepsilon$,  $\phi_\varepsilon$, $A_\varepsilon$, $q_\varepsilon$,  $\psi_\varepsilon$, respectively, with the understanding that all estimates are independent of $\varepsilon$. Then, $v$ satisfies equation
\begin{equation}\label{Equation-minimal rank is 2}
F(\nabla^2v,\nabla v,v)=R(x),
\end{equation}
where $F$ is defined as in \eqref{Equ-Dirichlet problem-real}.

Fix a point $x\in\mathcal{O}$ and choose the coordinates described in Section \ref{Section-auxiliary function}. We assume that $\nabla^2v(x)$ is diagonal, satisfies $v_{11}+v_{33}=v_{22}+v_{44}$, and that the eigenvalues are ordered as $v_{11}(x)\geq v_{22}(x)\geq v_{33}(x)\geq v_{44}(x)\geq 0$. By \eqref{Estimate-varepsilon}, we have $$\varepsilon\leq C\phi(x).$$ Furthermore, combining this relation with \eqref{Estimate-R}, we know that
\begin{equation}\label{Estimate-R1}
R(x)=O(\phi(x)),\quad \nabla R(x)=O(\phi(x))\ \ \mathrm{and}\ \ \nabla^2R(x)=O(\phi(x)).
\end{equation}
Below, all computations will be performed at the fixed point $x$. Since the minimal rank of $\nabla^2v$ is $2$, we denote by $G=\{1,2\}$ and $B=\{3,4\}$. Recalling the analysis in Section \ref{Section 3.1}, we know that $v_{ii}\geq c$ for $i\in G$, for some positive constant $c$. For convenience, we also restate the previously established estimates \eqref{Estimate-v_{33},v_{44}} and \eqref{Estimate-v_{22}} as
\begin{equation}\label{Estimate-v_{33},v_{44}'}
v_{33}=O(\phi),\quad v_{44}=O(\phi)
\end{equation}
and
\begin{equation}\label{Estimate-v_{22}'}
v_{22}=v_{11}+O(\phi).
\end{equation}

Taking the first derivatives of $\phi$, at the point $x$, by Lemma \ref{Lemma-first derivative of q}, we get
\begin{align*}
\phi_i
=&-2\cdot16^2(v_{11}+v_{33})^{-3}(v_{22}+v_{44})^{-2}(v_{11i}+v_{33i})\psi\\
&-2\cdot16^2(v_{11}+v_{33})^{-2}(v_{22}+v_{44})^{-3}(v_{22i}+v_{44i})\psi\\
&+16^2(v_{11}+v_{33})^{-2}(v_{22}+v_{44})^{-2}\psi_i+O\bigg(\phi+\sum_{\alpha,\beta\in B}|\nabla v_{\alpha\beta}|\bigg)
\end{align*}
Noting that $v_{11}+v_{33}$ and $v_{22}+v_{44}$ are both bounded, $\psi=O(\phi)$, and recalling \eqref{Equality-psi_alpha}, i.e., 
\begin{equation}\label{Equality-psi_i}
\psi_i=O\bigg(\phi+\sum_{\alpha,\beta\in B}|\nabla v_{\alpha\beta}|\bigg),
\end{equation}
we conclude that
\begin{equation}\label{Equality-phi_i-2}
\phi_i=O\bigg(\phi+\sum_{\alpha,\beta\in B}|\nabla v_{\alpha\beta}|\bigg).
\end{equation}

The second derivatives of $\phi$ evaluated at $x$ is
\begin{align*}
\phi_{ij}
=&~6\cdot16^2(v_{11}+v_{33})^{-4}(v_{22}+v_{44})^{-4}\\
&\cdot\big[(v_{11i}+v_{33i})(v_{22}+v_{44})+(v_{11}+v_{33})(v_{22i}+v_{44i})\big]\\
&\cdot\big[(v_{11j}+v_{33j})(v_{22}+v_{44})+(v_{11}+v_{33})(v_{22j}+v_{33j})\big]\psi\\
&-2\cdot16^2(v_{11}+v_{33})^{-3}(v_{22}+v_{44})^{-3}\\
&\cdot\big[(v_{11ij}+v_{33ij})(v_{22}+v_{44})+(v_{11}+v_{33})(v_{22ij}+v_{44ij})\\
&\hspace{0.4cm}+(v_{11i}v_{33i})(v_{22j}+v_{44j})+(v_{11j}+v_{33j})(v_{22i}+v_{44i})\\
&\hspace{0.4cm}-2(v_{12i}+v_{34i})(v_{12j}+v_{34j})-2(v_{14i}-v_{32i})(v_{14j}-v_{32j})\big]\psi\\
&-2\cdot16^2(v_{11}+v_{33})^{-3}(v_{22}+v_{44})^{-3}\\
&\cdot\big[(v_{11i}+v_{33i})(v_{22}+v_{44})+(v_{11}+v_{33})(v_{22i}+v_{44i})\big]\psi_j\\
&-2\cdot16^2(v_{11}+v_{33})^{-3}(v_{22}+v_{44})^{-3}\\
&\cdot\big[(v_{11j}+v_{33j})(v_{22}+v_{44})+(v_{11}+v_{33})(v_{22j}+v_{44j})\big]\psi_i\\
&+16^2(v_{11}+v_{33})^{-2}(v_{22}+v_{44})^{-2}\psi_{ij}+q_{ij}.
\end{align*}
For the first and second terms on the right-hand side of the above equation, using the boundedness of $v_{11}+v_{33}$ and $v_{22}+v_{44}$, together with $\psi=O(\phi)$, we see that both are $O(\phi)$. For the third and fourth terms, similarly, by \eqref{Equality-psi_i}, they are $O\big(\phi+\sum_{\alpha,\beta\in B}|\nabla v_{\alpha\beta}|\big)$. For the fifth term, applying the estimates \eqref{Estimate-v_{33},v_{44}'} and \eqref{Estimate-v_{22}'}, we obtain $$16^2(v_{11}+v_{33})^{-2}(v_{22}+v_{44})^{-2}=16^2v_{11}^{-4}+O(\phi).$$ Regarding $\psi_{ij}$, we have already computed it in the proof of Proposition \ref{Proposition-q_{ij}}, see \eqref{Equality-psi_{alphabeta}}. It follows that
\begin{align*}
16^2(v_{11}+v_{33})^{-2}(v_{22}+v_{44})^{-2}\psi_{ij}=&~4v_{11}^{-1}(v_{33ij}+v_{44ij})-8v_{11}^{-2}v_{13i}v_{13j}-8v_{11}^{-2}v_{14i}v_{14j}\\
&-8v_{11}^{-2}v_{32i}v_{32j}-8v_{11}^{-2}v_{24i}v_{24j}+O\bigg(\phi+\sum_{\alpha,\beta\in B}|\nabla v_{\alpha\beta}|\bigg).
\end{align*}
Combining these results with Proposition \ref{Proposition-q_{ij}}, we finally obtain
\begin{align*}
\phi_{ij}
=&~4v_{11}^{-1}\left(1+16\frac{v_{44}^2}{\sigma_1^2(B)}\right)v_{33ij}+4v_{11}^{-1}\left(1+16\frac{v_{33}^2}{\sigma_1^2(B)}\right)v_{44ij}\\
&-8v_{11}^{-2}\left(1+16\frac{v_{44}^2}{\sigma_1^2(B)}\right)v_{13i}v_{13j}-8v_{11}^{-2}\left(1+16\frac{v_{33}^2}{\sigma_1^2(B)}\right)v_{14i}v_{14j}\\
&-8v_{11}^{-2}\left(1+16\frac{v_{44}^2}{\sigma_1^2(B)}\right)v_{23i}v_{23j}-8v_{11}^{-2}\left(1+16\frac{v_{33}^2}{\sigma_1^2(B)}\right)v_{24i}v_{24j}\\
&-4\cdot16v_{11}^{-1}\frac{1}{\sigma_1^3(B)}(V_{3i}V_{3j}+V_{4i}V_{4j})-8\cdot16v_{11}^{-1}\frac{1}{\sigma_1(B)}v_{34i}v_{34j}\\
&+O\bigg(\phi+\sum_{\alpha,\beta\in B}|\nabla v_{\alpha\beta}|\bigg),
\end{align*}
where $$V_{3i}=v_{33i}\sigma_1(B)-v_{33}\bigg(\sum_{\gamma\in B}v_{\gamma\gamma i}\bigg),\quad\mathrm{and}\ V_{4i}=v_{44i}\sigma_1(B)-v_{44}\bigg(\sum_{\gamma\in B}v_{\gamma\gamma i}\bigg).$$ Thus,
\begin{equation}\label{Equality-F1-2}
\begin{aligned}
\sum_{i,j}F^{ij}\phi_{ij}
=&~4v_{11}^{-1}\sum_{i,j}\sum_{\alpha\in B}\left(1+16\frac{\sigma_1^2(B|\alpha)}{\sigma_1^2(B)}\right)F^{ij}v_{\alpha\alpha ij}\\
&-8v_{11}^{-2}\sum_{i,j}\sum_{\alpha\in B,\beta\in G}\left(1+16\frac{\sigma_1^2(B|\alpha)}{\sigma_1^2(B)}\right)F^{ij}v_{\alpha\beta i}v_{\alpha\beta j}\\
&-4\cdot 16v_{11}^{-1}\frac{1}{\sigma_1^3(B)}\sum_{i,j}\sum_{\alpha\in B}F^{ij}V_{\alpha i}V_{\alpha j}\\
&-4\cdot16v_{11}^{-1}\frac{1}{\sigma_1(B)}\sum_{i,j}\sum_{\substack{\alpha,\beta\in B\\ \alpha\neq\beta}}F^{ij}v_{\alpha\beta i}v_{\alpha\beta j}\\
&+O\bigg(\phi+\sum_{\alpha,\beta\in B}|\nabla v_{\alpha\beta}|\bigg),
\end{aligned}
\end{equation}
where $$V_{\alpha i}=v_{\alpha\alpha i}\sigma_1(B) -v_{\alpha\alpha}\bigg( \sum_{\gamma\in B}v_{\gamma\gamma i}\bigg).$$

To deal with the terms involving fourth-order derivatives in \eqref{Equality-F1-2}, we take the first derivatives of the equation \eqref{Equation-minimal rank is 2} to obtain
\begin{equation}\label{Relation-1st condition-2}
\sum_{i,j}F^{ij}v_{ij\alpha}+\sum_kF^{v_k}v_{k\alpha}+F^vv_\alpha=O(\phi),
\end{equation}
and the second derivatives to obtain
\begin{align*}
O(\phi)=&\sum_{i, j, k, l}F^{ij, kl}v_{ij\alpha}v_{kl\alpha} +\sum_{i, j, k}F^{ij, v_k}v_{ij\alpha}v_{k\alpha} +\sum_{i, j}F^{ij, v}v_{ij\alpha}v_\alpha +\sum_{i, j}F^{ij}v_{ij\alpha\alpha}\\
&+\sum_{k, i, j}F^{v_k, ij}v_{k\alpha}v_{ij\alpha}+\sum_{k, l}F^{v_k, v_l}v_{k\alpha}v_{l\alpha} +\sum_kF^{v_k,v}v_{k\alpha}v_\alpha +\sum_kF^{v_k}v_{k\alpha\alpha}\\
&+\sum_{i, j}F^{v, ij}v_\alpha v_{ij\alpha}+\sum_k F^{v, v_k}v_\alpha v_{k\alpha} +F^{v, v}v_\alpha v_\alpha +F^vv_{\alpha\alpha},
\end{align*}
where we used the estimates in \eqref{Estimate-R1}. Therefore, combining the estimates \eqref{Estimate-v_{33},v_{44}'}, \eqref{Estimate-v_{22}'} and Lemma \ref{Lemma-v_{ijalpha}}, we deduce that
\begin{equation}\label{Equality-4th order derivative}
\begin{aligned}
&~4v_{11}^{-1}\sum_{i,j}\sum_{\alpha\in B}\left(1+16\frac{\sigma_1^2(B|\alpha)}{\sigma_1^2(B)}\right)F^{ij}v_{\alpha\alpha ij}\\
=&-4v_{11}^{-1}\sum_{\alpha\in B}\left(1+16\frac{\sigma_1^2(B|\alpha)}{\sigma_1^2(B)}\right)\\
&\cdot\bigg(\sum_{i,j,k,l}F^{ij,kl}v_{ij\alpha}v_{kl\alpha}+2\sum_{i,j,k}F^{ij,v_k}v_{ij\alpha}v_{k\alpha}+2\sum_{i,j}F^{ij,v}v_{ij\alpha}v_\alpha\\
&\hspace{0.5cm}+\sum_{k,l}F^{v_k,v_l}v_{k\alpha}v_{l\alpha}+2\sum_kF^{v_k,v}v_{k\alpha}v_\alpha+F^{v,v}v_\alpha^2+\sum_kF^{v_k}v_{k\alpha\alpha}\\
&\hspace{0.5cm}+F^vv_{\alpha\alpha}+O(\phi)\bigg)\\
=&-4v_{11}^{-1}\sum_{\alpha\in B}\left(1+16\frac{\sigma_1^2(B|\alpha)}{\sigma_1^2(B)}\right)\\
&\cdot\bigg(\sum_{i,j,k,l\in G}F^{ij,kl}v_{ij\alpha}v_{kl\alpha}+2\sum_{i,j\in G}F^{ij,v}v_{ij\alpha}v_\alpha+F^{v,v}v_\alpha^2\bigg)\\
&+O\bigg(\phi+\sum_{\alpha,\beta\in B}|\nabla v_{\alpha\beta}|\bigg).
\end{aligned}
\end{equation}
Substituting \eqref{Equality-4th order derivative} into \eqref{Equality-F1-2} yields
\begin{equation}\label{Equality-F2-2}
\begin{aligned}
\sum_{i,j}F^{ij}\phi_{ij}
=&~-4v_{11}^{-1}\sum_{\alpha\in B}\left(1+16\frac{\sigma_1^2(B|\alpha)}{\sigma_1^2(B)}\right)\\
&\cdot\bigg(\sum_{i,j,k,l\in G}F^{ij,kl}v_{ij\alpha}v_{kl\alpha}+2\sum_{i,j\in G}F^{ij,v}v_{ij\alpha}v_\alpha+F^{v,v}v_\alpha^2\\
&\hspace{0.5cm}+2v_{11}^{-1}\sum_{i,j,\beta\in G}F^{ij}v_{\alpha\beta i}v_{\alpha\beta j}\bigg)\\
&-4\cdot16v_{11}^{-1}\bigg(\frac{1}{\sigma_1^3(B)}\sum_{i,j}\sum_{\alpha\in B}F^{ij}V_{\alpha i}V_{\alpha j}+\frac{1}{\sigma_1(B)}\sum_{i,j}\sum_{\substack{\alpha,\beta\in B\\ \alpha\neq\beta}}F^{ij}v_{\alpha\beta i}v_{\alpha\beta j}\bigg)\\
&+O\bigg(\phi+\sum_{\alpha,\beta\in B}|\nabla v_{\alpha\beta}|\bigg).
\end{aligned}
\end{equation}

To establish the differential inequality \eqref{Inequ-Theorem-minimal rank=2}, we present the following claim.

\textbf{Claim 1. } For each $\alpha\in B$, we have
\begin{equation}\label{Claim}
\begin{aligned}
&\sum_{i, j, k, l\in G}F^{ij, kl}v_{ij\alpha}v_{kl\alpha} +2\sum_{i, j\in G}F^{ij, v} v_{ij\alpha}v_\alpha +F^{v, v}v_\alpha^2 +2\sum_{i, j, \beta\in G}\frac{1}{v_{\beta\beta}}F^{ij}v_{\alpha\beta i}v_{\alpha\beta j}\\ 
\geq&~ -C\bigg(\phi+\sum_{\alpha, \beta\in B}|\nabla v_{\alpha\beta}|\bigg).
\end{aligned}
\end{equation}
If we have proved \textbf{Claim 1}, since $v_{22}=v_{11}+O(\phi)$, then \eqref{Equality-F2-2} gives
\begin{align*}
\sum_{i, j}F^{ij}\phi_{ij} \leq&-64v_{11}^{-1}\bigg(\frac{1}{\sigma_1^3(B)}\sum_{i, j} \sum_{\alpha\in B}F^{ij}V_{\alpha i}V_{\alpha j} +\frac{1}{\sigma_1(B)} \sum_{i,j} \sum_{\substack{\alpha, \beta\in B\\ \alpha\neq\beta}} F^{ij}v_{\alpha\beta i}v_{\alpha\beta j}\bigg)\\
&+C\bigg(\phi+\sum_{\alpha, \beta\in B}|\nabla v_{\alpha\beta}|\bigg),
\end{align*}
where $V_{\alpha i}=v_{\alpha\alpha i}\sigma_1(B) -v_{\alpha\alpha}\left( \sum_{\gamma\in B} v_{\gamma\gamma i}\right)$. Since $v\in C^4(\Omega)$ and $\bar{\mathcal{O}}\subset\Omega$, there exists a constant $\delta_0>0$ such that $$(F^{ij})\geq \delta_0I_4 \ \ \mathrm{in}\ \ \mathcal{O}. $$ Consequently, we can deduce $$\sum_{i, j}F^{ij}V_{\alpha i}V_{\alpha j}\geq \delta_0\sum_i V_{\alpha i}^2 \ \ \mathrm{and} \ \ \sum_{i, j} F^{ij} v_{\alpha\beta i}v_{\beta\alpha j}\geq \delta_0\sum_i v_{\alpha\beta i}^2. $$
It follows that
\begin{equation*}
\sum_{i, j}F^{ij}\phi_{ij}\leq -\frac{\delta'_0}{\sigma_1^3(B)}\sum_i\sum_{\alpha\in B}V_{\alpha i}^2 -\frac{\delta'_0}{\sigma_1(B)}\sum_i\sum_{\substack{\alpha, \beta\in B\\ \alpha\neq\beta}}v_{\alpha\beta i}^2 +C\bigg(\phi+\sum_{\alpha,\beta\in B}|\nabla v_{\alpha\beta}|\bigg).
\end{equation*}
By Lemma \ref{Lemma-v_{ijk}} and equation \eqref{Equality-phi_i-2}, we finally obtain
$$\sum_{i, j}F^{ij}\phi_{ij}\leq C\left(\phi+|\nabla\phi|\right). $$
Since the constant $C$ is independent of $\varepsilon$, letting $\varepsilon\rightarrow 0$ yields \eqref{Inequ-Theorem-minimal rank=2} for $v$.
\end{proof}

\textbf{Proof of Claim 1. }  Recall the definition of $F$, namely, 
\begin{align*}
F(\nabla^2v, \nabla v, v) = &~v^2\left[(v_{11}+v_{33})(v_{22}+v_{44}) -(v_{12}+v_{34})(v_{21}+v_{43}) -(v_{14}-v_{32})(v_{41}-v_{23})\right]\\
&+v\left[(v_2^2+v_4^2)(v_{11}+v_{33}) +(v_1^2+v_3^2)(v_{22}+v_{44})\right.\\
&\hspace{8mm}\left.-(v_1v_2+v_3v_4)(v_{12}+v_{34}+v_{21}+v_{43})\right.\\
&\hspace{8mm}\left.-(v_1v_4-v_2v_3)(v_{14}-v_{32}+v_{41}-v_{23})\right]-1.
\end{align*}
It is clear that
\begin{align*}
&F^{11, 22}=F^{22, 11}=v^2, \quad F^{12, 21}=F^{21, 12}=-v^2, \\
&F^{11, v}=2v(v_{22}+v_{44}) +v_2^2+v_4^2 =\frac{F^{11}}{v} +v(v_{22}+v_{44}), \\[5pt]
&F^{22, v}=2v(v_{11}+v_{33}) +v_1^2+v_3^2 =\frac{F^{22}}{v} +v(v_{11}+v_{33}),
\intertext{and}
&F^{12, v}=-(v_1v_2+v_3v_4) =\frac{F^{12}}{v}, \quad F^{v, v}=2(v_{11}+v_{33})(v_{22}+v_{44}).
\end{align*}
Note that $G=\{1, 2\}$ and $B=\{3, 4\}$. Expanding and collecting the terms on the left-hand side of \eqref{Claim}, we obtain
\begin{equation}\label{Equality-F3-2}
\begin{aligned}
&\sum_{i, j, k, l\in G}F^{ij, kl}v_{ij\alpha}v_{kl\alpha} +2\sum_{i, j\in G}F^{ij, v}v_{ij\alpha}v_\alpha +F^{v, v}v_\alpha^2 +2\sum_{i, j, \beta\in G} \frac{1}{v_{\beta\beta}} F^{ij}v_{\alpha\beta i}v_{\beta\alpha j}\\
=&~\frac{2F^{11}}{v_{11}}v_{11\alpha}^2 +\frac{2}{v_{11}v_{22}} (F^{11}v_{11} +F^{22}v_{22} -v^2v_{11}v_{22}) v_{12\alpha}^2 +\frac{2F^{22}}{v_{22}} v_{22\alpha}^2 +4\frac{F^{12}}{v_{11}}v_{11\alpha} v_{12\alpha}\\
&+2v^2v_{11\alpha}v_{22\alpha} +4\frac{F^{12}}{v_{22}} v_{12\alpha}v_{22\alpha} +2\bigg(\frac{F^{11}}{v} +v(v_{22}+v_{44})\bigg)  v_{11\alpha} v_\alpha\\
&+4\frac{F^{12}}{v} v_{12\alpha}v_\alpha +2\bigg(\frac{F^{22}}{v} +v(v_{11}+v_{33})\bigg) v_{22\alpha}v_\alpha  +2(v_{11}+v_{33}) (v_{22}+v_{44}) v_\alpha^2.
\end{aligned}
\end{equation}
By \eqref{Estimate-v_{33},v_{44}'}, equation \eqref{Relation-1st condition-2} is equivalent to $$\sum_{i, j\in G}F^{ij}v_{ij\alpha} +F^vv_\alpha =O\bigg(\phi +\sum_{\alpha,\beta\in B}|\nabla v_{\alpha\beta}|\bigg), $$
namely, 
\begin{equation}\label{Relation-v_alpha-2}
v_\alpha=-\frac{F^{11}}{F^v}v_{11\alpha} -2\hspace{1pt}\frac{F^{12}}{F^v} v_{12\alpha} -\frac{F^{22}}{F^v}v_{22\alpha} +O\bigg(\phi+\sum_{\alpha, \beta\in B}|\nabla v_{\alpha\beta}|\bigg).
\end{equation}
In equation \eqref{Relation-v_alpha-2}, we used the fact that $F^v\neq0$. This can be derived from the identity 
\begin{align}\label{Relation-Fv-not0}
F^{11}v_{11}+F^{22}v_{22} +F^{33}v_{33} +F^{44}v_{44}=vF^v. 
\end{align}
Putting \eqref{Relation-v_alpha-2} into \eqref{Equality-F3-2}, we have
\begin{align}\label{Equality-F3-2a}
\begin{split}
&\sum_{i, j, k, l\in G}F^{ij, kl}v_{ij\alpha}v_{kl\alpha} +2\sum_{i, j\in G}F^{ij, v} v_{ij\alpha}v_\alpha +F^{v, v}v_\alpha^2 +2\sum_{i, j, \beta\in G}\frac{1}{v_{\beta\beta}}F^{ij}v_{\alpha\beta i}v_{\beta\alpha j}\\
=&\bigg(\frac{2F^{11}}{v_{11}} -2\bigg(\frac{F^{11}}{v}+vv_{22}\bigg) \frac{F^{11}}{F^v} +2v_{11}v_{22}\bigg(\frac{F^{11}}{F^v} \bigg)^2\bigg)v_{11\alpha}^2\\
&+\bigg(\frac{2}{v_{11}v_{22}}(F^{11}v_{11}+F^{22}v_{22}-v^2v_{11}v_{22}) -\frac{8(F^{12})^2}{vF^v} +8v_{11}v_{22}\bigg(\frac{F^{12}}{F^v}\bigg)^2 \bigg) v_{12\alpha}^2\\
&+\bigg(\frac{2F^{22}}{v_{22}} -2\bigg(\frac{F^{22}}{v}+vv_{11}\bigg) \frac{F^{22}}{F^v} +2v_{11}v_{22}\bigg(\frac{F^{22}}{F^v}\bigg)^2 \bigg) v_{22\alpha}^2\\
&+2\bigg(\frac{2F^{12}}{v_{11}}-2\bigg(\frac{F^{11}}{v}+vv_{22}\bigg)\frac{F^{12}}{F^v}-\frac{2F^{11}F^{12}}{vF^v}+4v_{11}v_{22}\hspace{1pt}\frac{F^{11}F^{12}}{(F^v)^2}\bigg)v_{11\alpha}v_{12\alpha}\\
&+2\bigg(v^2-\bigg(\frac{F^{11}}{v}+vv_{22}\bigg)\frac{F^{22}}{F^v} -\bigg(\frac{F^{22}}{v}+vv_{11}\bigg)\frac{F^{11}}{F^v}+2v_{11}v_{22} \hspace{1pt}\frac{F^{11}F^{22}}{(F^v)^2}\bigg)v_{11\alpha}v_{22\alpha}\\
&+2\bigg(\frac{2F^{12}}{v_{22}}-2\bigg(\frac{F^{22}}{v}+vv_{11}\bigg)\frac{F^{12}}{F^v}-\frac{2F^{12}F^{22}}{vF^v}+4v_{11}v_{22}\hspace{1pt}\frac{F^{12}F^{22}}{(F^v)^2}\bigg)v_{12\alpha}v_{22\alpha}\\
&+O\bigg(\phi+\sum_{\alpha,\beta\in B}|\nabla v_{\alpha\beta}|\bigg).
\end{split}
\end{align}
 
Recalling our equation $F(\nabla^2v, \nabla v, v)=R(x)$ in \eqref{Equation-minimal rank is 2}, and utilizing the estimates $v_{33}=O(\phi)$ and $v_{44}=O(\phi)$, we obtain 
\begin{align*}
&v^2v_{11}v_{22}+v(v_2^2+v_4^2)v_{11}+v(v_1^2+v_3^2)v_{22}=1+O(\phi), 
\intertext{and}
&F^{11}=v^2v_{22}+v(v_2^2+v_4^2)+O(\phi), \quad F^{22}=v^2v_{11}+v(v_1^2+v_3^2)+O(\phi), \\
&~F^v=2vv_{11}v_{22}+(v_2^2+v_4^2)v_{11}+(v_1^2+v_3^2)v_{22}+O(\phi).
\end{align*}
Clearly, $F^{11}$, $F^{22}$ and $F^v$ satisfy the following relations
\begin{align}
vF^v=&~F^{11}v_{11}+F^{22}v_{22}+O(\phi), \label{Relation-vF^v1-2} \\
vF^v=&~v^2v_{11}v_{22}+1+O(\phi). \label{Relation-vF^v2-2}
\end{align}
Let us simplify the six coefficients of the third-order derivatives on the right-hand side of \eqref{Equality-F3-2a}. Throughout these calculations, we will frequently use \eqref{Relation-vF^v1-2} and \eqref{Relation-vF^v2-2}.
\begin{itemize}
\item The coefficient of $v_{11\alpha}^2$ is 
\begin{align*}
&\frac{2F^{11}}{v_{11}}-2\bigg(\frac{F^{11}}{v}+vv_{22}\bigg)\frac{F^{11}}{F^v}+2v_{11}v_{22}\bigg(\frac{F^{11}}{F^v}\bigg)^2\\
=&~\frac{2F^{11}}{v_{11}(vF^v)^2}((vF^v)^2-F^{11}v_{11}\cdot vF^v -v^2v_{11}v_{22}\cdot vF^v +v^2v_{11}v_{22}\cdot F^{11}v_{11})\\
=&~\frac{2}{(vF^v)^2}\frac{F^{11}}{v_{11}}(F^{22}v_{22}\cdot vF^v -v^2v_{11}v_{22}\cdot F^{22}v_{22}+O(\phi))\\
=&~\frac{2}{(vF^v)^2}\frac{v_{22}}{v_{11}}F^{11}F^{22}+O(\phi).
\end{align*}

\item The coefficient of $v_{12\alpha}^2$ is
\begin{align*}
&\frac{2}{v_{11}v_{22}}(F^{11}v_{11}+F^{22}v_{22}-v^2v_{11}v_{22})-\frac{8(F^{12})^2}{vF^v}+8v_{11}v_{22}\bigg(\frac{F^{12}}{F^v}\bigg)^2\\
=&~\frac{2}{v_{11}v_{22}}-\frac{8(F^{12})^2}{vF^v}+8v_{11}v_{22}\bigg(\frac{F^{12}}{F^v}\bigg)^2+O(\phi)\\
=&~\frac{2}{(vF^v)^2}\frac{1}{v_{11}v_{22}}((vF^v)^2-4v_{11}v_{22}(F^{12})^2\cdot vF^v+4v_{11}v_{22}(F^{12})^2\cdot v^2v_{11}v_{22})+O(\phi)\\
=&~\frac{2}{(vF^v)^2}\frac{1}{v_{11}v_{22}}((vF^v)^2-4v_{11}v_{22}(F^{12})^2)+O(\phi).
\end{align*}

\item The coefficient of $v_{22\alpha}^2$ is
\begin{align*}
&\frac{2F^{22}}{v_{22}}-2\bigg(\frac{F^{22}}{v}+vv_{11}\bigg)\frac{F^{22}}{F^v}+2v_{11}v_{22}\bigg(\frac{F^{22}}{F^v}\bigg)^2\\
=&~\frac{2}{(vF^v)^2}\frac{v_{11}}{v_{22}}F^{11}F^{22}+O(\phi).
\end{align*}

\item The coefficient of $2v_{11\alpha}v_{12\alpha}$ is
\begin{align*}
&\frac{2F^{12}}{v_{11}}-2\bigg(\frac{F^{11}}{v}+vv_{22}\bigg)\frac{F^{12}}{F^v}-\frac{2F^{11}F^{12}}{vF^v}+4v_{11}v_{22}\hspace{1pt}\frac{F^{11}F^{12}}{(F^v)^2}\\
=&~\frac{2}{(vF^v)^2}\frac{F^{12}}{v_{11}}((vF^v)^2-2F^{11}v_{11}\cdot vF^v -v^2v_{11}v_{22}\cdot vF^v+2F^{11}v_{11}\cdot v^2v_{11}v_{22})\\
=&~\frac{2}{(vF^v)^2}\frac{F^{12}}{v_{11}}(vF^v-2F^{11}v_{11})+O(\phi).
\end{align*}

\item The coefficient of $2v_{11\alpha}v_{22\alpha}$ is
\begin{align*}
&v^2-\bigg(\frac{F^{11}}{v}+vv_{22}\bigg)\frac{F^{22}}{F^v}-\bigg(\frac{F^{22}}{v}+vv_{11}\bigg)\frac{F^{11}}{F^v}+2v_{11}v_{22}\hspace{1pt}\frac{F^{11}F^{22}}{(F^v)^2}\\
=&~v^2-\frac{2F^{11}F^{22}}{vF^v}-\frac{v}{F^v}(F^{11}v_{11}+F^{22}v_{22})+2v_{11}v_{22}\hspace{1pt}\frac{F^{11}F^{22}}{(F^v)^2}\\
=&~\frac{2}{(vF^v)^2}F^{11}F^{22}(-vF^v+v^2v_{11}v_{22})+O(\phi)\\
=&-\frac{2}{(vF^v)^2}F^{11}F^{22}+O(\phi).
\end{align*}

\item The coefficient of $2v_{12\alpha}v_{22\alpha}$ is
\begin{align*}
&\frac{2F^{12}}{v_{22}}-2\bigg(\frac{F^{22}}{v}+vv_{11}\bigg)\frac{F^{12}}{F^v}-\frac{2F^{12}F^{22}}{vF^v}+4v_{11}v_{22}\hspace{1pt}\frac{F^{12}F^{22}}{(F^v)^2}\\
=&~\frac{2}{(vF^v)^2}\frac{F^{12}}{v_{22}}(vF^v-2v_{22}F^{22})+O(\phi).
\end{align*}
\end{itemize}

Therefore, formula \eqref{Equality-F3-2a} is reduced to 
\begin{align*}
&\frac{(vF^v)^2}{2}\bigg(\sum_{i, j, k, l\in G}F^{ij, kl}v_{ij\alpha}v_{kl\alpha} +2\sum_{i, j\in G}F^{ij, v} v_{ij\alpha}v_\alpha +F^{v, v}v_\alpha^2+2\sum_{i, j, \beta\in G}\frac{1}{v_{\beta\beta}} F^{ij}v_{\alpha\beta i}v_{\beta\alpha j}\bigg)\\
=&~\frac{v_{22}}{v_{11}}F^{11}F^{22}v_{11\alpha}^2+\frac{1}{v_{11}v_{22}}((vF^v)^2-4v_{11}v_{22}(F^{12})^2)v_{12\alpha}^2+\frac{v_{11}}{v_{22}}F^{11}F^{22}v_{22\alpha}^2\\
&+2\hspace{1pt}\frac{F^{12}}{v_{11}}(vF^v-2F^{11}v_{11})v_{11\alpha}v_{12\alpha}-2F^{11}F^{22}v_{11\alpha}v_{22\alpha}+2\hspace{1pt}\frac{F^{12}}{v_{22}}(vF^v-2F^{22}v_{22})v_{12\alpha}v_{22\alpha}\\
&+O\bigg(\phi+\sum_{\alpha,\beta\in B}|\nabla v_{\alpha\beta}|\bigg).
\end{align*}
Let us define a quadratic form in three variables $$f(X)=X^TMX, $$ where $X=(v_{11\alpha}, v_{12\alpha}, v_{22\alpha})^T$ and the symmetric matrix
$$M=\begin{pmatrix}
\frac{v_{22}}{v_{11}}F^{11}F^{22} & \frac{F^{12}}{v_{11}}(vF^v-2F^{11}v_{11}) & -F^{11}F^{22}\\[7pt]
\frac{F^{12}}{v_{11}}(vF^v-2F^{11}v_{11}) & \frac{1}{v_{11}v_{22}}((vF^v)^2-4v_{11}v_{22}(F^{12})^2) & \frac{F^{12}}{v_{22}}(vF^v-2F^{22}v_{22})\\[7pt]
-F^{11}F^{22} & \frac{F^{12}}{v_{22}}(vF^v-2F^{22}v_{22}) & \frac{v_{11}}{v_{22}}F^{11}F^{22}
\end{pmatrix}. $$ Then, it suffices to establish the semi-positivity of $f(X)$ in the sense that $$f(X)\geq -C\phi. $$ This is equivalent to saying that $A$ is semi-positive modulo $\phi$.

By Lemma \ref{Lemma-linear algebra}, we need to compute the leading principal minors of the matrix $M$, denoted by $P_1(M)$, $P_2(M)$ and $P_3(M)$.

$\bullet$ $P_1(M)$. Clearly, $$P_1(M)=\frac{v_{22}}{v_{11}}F^{11}F^{22}. $$

$\bullet$ $P_2(M)$. By \eqref{Relation-vF^v1-2} and \eqref{Relation-vF^v2-2}, we have
\begin{align*}
P_2(M)=&~\frac{F^{11}F^{22}}{v_{11}^2}((vF^v)^2-4v_{11}v_{22}(F^{12})^2)-\frac{(F^{12})^2}{v_{11}^2}(vF^v-2F^{11}v_{11})^2\\
=&~\frac{(vF^v)^2}{v_{11}^2}(F^{11}F^{22}-(F^{12})^2)+O(\phi).
\end{align*}

$\bullet$ $P_3(M)$. Observe that the first and third columns of $M$ are proportional, which immediately implies that $$P_3(M)=\det(M)=0. $$

So we conclude that $A$ is semi-positive modulo $\phi$ by Lemma \ref{Lemma-linear algebra}, completing the proof of \textbf{Claim 1}. \qed

An application of Theorem \ref{Theorem-rank2} shows that if $\nabla^2v(x)$ attains its minimal rank $2$ at some point $x_0\in\Omega$, then $\nabla^2v$ must have rank $2$ everywhere in $\Omega$. The detailed proof of this consequence will be given in Corollary \ref{Corollary-CRT} of Section \ref{Section 4}.

\section{Constant rank theorem: the case of minimal rank $3$}\label{Section 4}

In this section, we study the case where the real Hessian $\nabla^2v$ attains its minimal rank $3$ in $\Omega$, equivalently, $K$ has rank $1$. We define $$\phi(x)=\det(K(x)).$$ Analogous to the previous section, our goal is to prove a differential inequality for $\phi$. Finally, combining Theorem \ref{Theorem-rank2} with the strong maximum principle, we complete the proof of the constant rank theorem.

\begin{theorem}\label{Theorem-rank3}
Let $\Omega\subset\mathbb{R}^4$ be a domain and let $v$ be a convex solution to equation \eqref{Equ-real equation(v)}. If the Hessian $\nabla^2v(x)$ attains its minimal rank $3$ at some point $x_0\in\Omega$, then there exist a neighborhood $\mathcal{O}$ of $x_0$ and a positive constant $C$ (independent of $\phi$) such that
\begin{equation}\label{Inequ-lemma(minimal rank=3)}
\sum_{i, j}F^{ij}\phi_{ij}\leq C(\phi+|\nabla\phi|) \quad \mathrm{in}\ \mathcal{O}.
\end{equation}
\end{theorem}

\begin{proof}
Following the notations in \cite{Caffarelli-Friedman1985, Korevaar-Lewis1987}, for two functions $h(x)$ and $k(x)$ defined in $\mathcal{O}$, we say that $h(x)\lesssim k(x)$ if there exists a positive constant $C_1$ such that $$h(x)-k(x)\leq C_1(\phi(x)+|\nabla\phi|(x)). $$
We also write $h(x)\sim k(x)$ if $h(x)\lesssim k(x)$ and $k(x)\lesssim h(x)$. Next, we write $h\lesssim k$ if the above inequality holds for any $x\in\mathcal{O}$, with constant $C_1$ independent of $x$. Finally, $h\sim k$ means $h\lesssim k$ and $k\lesssim h$. In a similar way, we will also use the notation $h\gtrsim k$ which means $-h\lesssim -k$, namely, $$-h+k\leq C_2(\phi+|\nabla\phi|) \quad \mathrm{in}\ \mathcal{O}, $$
or equivalently, $$h-k\geq -C_2(\phi+|\nabla\phi|) \quad \mathrm{in}\ \mathcal{O}. $$ Thus, $h\sim k$ can also be interpreted as $h\lesssim k$ and $h\gtrsim k$, meaning that there exist positive constants $C_3$ and $C_4$ such that $$-C_3(\phi+|\nabla\phi|) \leq h-k \leq C_4(\phi+|\nabla\phi|) \quad \mathrm{in} \ \mathcal{O}. $$ We can also write this as $$h=k+O(\phi+|\nabla\phi|). $$
 
For each $x\in\mathcal{O}$ fixed, we choose coordinates introduced in Section \ref{Section-auxiliary function} so that $\nabla^2v(x)$ is diagonal. In the following, all calculations will be carried out at the fixed point $x$. When we use the relation $\lesssim$ or $\gtrsim$, all constants involved are under controlled. In this subsection, we set $G=\{1, 2, 3\}$ and $B=\{4\}$. From \eqref{Equality-det K and det D^2v}, the auxiliary function $\phi$ takes the form $$\phi=(\det(A))^{-2}\sigma_4(\nabla^2v),$$ with $A$ as defined in \eqref{Matrix-A}. Then we have
$$0\sim\phi=(\det(A))^{-2}\sigma_3(G)v_{44}, $$ and
\begin{equation}\label{Relation-v_44-3}
v_{44}\sim 0.
\end{equation}
Taking the first derivatives of $\phi$, by Lemma \ref{Lemma-sigma_k}, we get
\begin{align*}
0\sim\phi_i
=&~\big((\det(A))^{-2}\big)_i\sigma_4+(\det(A))^{-2}\sum_{p,q}\frac{\partial\sigma_4}{\partial v_{pq}}v_{pqi}\\
\sim&~(\det(A))^{-2}\sum_p\sigma_3(\lambda|p)v_{ppi}\sim(\det(A))^{-2}\sigma_3(G)v_{44i},\quad i=1,2,3,4.
\end{align*} 
namely, 
\begin{equation}\label{Relation-v_44i-3}
v_{44i}\sim 0, \quad i=1, 2, 3, 4.
\end{equation}
Using relations \eqref{Relation-v_44-3} and \eqref{Relation-v_44i-3}, by Lemma \ref{Lemma-sigma_k} again, we compute the second derivatives of $\phi$ as
\begin{align*}
\phi_{ij}
=&~\big((\det(A))^{-2}\big)_{ij}\sigma_4+\big((\det(A))^{-2}\big)_i\sum_{p,q}\frac{\partial\sigma_4}{\partial v_{pq}}v_{pqj}+\big((\det(A))^{-2}\big)_j\sum_{p,q}\frac{\partial\sigma_4}{\partial v_{pq}}v_{pqi}\\
&+(\det(A))^{-2}\sum_{p,q}\frac{\partial\sigma_4}{\partial v_{pq}}v_{pqij}+(\det(A))^{-2}\sum_{p,q,r,s}\frac{\partial^2\sigma_4}{\partial v_{pq}\partial v_{rs}}v_{pqi}v_{rsj}\\
=&~\big((\det(A))^{-2}\big)_{ij}\sigma_3(G)v_{44}+\big((\det(A))^{-2}\big)_i\sum_p\sigma_3(\lambda|p)v_{ppj}+\big((\det(A))^{-2}\big)_j\sum_p\sigma_3(\lambda|p)v_{ppi}\\
&+(\det(A))^{-2}\sum_p\sigma_3(\lambda|p)v_{ppij}\\
&+(\det(A))^{-2}\sum_{p\neq r}\sigma_2(\lambda|pr) v_{ppi}v_{rrj} -(\det(A))^{-2}\sum_{p\neq r}\sigma_2(\lambda|pr)v_{pri}v_{rpj}\\
\sim&~(\det(A))^{-2}\sigma_3(G)v_{44ij}-2(\det(A))^{-2}\sigma_3(G)\sum_{p\in G}\frac{1}{v_{pp}}v_{4pi}v_{4pj}
\end{align*}
It follows that
\begin{equation}\label{Equality-F1-3}
-\frac{(\det(A))^2}{\sigma_3(G)}\sum_{i,j}F^{ij}\phi_{ij}\sim-\sum_{i,j}F^{ij}v_{44ij}+2\sum_{p\in G}\sum_{i,j}\frac{1}{v_{pp}}F^{ij}v_{4pi}v_{4pj}.
\end{equation}

To handle the fourth-order derivative terms in \eqref{Equality-F1-3}, we first take the derivative of the equation $F(\nabla^2v, \nabla v, v)=0$ with respect to the variable $x_4$, resulting in
\begin{equation}\label{Equality-1st order derivative of F-3}
0=\sum_{i, j}F^{ij}v_{ij4}+\sum_k F^{v_k}v_{k4}+F^vv_4.
\end{equation}
Then we take the second derivative, which gives
\begin{equation}\label{Equality-2nd order derivative of F-3}
\begin{aligned}
0=&\sum_{i, j, k, l}F^{ij, kl}v_{ij4}v_{kl4} +\sum_{i, j, k}F^{ij, v_k}v_{ij4}v_{k4}+\sum_{i, j}F^{ij, v}v_{ij4}v_4 +\sum_{i, j}F^{ij}v_{ij44}\\
&+\sum_{k, i, j}F^{v_k, ij}v_{k4}v_{ij4} +\sum_{k, l}F^{v_k, v_l}v_{k4}v_{l4} +\sum_k F^{v_k, v}v_{k4}v_4 +\sum_k F^{v_k}v_{k44}\\
&+\sum_{i, j}F^{v, ij}v_4v_{ij4}+\sum_k F^{v, v_k}v_4v_{k4} +F^{v, v}v_4^2+F^vv_{44}.
\end{aligned}
\end{equation}
By \eqref{Relation-v_44-3} and \eqref{Relation-v_44i-3}, equations \eqref{Equality-1st order derivative of F-3} and \eqref{Equality-2nd order derivative of F-3} can be simplified as
\begin{equation}\label{Relation-1st condition-3}
0\sim \sum_{i, j\in G}F^{ij}v_{ij4}+F^vv_4
\end{equation}
and
\begin{equation}\label{Relation-2nd condition-3}
-\sum_{i, j}F^{ij}v_{ij44}\sim \sum_{i, j, k, l\in G}F^{ij, kl}v_{ij4}v_{kl4} +2\sum_{i, j\in G}F^{ij, v}v_{ij4}v_4 +F^{v, v}v_4^2.
\end{equation}
Substituting \eqref{Relation-2nd condition-3} into \eqref{Equality-F1-3}, we have
\begin{equation}\label{Equality-F2-3}
\begin{aligned}
-\frac{(\det(A))^2}{\sigma_3(G)}\sum_{i, j}F^{ij}\phi_{ij} \sim&\sum_{i, j, k, l\in G}F^{ij, kl}v_{ij4}v_{kl4} +2\sum_{i, j\in G}F^{ij, v}v_{ij4}v_4 +F^{v, v}v_4^2\\
&+2\sum_{p\in G}\sum_{i, j}\frac{1}{v_{pp}}F^{ij}v_{4pi}v_{4pj}.
\end{aligned}
\end{equation}
 
Next, we carefully analyze each term on the right-hand side of \eqref{Equality-F2-3}. Recall that
\begin{align*}
F(\nabla^2v, \nabla v, v) =~&v^2\left[(v_{11}+v_{33})(v_{22}+v_{44}) -(v_{12}+v_{34})(v_{21}+v_{43})-(v_{14}-v_{32})(v_{41}-v_{23})\right]\\
&+v\left[(v_2^2+v_4^2)(v_{11}+v_{33})+(v_1^2+v_3^2)(v_{22}+v_{44})\right. \\
&\hspace{8mm}-(v_1v_2+v_3v_4)(v_{12}+v_{34}+v_{21}+v_{43})\\
&\hspace{8mm}\left.-(v_1v_4-v_2v_3)(v_{14}-v_{32}+v_{41}-v_{23})\right]-1
\end{align*}
Direct calculation implies
\begin{align*}
F^{11, 22}&=F^{22, 11}=F^{22, 33}=F^{33, 22}=v^2, \\
F^{12, 21}&=F^{21, 12}=F^{23, 32}=F^{32, 23}=-v^2.
\end{align*}
By collecting the terms in equation \eqref{Equality-F2-3}, we obtain
\begin{equation}\label{Equality-F3-3}
\begin{aligned}
-\frac{(\det(A))^2}{\sigma_3(G)}\sum_{i,j}F^{ij}\phi_{ij}\sim&~\frac{2F^{11}}{v_{11}}v_{114}^2+\frac{2F^{22}}{v_{22}}v_{224}^2+\frac{2F^{33}}{v_{33}}v_{334}^2+\bigg(\frac{2F^{33}}{v_{11}}+\frac{2F^{11}}{v_{33}}\bigg)v_{134}^2\\
&+\bigg(-2v^2+\frac{2F^{22}}{v_{11}}+\frac{2F^{11}}{v_{22}}\bigg)v_{124}^2+\bigg(-2v^2+\frac{2F^{33}}{v_{22}}+\frac{2F^{22}}{v_{33}}\bigg)v_{234}^2\\
&+2v^2v_{114}v_{224}+\frac{4F^{12}}{v_{11}}v_{114}v_{124}+2v^2v_{224}v_{334}+\frac{4F^{12}}{v_{22}}v_{224}v_{124}\\
&+\frac{4F^{23}}{v_{22}}v_{224}v_{234}+\frac{4F^{23}}{v_{33}}v_{334}v_{234}+\frac{4F^{23}}{v_{11}}v_{134}v_{124}+\frac{4F^{12}}{v_{33}}v_{134}v_{234}\\
&+2F^{11, v}v_{114}v_4+2F^{22, v}v_{224}v_4+2F^{33, v}v_{334}v_4+4F^{12, v}v_{124}v_4\\
&+4F^{23, v}v_{234}v_4+F^{v, v}v_4^2.
\end{aligned}
\end{equation}
It is easy to see that
\begin{equation}\label{Relation-vF^v-3}
vF^v=v^2(v_{11}+v_{33})(v_{22}+v_{44})+1.
\end{equation}
This implies $F^v\neq 0$. Since $F^{13}=0$, by \eqref{Relation-1st condition-3}, we have
\begin{equation}\label{Relation-v_4-3}
v_4\sim-\frac{F^{11}}{F^v}v_{114}-\frac{F^{22}}{F^v}v_{224}-\frac{F^{33}}{F^v}v_{334}-\frac{2F^{12}}{F^v}v_{124}-\frac{2F^{23}}{F^v}v_{234}.
\end{equation}
Putting \eqref{Relation-v_4-3} into \eqref{Equality-F3-3}, we obtain a quadratic form in terms of $v_{114}$, $v_{224}$, $v_{334}$, $v_{134}$, $v_{124}$ and $v_{234}$, given by
\begin{align*}
-\frac{(\det(A))^2}{\sigma_3(G)}\sum_{i,j}F^{ij}\phi_{ij}
\sim&\left(\frac{2F^{11}}{v_{11}}-\frac{2F^{11}F^{11,v}}{F^v}+\left(\frac{F^{11}}{F^v}\right)^2F^{v,v}\right)v_{114}^2\\
&+\left(\frac{2F^{22}}{v_{22}}-\frac{2F^{22}F^{22,v}}{F^v}+\left(\frac{F^{22}}{F^v}\right)^2F^{v,v}\right)v_{224}^2\\
&+\left(\frac{2F^{33}}{v_{33}}-\frac{2F^{33}F^{33,v}}{F^v}+\left(\frac{F^{33}}{F^v}\right)^2F^{v,v}\right)v_{334}^2\\
&+\left(\frac{2F^{33}}{v_{11}}+\frac{2F^{11}}{v_{33}}\right)v_{134}^2\\
&+\left(-2v^2+\frac{2F^{22}}{v_{11}}+\frac{2F^{11}}{v_{22}}-\frac{8F^{12}F^{12,v}}{F^v}+4\left(\frac{F^{12}}{F^v}\right)^2F^{v,v}\right)v_{124}^2\\
&+\left(-2v^2+\frac{2F^{33}}{v_{22}}+\frac{2F^{22}}{v_{33}}-\frac{8F^{23}F^{23,v}}{F^v}+4\left(\frac{F^{23}}{F^v}\right)^2F^{v,v}\right)v_{234}^2\\
&+2\left(v^2-\frac{F^{22}F^{11,v}}{F^v}-\frac{F^{11}F^{22,v}}{F^v}+\frac{F^{11}F^{22}}{(F^v)^2}F^{v,v}\right)v_{114}v_{224}\\
&+2\left(-\frac{F^{33}F^{11,v}}{F^v}-\frac{F^{11}F^{33,v}}{F^v}+\frac{F^{11}F^{33}}{(F^v)^2}F^{v,v}\right)v_{114}v_{334}\\
&+2\left(\frac{2F^{12}}{v_{11}}-\frac{2F^{12}F^{11,v}}{F^v}-\frac{2F^{11}F^{12,v}}{F^v}+\frac{2F^{11}F^{12}}{(F^v)^2}F^{v,v}\right)v_{114}v_{124}\\
&+2\left(-\frac{2F^{23}F^{11,v}}{F^v}-\frac{2F^{11}F^{23,v}}{F^v}+\frac{2F^{11}F^{23}}{(F^v)^2}F^{v,v}\right)v_{114}v_{234}\\
&+2\left(v^2-\frac{F^{33}F^{22,v}}{F^v}-\frac{F^{22}F^{33,v}}{F^v}+\frac{F^{22}F^{33}}{(F^v)^2}F^{v,v}\right)v_{224}v_{334}\\
&+2\left(\frac{2F^{12}}{v_{22}}-\frac{2F^{12}F^{22,v}}{F^v}-\frac{2F^{22}F^{12,v}}{F^v}+\frac{2F^{12}F^{22}}{(F^v)^2}F^{v,v}\right)v_{224}v_{124}\\
&+2\left(\frac{2F^{23}}{v_{22}}-\frac{2F^{23}F^{22,v}}{F^v}-\frac{2F^{22}F^{23,v}}{F^v}+\frac{2F^{22}F^{23}}{(F^v)^2}F^{v,v}\right)v_{224}v_{234}\\
&+2\left(-\frac{2F^{12}F^{33,v}}{F^v}-\frac{2F^{33}F^{12,v}}{F^v}+\frac{2F^{12}F^{33}}{(F^v)^2}F^{v,v}\right)v_{334}v_{124}\\
&+2\left(\frac{2F^{23}}{v_{33}}-\frac{2F^{23}F^{33,v}}{F^v}-\frac{2F^{33}F^{23,v}}{F^v}+\frac{2F^{23}F^{33}}{(F^v)^2}F^{v,v}\right)v_{334}v_{234}\\
&+2\left(\frac{2F^{23}}{v_{11}}\right)v_{134}v_{124}+2\left(\frac{2F^{12}}{v_{33}}\right)v_{134}v_{234}\\
&+2\left(-\frac{4F^{23}F^{12,v}}{F^v}-\frac{4F^{12}F^{23,v}}{F^v}+\frac{4F^{12}F^{23}}{(F^v)^2}F^{v,v}\right)v_{124}v_{234}.
\end{align*}
It suffices to verify that the above quadratic form in the variables $v_{114}$, $v_{224}$, $v_{334}$, $v_{134}$, $v_{124}$ and $v_{234}$ is semi-positive definite modulo $\phi$. To achieve this, we must arrange each coefficients in a sophisticated manner.

Before arranging these coefficients, we recall the equation $F(\nabla^2v, \nabla v, v)=0$. It is equivalent to
\begin{equation}\label{Relation-Equation-3}
0\sim v^2(v_{11}+v_{33})v_{22}+v(v_2^2+v_4^2)(v_{11}+v_{33}) +v(v_1^2+v_3^2)v_{22}-1.
\end{equation}
Clearly, the matrix $(F^{ij})$ is given by
\begin{align}\label{Relation-Fij-a}
(F^{ij})=\begin{pmatrix}
F^{11} & F^{12} & 0 & -F^{23}\\[4pt]
F^{12} & F^{22} & F^{23} & 0\\[4pt]
0 & F^{23} & F^{11} & F^{12}\\[4pt]
-F^{23} & 0 & F^{12} & F^{22}
\end{pmatrix},
\end{align}
where
\begin{align}\label{Relation-Fij-b}
\begin{split}
&F^{11}\sim v^2v_{22}+v(v_2^2+v_4^2), \hspace{1.7cm} F^{12}=-v(v_1v_2+v_3v_4), \\
&F^{22}= v^2(v_{11}+v_{33})+v(v_1^2+v_3^2), \quad F^{23}= v(v_1v_4-v_2v_3).
\end{split}
\end{align}
Since $v_{44}\sim0$, by \eqref{Relation-vF^v-3}, we have
\begin{equation}\label{Relation-vF^v1-3}
vF^v\sim v^2(v_{11}+v_{33})v_{22}+1.
\end{equation}
Combining this with \eqref{Relation-Fij-a}--\eqref{Relation-Fij-b}, we also have
\begin{equation}\label{Relation-vF^v2-3}
vF^v\sim F^{11}v_{11}+F^{22}v_{22}+F^{33}v_{33}.
\end{equation}
Furthermore, a simple computation gives
\begin{align}\label{Relation-Fijv}
\begin{split}
&F^{11, v}\sim 2vv_{22}+(v_2^2+v_4^2)\sim F^{33, v}, \\
&F^{22, v}\sim 2v(v_{11}+v_{33})+(v_1^2+v_3^2), \\
&F^{12, v}\sim -(v_1v_2+v_3v_4), \quad F^{23, v}\sim v_1v_4-v_2v_3, \\
&F^{v, v}\sim 2(v_{11}+v_{33})v_{22}.
\end{split}
\end{align}
We now begin to simplify the coefficients of all terms in the final quadratic form. Throughout this process, we will repeatedly use \eqref{Relation-Fij-a}--\eqref{Relation-Fijv}.

$\bullet$ The coefficient of $v_{114}^2$ is 
\begin{align*}
&\frac{2F^{11}}{v_{11}}-\frac{2F^{11}F^{11,v}}{F^v}+\bigg(\frac{F^{11}}{F^v}\bigg)^2F^{v,v}\\
=&~\frac{1}{(vF^v)^2}\frac{F^{11}}{v_{11}}\left(2(vF^v)^2-2v_{11}\cdot vF^{11,v}\cdot vF^v+F^{11}v_{11}\cdot v^2F^{v,v}\right)\\
\sim&~\frac{1}{(vF^v)^2}\frac{F^{11}}{v_{11}}\left(2(v^2(v_{11}+v_{33})v_{22}+1)\cdot vF^v-2v_{11}\cdot(F^{11}+v^2v_{22})\cdot vF^v\right.\\
&\hspace{2.3cm}\left.+F^{11}v_{11}\cdot2v^2(v_{11}+v_{33})v_{22}\right)\\
=&~\frac{2}{(vF^v)^2}\frac{F^{11}}{v_{11}}\left(v^2v_{22}v_{33}\cdot vF^v+vF^v-F^{11}v_{11}\cdot vF^v+F^{11}v_{11}\cdot v^2(v_{11}+v_{33})v_{22} \right)\\
\sim&~\frac{2}{(vF^v)^2}\frac{F^{11}}{v_{11}}\left(v^2v_{22}v_{33}\cdot vF^v+vF^v-F^{11}v_{11}\right)\\
\sim&~\frac{2}{(vF^v)^2}\frac{F^{11}}{v_{11}}\left(v^2v_{22}v_{33}\cdot vF^v+F^{22}v_{22}+F^{33}v_{33}\right).
\end{align*}

$\bullet$ The coefficient of $v_{224}^2$ is 
\begin{align*}
&\frac{2F^{22}}{v_{22}}-\frac{2F^{22}F^{22,v}}{F^v}+\bigg(\frac{F^{22}}{F^v}\bigg)^2F^{v,v}\\
=&~\frac{1}{(vF^v)^2}\frac{F^{22}}{v_{22}}\left(2(vF^v)^2-2v_{22}\cdot vF^{22,v}\cdot vF^v+F^{22}v_{22}\cdot v^2F^{v,v}\right)\\
\sim&~\frac{1}{(vF^v)^2}\frac{F^{22}}{v_{22}}\left(2(v^2(v_{11}+v_{33})v_{22}+1)\cdot vF^v-2v_{22}\cdot(F^{22}+v^2(v_{11}+v_{33}))\cdot vF^v\right.\\
&\hspace{2.3cm}\left.+F^{22}v_{22}\cdot2v^2(v_{11}+v_{33})v_{22}\right)\\
=&~\frac{2}{(vF^v)^2}\frac{F^{22}}{v_{22}}\left(vF^v-F^{22}v_{22}\cdot vF^v+F^{22}v_{22}\cdot v^2(v_{11}+v_{33})v_{22}\right)\\
\sim&~\frac{2}{(vF^v)^2}\frac{F^{22}}{v_{22}}(vF^v-F^{22}v_{22})\\
\sim&~\frac{2}{(vF^v)^2}\frac{F^{22}}{v_{22}}(F^{11}v_{11}+F^{33}v_{33}).
\end{align*}

$\bullet$ The coefficient of $v_{334}^2$ is 
\begin{align*}
&\frac{2F^{33}}{v_{33}}-\frac{2F^{33}F^{33,v}}{F^v}+\bigg(\frac{F^{33}}{F^v}\bigg)^2F^{v,v}\\
\sim&~\frac{2}{(vF^v)^2}\frac{F^{33}}{v_{33}}\left(v^2v_{11}v_{22}\cdot vF^v+F^{11}v_{11}+F^{22}v_{22}\right),
\end{align*}
where we exploited the symmetry between indices $1$ and $3$.

$\bullet$ The coefficient of $v_{134}^2$ is $$\frac{2F^{33}}{v_{11}} +\frac{2F^{11}}{v_{33}}. $$

$\bullet$ The coefficient of $v_{124}^2$ is
\begin{align*}
&-2v^2+\frac{2F^{22}}{v_{11}}+\frac{2F^{11}}{v_{22}}-\frac{8F^{12}F^{12,v}}{F^v}+4\bigg(\frac{F^{12}}{F^v}\bigg)^2F^{v,v}\\
=&~\frac{2}{(vF^v)^2}\bigg(\bigg(-v^2+\frac{F^{22}}{v_{11}}+\frac{F^{11}}{v_{22}}\bigg)\cdot(vF^v)^2-4F^{12}\cdot vF^{12,v}\cdot vF^v+2(F^{12})^2\cdot v^2F^{v,v}\bigg)\\
\sim&~\frac{2}{(vF^v)^2}\bigg(\frac{-v^2v_{11}v_{22}+F^{22}v_{22}+F^{11}v_{11}}{v_{11}v_{22}}\hspace{1pt}(vF^v)^2\\
&\hspace{1.7cm}-4(F^{12})^2\bigg(v^2(v_{11}+v_{33})v_{22}+1\bigg)+4(F^{12})^2\cdot v^2(v_{11}+v_{33})v_{22}\bigg)\\
\sim&~\frac{2}{(vF^v)^2}\bigg(\frac{v^2v_{22}v_{33}+1-F^{33}v_{33}}{v_{11}v_{22}}\hspace{1pt}(vF^v)^2-4(F^{12})^2\bigg).
\end{align*}

$\bullet$ The coefficient of $v_{234}^2$ is
\begin{align*}
&-2v^2+\frac{2F^{33}}{v_{22}}+\frac{2F^{22}}{v_{33}}-\frac{8F^{23}F^{23,v}}{F^v}+4\bigg(\frac{F^{23}}{F^v}\bigg)^2F^{v,v}\\
\sim&~\frac{2}{(vF^v)^2}\bigg(\frac{v^2v_{11}v_{22}+1-F^{11}v_{11}}{v_{22}v_{33}}\hspace{1pt}(vF^v)^2-4(F^{23})^2\bigg).
\end{align*}

$\bullet$ The coefficient of $2v_{114}v_{224}$ is
\begin{align*}
&v^2-\frac{F^{22}F^{11,v}}{F^v}-\frac{F^{11}F^{22,v}}{F^v}+\frac{F^{11}F^{22}}{(F^v)^2}F^{v,v}\\
=&~\frac{1}{(vF^v)^2}\left(v^2\cdot(vF^v)^2-F^{22}\cdot vF^{11,v}\cdot vF^v-F^{11}\cdot vF^{22,v}\cdot vF^v+F^{11}F^{22}\cdot v^2F^{v,v}\right)\\
\sim&~\frac{1}{(vF^v)^2}\left(v^2\cdot(vF^v)^2-F^{22}\cdot(F^{11}+v^2v_{22})\cdot vF^v\right)-F^{11}\cdot(F^{22}+v^2(v_{11}+v_{33}))\cdot vF^v\\
&\hspace{1.5cm}\left.+F^{11}F^{22}\cdot2v^2(v_{11}+v_{33})v_{22}\right)\\
=&~\frac{1}{(vF^v)^2}\left(v^2\cdot(vF^v)^2-v^2\cdot(F^{22}v_{22}+F^{11}v_{11}+F^{11}v_{33})\cdot vF^v-2F^{11}F^{22}\right.\\
&\hspace{1.5cm}\left.+2F^{11}F^{22}\cdot v^2(v_{11}+v_{33})v_{22}\right)\\
\sim&-\frac{2}{(vF^v)^2}\hspace{1pt}F^{11}F^{22}.
\end{align*}

$\bullet$ The coefficient of $2v_{114}v_{334}$ is
\begin{align*}
&-\frac{F^{33}F^{11,v}}{F^v}-\frac{F^{11}F^{33,v}}{F^v}+\frac{F^{11}F^{33}}{(F^v)^2}F^{v,v}\\
=&-\frac{1}{(vF^v)^2}\left(F^{33}\cdot vF^{11,v}\cdot vF^v+F^{11}\cdot vF^{33,v}\cdot vF^v-F^{11}F^{33}\cdot v^2F^{v,v}\right)\\
\sim&-\frac{1}{(vF^v)^2}\hspace{1pt}(F^{33}\cdot(F^{11}+v^2v_{22})\cdot vF^v+F^{11}\cdot(F^{33}+v^2v_{22})\cdot vF^v\\
&\hspace{2cm}-F^{11}F^{33}\cdot2v^2(v_{11}+v_{33})v_{22})\\
=&-\frac{1}{(vF^v)^2}\left(2F^{11}F^{33}-2F^{11}F^{33}\cdot v^2(v_{11}+v_{33})v_{22}+(F^{33}+F^{11})\cdot v^2v_{22}\cdot vF^v\right)\\
\sim&-\frac{2F^{11}}{(vF^v)^2}\left(F^{11}+v^2v_{22}\cdot vF^v\right).
\end{align*}

$\bullet$ The coefficient of $2v_{114}v_{124}$ is
\begin{align*}
&\frac{2F^{12}}{v_{11}}-\frac{2F^{12}F^{11,v}}{F^v}-\frac{2F^{11}F^{12,v}}{F^v}+\frac{2F^{11}F^{12}}{(F^v)^2}F^{v,v}\\
=&~\frac{2}{(vF^v)^2}\bigg(\frac{F^{12}}{v_{11}}(vF^v)^2-F^{12}\cdot vF^{11,v}\cdot vF^v-F^{11}\cdot vF^{12,v}\cdot vF^v +F^{11}F^{12}\cdot v^2F^{v,v}\bigg)\\
\sim&~\frac{2}{(vF^v)^2}\bigg(\frac{F^{12}}{v_{11}}\left(v^2(v_{11}+v_{33})v_{22}+1\right)\cdot vF^v -F^{12}\cdot(F^{11}+v^2v_{22})\cdot vF^v\\
&\hspace{1.8cm}-F^{11}F^{12}\cdot vF^v+F^{11}F^{12}\cdot 2v^2(v_{11}+v_{33})v_{22}\bigg)\\
=&~\frac{2}{(vF^v)^2}\bigg(\frac{F^{12}}{v_{11}}(v^2v_{22}v_{33}+1)\cdot vF^v-2F^{11}F^{12}(vF^v-v^2(v_{11}+v_{33})v_{22})\bigg)\\
\sim&~\frac{2}{(vF^v)^2}\frac{F^{12}}{v_{11}}(v^2v_{22}v_{33}\cdot vF^v+vF^v-2F^{11}v_{11}).
\end{align*}

$\bullet$ The coefficient of $2v_{114}v_{234}$ is
\begin{align*}
&-\frac{2F^{23}F^{11,v}}{F^v}-\frac{2F^{11}F^{23,v}}{F^v}+\frac{2F^{11}F^{23}}{(F^v)^2}F^{v,v}\\
=&-\frac{2}{(vF^v)^2}(F^{23}\cdot vF^{11,v}\cdot vF^v+F^{11}\cdot vF^{23,v}\cdot vF^v-F^{11}F^{23}\cdot v^2F^{v,v})\\
\sim&-\frac{2}{(vF^v)^2}\left(F^{23}\cdot(F^{11}+v^2v_{22})\cdot vF^v +F^{11}F^{23}\cdot vF^v-F^{11}F^{23}\cdot 2v^2(v_{11}+v_{33})v_{22}\right)\\
=&-\frac{2}{(vF^v)^2}\left(F^{23}\cdot v^2v_{22}\cdot vF^v+2F^{11}F^{23}\cdot(vF^v-v^2(v_{11}+v_{33})v_{22})\right)\\
\sim&-\frac{2F^{23}}{(vF^v)^2}\left(v^2v_{22}\cdot vF^v+2F^{11}\right).
\end{align*}

$\bullet$ The coefficient of $2v_{224}v_{334}$ is $$v^2-\frac{F^{33}F^{22,v}}{F^v}-\frac{F^{22}F^{33,v}}{F^v}+\frac{F^{22}F^{33}}{(F^v)^2}F^{v,v}\sim-\frac{2}{(vF^v)^2}\hspace{1pt}F^{22}F^{33}. $$

$\bullet$ The coefficient of $2v_{224}v_{124}$ is
\begin{align*}
&\frac{2F^{12}}{v_{22}}-\frac{2F^{12}F^{22,v}}{F^v}-\frac{2F^{22}F^{12,v}}{F^v}+\frac{2F^{12}F^{22}}{(F^v)^2}F^{v,v}\\
=&~\frac{2}{(vF^v)^2}\bigg(\frac{F^{12}}{v_{22}}(vF^v)^2-F^{12}\cdot vF^{22,v}\cdot vF^v-F^{22}\cdot vF^{12,v}\cdot vF^v+F^{12}F^{22}\cdot v^2F^{v,v}\bigg)\\
\sim&~\frac{2}{(vF^v)^2}\bigg(\frac{F^{12}}{v_{22}}(vF^v)^2-F^{12}\cdot(F^{22}+v^2(v_{11}+v_{33}))\cdot vF^v-F^{12}F^{22}\cdot vF^v\\
&\hspace{1.8cm}+F^{12}F^{22}\cdot2v^2(v_{11}+v_{33})v_{22}\bigg)\\
\sim&~\frac{2}{(vF^v)^2}\frac{F^{12}}{v_{22}}\left((vF^v)^2-v^2(v_{11}+v_{33})v_{22}\cdot vF^v-2F^{22}v_{22}\cdot (vF^v-v^2(v_{11}+v_{33})v_{22})\right)\\
\sim&~\frac{2}{(vF^v)^2}\frac{F^{12}}{v_{22}}\left(vF^v-2F^{22}v_{22}\right).
\end{align*}

$\bullet$ The coefficient of $2v_{224}v_{234}$ is
\begin{align*}
&\frac{2F^{23}}{v_{22}}-\frac{2F^{23}F^{22,v}}{F^v}-\frac{2F^{22}F^{23,v}}{F^v}+\frac{2F^{22}F^{23}}{(F^v)^2}F^{v,v}\\
\sim&~\frac{2}{(vF^v)^2}\frac{F^{23}}{v_{22}}\left(vF^v-2F^{22}v_{22}\right).
\end{align*}

$\bullet$ The coefficient of $2v_{334}v_{124}$ is
\begin{align*}
&-\frac{2F^{12}F^{33,v}}{F^v}-\frac{2F^{33}F^{12,v}}{F^v}+\frac{2F^{12}F^{33}}{(F^v)^2}F^{v,v}\\
\sim&-\frac{2F^{12}}{(vF^v)^2}(v^2v_{22}\cdot vF^v+2F^{33}).
\end{align*}

$\bullet$ The coefficient of $2v_{334}v_{234}$ is
\begin{align*}
&\frac{2F^{23}}{v_{33}}-\frac{2F^{23}F^{33,v}}{F^v}-\frac{2F^{33}F^{23,v}}{F^v}+\frac{2F^{23}F^{33}}{(F^v)^2}F^{v,v}\\
\sim&~\frac{2}{(vF^v)^2}\frac{F^{23}}{v_{33}}(v^2v_{11}v_{22}\cdot vF^v+vF^v-2F^{33}v_{33}).
\end{align*}

$\bullet$ The coefficient of $2v_{124}v_{234}$ is
\begin{align*}
&-\frac{4F^{23}F^{12,v}}{F^v}-\frac{4F^{12}F^{23,v}}{F^v}+\frac{4F^{12}F^{23}}{(F^v)^2}F^{v,v}\\
=&-\frac{4}{(vF^v)^2}\hspace{1pt}(F^{23}\cdot vF^{12,v}\cdot vF^v+F^{12}\cdot vF^{23,v}\cdot vF^v-F^{12}F^{23}\cdot v^2F^{v,v})\\
\sim&-\frac{8F^{12}F^{23}}{(vF^v)^2}\left(vF^v-v^2(v_{11}+v_{33})v_{22}\right)\\
\sim&-\frac{8F^{12}F^{23}}{(vF^v)^2}.
\end{align*}

$\bullet$ The coefficients of $2v_{114}v_{134}$, $2v_{224}v_{134}$, $2v_{334}v_{134}$, $2v_{134}v_{124}$ and $2v_{134}v_{234}$ are $$0, \quad 0, \quad 0, \quad \frac{2F^{23}}{v^{11}} \ \ \mathrm{and} \ \  \frac{2F^{12}}{v^{33}}, $$ respectively.

Therefore, we obtain 
\begin{align}\label{Quadratic form-f}
-\frac{(\det(A))^2(vF^v)^2}{2\sigma_3(G)} \sum_{i, j}F^{ij}\phi_{ij}\sim f(X),
\end{align}
where $ f(X)=X^TMX$ is a quadratic form with $X=(v_{114}, v_{224}, v_{334}, v_{134}, v_{124}, v_{234})^T$ as the associated vector and $M$ as the symmetric matrix defined below. 
\begin{center}
\scalebox{0.75}{
$\begin{matrix}
  & v_{114} & v_{224} & v_{334} & v_{134} & v_{124} & v_{234}\\[13pt]
v_{114} & \makecell{\frac{F^{11}}{v_{11}}(v^2v_{22}v_{33}\cdot vF^v\\ +F^{22}v_{22}+F^{33}v_{33})} & -F^{11}F^{22} & \makecell{-F^{11}(v^2v_{22}\cdot vF^v\\ +F^{11})} & 0 & \makecell{\frac{F^{12}}{v_{11}} (v^2v_{22}v_{33}\cdot vF^v\\ +vF^v-2F^{11}v_{11})} & \makecell{-F^{23}(v^2v_{22}\cdot vF^v\\ +2F^{11})}\\[18pt]
v_{224} & -F^{11}F^{22} & \makecell{\frac{F^{22}}{v_{22}}(F^{11}v_{11}\\ +F^{33}v_{33})} & -F^{22}F^{33} & 0 & \frac{F^{12}}{v_{22}}(vF^v -2F^{22}v_{22}) & \frac{F^{23}}{v_{22}}(vF^v-2F^{22}v_{22})\\[18pt]
v_{334} & \makecell{-F^{11}(v^2v_{22}\cdot vF^v\\ +F^{11})} & -F^{22}F^{33} & \makecell{\frac{F^{33}}{v_{33}}(v^2v_{11}v_{22}\cdot vF^v\\ +F^{11}v_{11}+F^{22}v_{22})} & 0 & \makecell{-F^{12}(v^2v_{22}\cdot vF^v\\ +2F^{33})} & \makecell{\frac{F^{23}}{v_{33}}(v^2v_{11}v_{22}\cdot vF^v\\ +vF^v-2F^{33}v_{33})}\\[18pt]
v_{134} & 0 & 0 & 0 & \makecell{\left(\frac{F^{33}}{v_{11}}+\frac{F^{11}}{v_{33}}\right)\\ \cdot (vF^v)^2} & \frac{F^{23}}{v_{11}}(vF^v)^2 & \frac{F^{12}}{v_{33}}(vF^v)^2\\[18pt]
v_{124} & \makecell{\frac{F^{12}}{v_{11}} (v^2v_{22}v_{33}\cdot vF^v\\ +vF^v-2F^{11}v_{11})} & \frac{F^{12}}{v_{22}}(vF^v -2F^{22}v_{22}) & \makecell{-F^{12}(v^2v_{22}\cdot vF^v\\ +2F^{33})} & \frac{F^{23}}{v_{11}}(vF^v)^2 & \makecell{\frac{v^2v_{22}v_{33}+1-F^{33}v_{33}}{v_{11}v_{22}}\\ \cdot (vF^v)^2 -4(F^{12})^2} & -4F^{12}F^{23}\\[18pt]
v_{234} & \makecell{-F^{23}(v^2v_{22}\cdot vF^v\\ +2F^{11})} & \frac{F^{23}}{v_{22}}(vF^v-2F^{22}v_{22}) & \makecell{\frac{F^{23}}{v_{33}}(v^2v_{11}v_{22}\cdot vF^v\\ +vF^v-2F^{33}v_{33})} & \frac{F^{12}}{v_{33}}(vF^v)^2 & -4F^{12}F^{23} & \makecell{\frac{v^2v_{11}v_{22}+1-F^{11}v_{11}}{v_{22}v_{33}}\\ \cdot (vF^v)^2-4(F^{23})^2}
\end{matrix}$
}
\end{center}

\vspace{0.4cm} Let us verify that the matrix $M$ is semi-positive modulo $\phi$, i.e., the corresponding quadratic form $f(X)$ satisfies $f(X)\gtrsim0$. This is highly technical and challenging.

First, we perform congruence transformations on the matrix $M$, which preserve the semi-positivity of $M$, to convert it into a block diagonal form. The specific transformations are as follows.
\begin{enumerate}
\item The first column is multiplied by $-\frac{F^{12}}{F^{11}}$ and added to the fifth column, while the first row is multiplied by $-\frac{F^{12}} {F^{11}}$ and added to the fifth row. Then, the second column is multiplied by $-\frac{F^{12}}{F^{22}}$ and added to the fifth column, while the second row is multiplied by $-\frac{F^{12}}{F^{22}}$ and added to the fifth row.

\item The third column is multiplied by $-\frac{F^{23}}{F^{11}}$ and added to the sixth column, while the third row is multiplied by $-\frac{F^{23}} {F^{11}}$ and added to the sixth row. Then, the second column is multiplied by $-\frac{F^{23}}{F^{22}}$ and added to the sixth column, while the second row is multiplied by $-\frac{F^{23}}{F^{22}}$ and added to the sixth row.
\end{enumerate}
It yields the following block diagonal matrix.
\begin{center}
\scalebox{0.95}{
$\begin{matrix}
& v_{114} & v_{224} & v_{334} & v_{134} & v_{124} & v_{234}\\[8pt]
v_{114} & \makecell{\frac{F^{11}}{v_{11}}(v^2v_{22}v_{33}\cdot vF^v\\ +F^{22}v_{22}+F^{33}v_{33})} & -F^{11}F^{22} & \makecell{-F^{11}(v^2v_{22}\cdot vF^v\\ +F^{11})} & 0 & 0 & 0\\[13pt]
v_{224} & -F^{11}F^{22} & \makecell{\frac{F^{22}}{v_{22}}(F^{11}v_{11}\\ +F^{33}v_{33})} & -F^{22}F^{33} & 0 & 0 & 0\\[13pt]
v_{334} & \makecell{-F^{11}(v^2v_{22}\cdot vF^v\\ +F^{11})} & -F^{22}F^{33} & \makecell{\frac{F^{33}}{v_{33}}(v^2v_{11}v_{22}\cdot vF^v\\ +F^{11}v_{11}+F^{22}v_{22})} & 0 & 0 & 0\\[13pt]
v_{134} & 0 & 0 & 0 & \makecell{\left(\frac{F^{33}}{v_{11}}+\frac{F^{11}}{v_{33}}\right)\\ \cdot (vF^v)^2} & \frac{F^{23}}{v_{11}}(vF^v)^2 & \frac{F^{12}}{v_{33}}(vF^v)^2\\[13pt]
v_{124} & 0 & 0 & 0 & \frac{F^{23}}{v_{11}}(vF^v)^2 & b_{11} & b_{12}\\[13pt]
v_{234} & 0 & 0 & 0 & \frac{F^{12}}{v_{33}}(vF^v)^2 & b_{21} & b_{22}
\end{matrix}$
}
\end{center}
Here
\begin{align*}
&b_{11}=\bigg(\frac{F^{11}}{v_{22}}+\frac{F^{22}}{v_{11}}-v^2\bigg)\cdot(vF^v)^2-\bigg(\frac{v^2v_{22}v_{33}+1}{F^{11}v_{11}}+\frac{1}{F^{22}v_{22}}\bigg)(F^{12})^2\cdot vF^v, \\
&b_{12}=b_{21}=F^{12}F^{23} \bigg(\frac{v^2v_{22}}{F^{11}} -\frac{1}{F^{22}v_{22}} \bigg)\cdot vF^v,
\intertext{and}
&b_{22}=\bigg(\frac{F^{22}}{v_{33}}+\frac{F^{33}}{v_{22}}-v^2\bigg)\cdot(vF^v)^2-\bigg(\frac{1}{F^{22}v_{22}}+\frac{v^2v_{11}v_{22}+1}{F^{33}v_{33}}\bigg)(F^{23})^2\cdot vF^v.
\end{align*}

\textbf{Claim 2. } The matrix $$M_1=
\begin{pmatrix}
\makecell{\frac{F^{11}}{v_{11}}(v^2v_{22}v_{33}\cdot vF^v\\ +F^{22}v_{22}+F^{33}v_{33})} & -F^{11}F^{22} & \makecell{-F^{11}(v^2v_{22}\cdot vF^v\\ +F^{11})}\\[13pt]
-F^{11}F^{22} & \makecell{\frac{F^{22}}{v_{22}}(F^{11}v_{11}\\ +F^{33}v_{33})} & -F^{22}F^{33}\\[13pt]
\makecell{-F^{11}(v^2v_{22}\cdot vF^v\\ +F^{11})} & -F^{22}F^{33} & \makecell{\frac{F^{33}}{v_{33}}(v^2v_{11}v_{22}\cdot vF^v\\ +F^{11}v_{11}+F^{22}v_{22})}
\end{pmatrix}$$ is semi-positive definite modulo $\phi$.

\textbf{Claim 3. } The matrix $$M_2=
\begin{pmatrix}
\makecell{\left(\frac{F^{33}}{v_{11}}+\frac{F^{11}}{v_{33}}\right)\\ \cdot (vF^v)^2} & \frac{F^{23}}{v_{11}}(vF^v)^2 & \frac{F^{12}}{v_{33}}(vF^v)^2\\[13pt]
\frac{F^{23}}{v_{11}}(vF^v)^2 & b_{11} & b_{12}\\[13pt]
\frac{F^{12}}{v_{33}}(vF^v)^2 & b_{21} & b_{22}
\end{pmatrix}$$ is positive definite modulo $\phi$.

Once \textbf{Claim 2} and \textbf{Claim 3} are proved, we will know that the symmetric matrix $M$ is semi-positive definite modulo $\phi$. That is, the corresponding quadratic form $f(X)$ satisfies $$f(X)=X^TMX\gtrsim0. $$ Combining this with formula \eqref{Quadratic form-f}, we conclude the proof of Theorem \ref{Theorem-rank3}.
\end{proof}

Our remaining task is to prove \textbf{Claim 2} and \textbf{Claim 3}.

\textbf{Proof of Claim 2.} Let us compute the leading principal minors of $M_1$, denoted by $P_1(M_1)$,  $P_2(M_1)$ and  $P_3(M_1)$.

$\bullet$ $P_1(M_1)$. It is evident that
$$P_1(A_1)=\frac{F^{11}}{v_{11}}(v^2v_{22}v_{33}\cdot vF^v +F^{22}v_{22} +F^{33}v_{33})>0. $$

$\bullet$ $P_2(M_1)$. Note that $F^{33}=F^{11}$. By \eqref{Relation-vF^v1-3} and \eqref{Relation-vF^v2-3}, we have
\begin{align*}
P_2(M_1)=&~\frac{F^{11}}{v_{11}}(v^2v_{22}v_{33}\cdot vF^v +F^{22}v_{22} +F^{33}v_{33})\cdot \frac{F^{22}}{v_{22}}(F^{11}v_{11}+F^{33}v_{33}) -(-F^{11}F^{22})^2\\
=&~\frac{(F^{11})^2F^{22}}{v_{11}v_{22}}\left((v^2v_{22}v_{33}\cdot vF^v +F^{22}v_{22}+F^{33}v_{33})\cdot (v_{11}+v_{33})-F^{22}v_{11}v_{22}\right)\\
=&~\frac{(F^{11})^2F^{22}}{v_{11}v_{22}}\left(v^2(v_{11}+v_{33})v_{22}v_{33}\cdot vF^v+F^{22}v_{22}v_{33}+F^{33}v_{11}v_{33}+F^{33}v_{33}^2\right)\\
\sim&~\frac{(F^{11})^2F^{22}v_{33}}{v_{11}v_{22}}(v^2(v_{11}+v_{33})v_{22}v_{33}\cdot vF^v+vF^v)\\
\sim&~\frac{(F^{11})^2F^{22}v_{33}}{v_{11}v_{22}}(vF^v)^2>0.
\end{align*}

$\bullet$ $P_3(M_1)$. To compute $\det(M_1)$, we perform elementary column operations on $M_1$. Specifically, we multiply the first column by $\frac{v_{11}+v_{33}}{v_{22}}$ and add it to the second column, then multiply the first column by $-1$ and add it to the third column. By \eqref{Relation-vF^v1-3} and \eqref{Relation-vF^v2-3}, we obtain the transformed matrix $$
\begin{pmatrix}
\makecell{\frac{F^{11}}{v_{11}}(v^2v_{22}v_{33}\cdot vF^v\\ +F^{22}v_{22}+F^{33}v_{33})} & \frac{F^{11}v_{33}}{v_{11}v_{22}}\hspace{1pt}(vF^v)^2 & -\frac{F^{11}}{v_{11}}\hspace{1pt}(vF^v)^2\\[13pt]
-F^{11}F^{22} & 0 & 0\\[13pt]
\makecell{-F^{11}(v^2v_{22}\cdot vF^v\\ +F^{11})} & -\frac{F^{11}}{v_{22}}\hspace{1pt}(vF^v)^2 & \frac{F^{11}}{v_{33}}\hspace{1pt}(vF^v)^2
\end{pmatrix}. $$
At this point, the second and third columns become proportional. Therefore, $$P_3(M)=\det(M_1)=0. $$ 

By Lemma \ref{Lemma-linear algebra}, the matrix $M_1$ is semi-positive definite modulo $\phi$. \qed 

\vspace{0.5cm}\textbf{Proof of Claim 3.} Along the same route, we compute the leading principal minors of $M_2$, denoted by $P_1(M_2)$,  $P_2(M_2)$ and  $P_3(M_2)$.

$\bullet$ $P_1(M_2)$. Clearly, $$P_1(M_2)=\bigg(\frac{F^{33}}{v_{11}}+\frac{F^{11}}{v_{33}}\bigg)\cdot(vF^v)^2>0. $$

$\bullet$ $P_2(M_2)$. By \eqref{Relation-Equation-3} and \eqref{Relation-Fij-b}, we have
\begin{equation}\label{Relation-elliptic coefficients-3}
\begin{aligned}
F^{11}F^{22}-(F^{12})^2-(F^{23})^2
\sim&~v^4(v_{11}+v_{33})v_{22}\\
&+v^3((v_2^2+v_4^2)(v_{11}+v_{33})+(v_1^2+v_3^2)v_{22})\\
&+v^2(v_1^2+v_3^2)(v_2^2+v_4^2)\\
&-v^2(v_1v_2+v_3v_4)^2-v^2(v_1v_4-v_2v_3)^2\\
\sim&~v^2,
\end{aligned}
\end{equation}
Using equations \eqref{Relation-vF^v1-3} and \eqref{Relation-vF^v2-3}, we obtain
\begin{align*}
P_2(M_2)=&\bigg(\frac{F^{33}}{v_{11}}+\frac{F^{11}}{v_{33}}\bigg)\cdot(vF^v)^2\cdot b_{11}-\bigg(\frac{F^{23}}{v_{11}}(vF^v)^2\bigg)^2\\
=&\bigg(\frac{F^{33}}{v_{11}}+\frac{F^{11}}{v_{33}}\bigg)\cdot(vF^v)^2\cdot
\bigg(\bigg(\frac{F^{11}}{v_{22}}+\frac{F^{22}}{v_{11}}-v^2\bigg)\cdot(vF^v)^2\\
&\hspace{1.2cm}-\bigg(\frac{v^2v_{22}v_{33}+1}{F^{11}v_{11}}+\frac{1}{F^{22}v_{22}}\bigg)(F^{12})^2\cdot vF^v\bigg) -\bigg(\frac{F^{23}}{v_{11}}(vF^v)^2\bigg)^2\\
=&~\frac{(v_{11}+v_{33})F^{11}}{v_{11}v_{33}}\frac{(vF^v)^3}{F^{11}v_{11}F^{22}v_{22}}\cdot\bigg((F^{11}v_{11}+F^{22}v_{22}-v^2v_{11}v_{22})\cdot F^{11}F^{22}\cdot vF^v\\
&\hspace{0.5cm}-((v^2v_{22}v_{33}+1)\cdot F^{22}v_{22}+F^{11}v_{11})\cdot(F^{12})^2\bigg) -\frac{(F^{23})^2}{v_{11}^2}\cdot(vF^v)^4\\
\sim&~\frac{(v_{11}+v_{33})F^{11}}{v_{11}v_{33}}\frac{(vF^v)^3}{F^{11}v_{11}F^{22}v_{22}}\cdot\bigg((v^2v_{22}v_{33}+1-F^{33}v_{33})\cdot F^{11}F^{22}\cdot vF^v\\
&\hspace{0.5cm}-(v^2v_{22}v_{33}\cdot F^{22}v_{22}+vF^v-F^{33}v_{33})\cdot(F^{12})^2\bigg) -\frac{(F^{23})^2}{v_{11}^2}\cdot(vF^v)^4\\
=&~\frac{F^{11}}{v_{11}^2}\cdot v^2(v_{11}+v_{33})\cdot (vF^v)^4+\frac{v_{11}+v_{33}}{v_{11}v_{33}}\frac{F^{11}}{F^{11}v_{11}F^{22}v_{22}}(F^{11}F^{22}-(F^{12})^2)\cdot(vF^v)^4\\
&-\frac{v_{11}+v_{33}}{v_{11}^2v_{22}}F^{11}F^{33}\cdot(vF^v)^4-\frac{1}{v_{11}^2}\cdot v^2(v_{11}+v_{33})v_{22}\cdot(F^{12})^2\cdot(vF^v)^3\\
&+\frac{v_{11}+v_{33}}{v_{11}}\frac{F^{11}F^{33}}{F^{11}v_{11}F^{22}v_{22}}\cdot(F^{12})^2\cdot(vF^v)^3-\frac{1}{v_{11}^2}(F^{23})^2\cdot(vF^v)^4.
\end{align*}
We continue to simplify $P_2(M_2)$. Using equations \eqref{Relation-vF^v1-3} and \eqref{Relation-vF^v2-3} again, we have
\begin{align*} P_2(M_2)\sim&~\frac{F^{11}}{v_{11}^2}\cdot\frac{vF^v-1}{v_{22}}\cdot(vF^v)^4+\frac{v_{11}+v_{33}}{v_{11}^2v_{22}v_{33}F^{22}}(F^{11}F^{22}-(F^{12})^2)\cdot(vF^v)^4\\
&-\frac{vF^v-F^{22}v_{22}}{v_{11}^2v_{22}}F^{11}\cdot(vF^v)^4-\frac{1}{v_{11}^2}\cdot(vF^v-1)\cdot(F^{12})^2\cdot(vF^v)^3\\
&+\frac{vF^v-F^{22}v_{22}}{v_{11}^2v_{22}F^{22}}\cdot(F^{12})^2\cdot(vF^v)^3-\frac{1}{v_{11}^2}(F^{23})^2\cdot(vF^v)^4\\
=&-\frac{1}{v_{11}^2v_{22}F^{22}}(F^{11}F^{22}-(F^{12})^2)\cdot(vF^v)^4+\frac{v_{11}+v_{33}}{v_{11}^2v_{22}v_{33}F^{22}}(F^{11}F^{22}-(F^{12})^2)\cdot(vF^v)^4\\
&+\frac{1}{v_{11}^2}(F^{11}F^{22}-(F^{12})^2-(F^{23})^2)\cdot(vF^v)^4\\
=&~\frac{1}{v_{11}v_{22}v_{33}F^{22}}(F^{11}F^{22}-(F^{12})^2)\cdot(vF^v)^4+\frac{1}{v_{11}^2}(F^{11}F^{22}-(F^{12})^2-(F^{23})^2)\cdot(vF^v)^4.
\end{align*}
From \eqref{Relation-elliptic coefficients-3}, it follows that
\begin{align*}
P_2(M_2)\sim&~\frac{v^2+(F^{23})^2}{v_{11}v_{22}v_{33}F^{22}}\cdot(vF^v)^4+\frac{v^2}{v_{11}^2}\cdot(vF^v)^4>0.
\end{align*}

$\bullet$ $P_3(M_2)$. Let us compute the determinant of $M_2$. As before, we perform elementary column operations on $M_2$, which do not change the determinant of $M_2$. Specifically, we multiply the first column by $-\frac{F^{23}}{F^{11}}\frac{v_{33}}{v_{11}+v_{33}}$ and add it to the second column, then multiply the first column by $-\frac{F^{12}}{F^{11}}\frac{v_{11}}{v_{11}+v_{33}}$ and add it to the third column. These operations yield the transformed matrix $$
\begin{pmatrix}
\frac{(v_{11}+v_{33})F^{11}}{v_{11}v_{33}}\cdot(vF^v)^2 & 0 & 0\\[13pt]
\frac{F^{23}}{v_{11}}\cdot(vF^v)^2 & b_{11}-\frac{v_{33}}{v_{11}+v_{33}} \frac{(F^{23})^2}{F^{11}v_{11}}\cdot(vF^v)^2 & b_{12} -\frac{1}{v_{11}+v_{33}}\frac{F^{12}F^{23}}{F^{11}}\cdot(vF^v)^2\\[13pt]
\frac{F^{12}}{v_{33}}(vF^v)^2 & b_{21}-\frac{1}{v_{11}+v_{33}} \frac{F^{12}F^{23}}{F^{11}}\cdot(vF^v)^2 & b_{22}-\frac{v_{11}} {v_{11}+v_{33}}\frac{(F^{12})^2}{F^{33}v_{33}}\cdot(vF^v)^2
\end{pmatrix}. $$ Note that $F^{11}=F^{33}$. By \eqref{Relation-vF^v1-3}, \eqref{Relation-vF^v2-3} and \eqref{Relation-elliptic coefficients-3}, we have
\begin{align*}
&b_{11}-\frac{v_{33}}{v_{11}+v_{33}}\frac{(F^{23})^2}{F^{11}v_{11}}\cdot(vF^v)^2\\
=&\bigg(\frac{F^{11}}{v_{22}}+\frac{F^{22}}{v_{11}}-v^2\bigg)\cdot(vF^v)^2-\bigg(\frac{v^2v_{22}v_{33}+1}{F^{11}v_{11}}+\frac{1}{F^{22}v_{22}}\bigg)(F^{12})^2\cdot vF^v\\
&-\frac{v_{33}}{v_{11}+v_{33}}\frac{(F^{23})^2}{F^{11}v_{11}}\cdot(vF^v)^2\\
\sim&~\frac{vF^v}{(v_{11}+v_{33})F^{11}v_{11}F^{22}v_{22}}\bigg((-F^{33}v_{33}+v^2v_{22}v_{33}+1)\cdot(v_{11}+v_{33})\cdot F^{11}F^{22}\cdot vF^v\\
&\hspace{4.5cm}-(v^2v_{22}v_{33}\cdot F^{22}v_{22}+vF^v-F^{33}v_{33})\cdot (v_{11}+v_{33})\cdot(F^{12})^2\\
&\hspace{4.5cm}-v_{22}v_{33}F^{22}\cdot(F^{23})^2\cdot vF^v\bigg)\\
=&~\frac{vF^v}{(v_{11}+v_{33})F^{11}v_{11}F^{22}v_{22}}\bigg(-F^{33}v_{33}\cdot F^{11}(v_{11}+v_{33})\cdot F^{22}\cdot vF^v\\
&\hspace{4.5cm}+v_{33}\cdot v^2(v_{11}+v_{33})v_{22}\cdot F^{11}F^{22}\cdot vF^v\\
&\hspace{4.5cm}+(v_{11}+v_{33})\cdot F^{11}F^{22}\cdot vF^v\\
&\hspace{4.5cm}-v_{22}v_{33}F^{22}\cdot v^2(v_{11}+v_{33})v_{22}\cdot(F^{12})^2\\
&\hspace{4.5cm}-(v_{11}+v_{33})\cdot(F^{12})^2\cdot vF^v+v_{33}\cdot F^{33}(v_{11}+v_{33})\cdot(F^{12})^2\\
&\hspace{4.5cm}-v_{22}v_{33}F^{22}\cdot(F^{23})^2\cdot vF^v\bigg)\\
\sim&~\frac{vF^v}{(v_{11}+v_{33})F^{11}v_{11}F^{22}v_{22}}\bigg(-F^{33}v_{33}\cdot(vF^v-F^{22}v_{22})\cdot F^{22}\cdot vF^v\\
&\hspace{4.5cm}+v_{33}\cdot(vF^v-1)\cdot F^{11}F^{22}\cdot vF^v\\
&\hspace{4.5cm}+(v_{11}+v_{33})\cdot(F^{11}F^{22}-(F^{12})^2)\cdot vF^v\\
&\hspace{4.5cm}-v_{22}v_{33}F^{22}\cdot(vF^v-1)\cdot(F^{12})^2+v_{33}\cdot(vF^v-F^{22}v_{22})\cdot(F^{12})^2\\
&\hspace{4.5cm}-v_{22}v_{33}F^{22}\cdot(F^{23})^2\cdot vF^v\bigg)\\
=&~\frac{(vF^v)^2}{(v_{11}+v_{33})F^{11}v_{11}F^{22}v_{22}}\bigg(v_{22}v_{33}F^{22}\cdot(F^{11}F^{11}-(F^{12})^2-(F^{23})^2)+v_{11}\cdot(F^{11}F^{22}-(F^{12})^2)\bigg)\\
\sim&~\frac{(vF^v)^2}{(v_{11}+v_{33})F^{11}v_{11}F^{22}v_{22}}\bigg(v^2v_{22}v_{33}F^{22}+v_{11}(v^2+(F^{23})^2)\bigg).
\end{align*}
Analogously, we have the following simplifications
\begin{align*}
&b_{12}-\frac{1}{v_{11}+v_{33}} \frac{F^{12}F^{23}}{F^{11}}\cdot(vF^v)^2 =b_{21}-\frac{1}{v_{11}+v_{33}} \frac{F^{12}F^{23}}{F^{11}}\cdot(vF^v)^2\\
=&~F^{12}F^{23} \bigg(\frac{v^2v_{22}}{F^{11}}-\frac{1}{F^{22}v_{22}} \bigg)\cdot vF^v-\frac{1}{v_{11}+v_{33}} \frac{F^{12}F^{23}}{F^{11}}\cdot (vF^v)^2\\
=&~\frac{F^{12}F^{23}}{(v_{11}+v_{33})v_{22}F^{11}F^{22}}\cdot vF^v\cdot \bigg(v^2(v_{11}+v_{33})v_{22}\cdot F^{22}v_{22}-F^{11}(v_{11}+v_{33})-F^{22}v_{22}\cdot vF^v\bigg)\\
\sim&~\frac{F^{12}F^{23}}{(v_{11}+v_{33})v_{22}F^{11}F^{22}}\cdot vF^v\cdot \bigg((vF^v-1)\cdot F^{22}v_{22}-F^{11}(v_{11}+v_{33})-F^{22}v_{22}\cdot vF^v\bigg)\\
\sim&-\frac{F^{12}F^{23}}{(v_{11}+v_{33})v_{22}F^{11}F^{22}}\cdot(vF^v)^2,
\end{align*}
and
\begin{align*}
&b_{22}-\frac{v_{11}}{v_{11}+v_{33}}\frac{(F^{12})^2}{F^{33}v_{33}}\cdot(vF^v)^2\\
\sim&~\frac{(vF^v)^2}{(v_{11}+v_{33})F^{22}v_{22}F^{33}v_{33}}\bigg(v^2v_{11}v_{22}F^{22}+v_{33}(v^2+(F^{12})^2)\bigg).
\end{align*}
Thus,
\begin{align*}
P_3(M_2)=\det(M_2)
\sim&~\frac{(v_{11}+v_{33})F^{11}}{v_{11}v_{33}}\cdot(vF^v)^2\\
&\cdot\left(\frac{(vF^v)^4}{(v_{11}+v_{33})^2F^{11}v_{11}(F^{22}v_{22})^2F^{33}v_{33}}\right.\\
&\hspace{0.4cm}\cdot\left(v^2v_{22}v_{33}F^{22}+v_{11}(v^2+(F^{23})^2)(v^2v_{11}v_{22}F^{22}+v_{33}(v^2+(F^{12})^2))\right)\\
&\hspace{0.4cm}\left.-\frac{(vF^v)^4}{(v_{11}+v_{33})^2v_{22}^2(F^{11}F^{22})^2}(F^{12})^2(F^{23})^2\right)\\
=&~\frac{(v_{11}+v_{33})F^{11}}{v_{11}v_{33}}\cdot(vF^v)^2\cdot\frac{(vF^v)^4}{(v_{11}+v_{33})^2F^{11}v_{11}(F^{22}v_{22})^2F^{33}v_{33}}\\
&\cdot\bigg(\left(v^2v_{22}v_{33}F^{22}+v_{11}(v^2+(F^{23})^2)(v^2v_{11}v_{22}F^{22}+v_{33}(v^2+(F^{12})^2))\right)\\
&\hspace{0.4cm}-v_{11}v_{33}(F^{12})^2(F^{23})^2\bigg)\\
=&~\frac{(vF^v)^6}{(v_{11}+v_{33})(v_{11}v_{22}v_{33})^2F^{11}(F^{22})^2}\\
&\cdot\bigg(v^4v_{11}v_{33}\cdot(F^{22}v_{22})^2+v^2v_{33}^2\cdot F^{22}v_{22}\cdot(v^2+(F^{12})^2)\\
&\hspace{0.4cm}+v^2v_{11}^2\cdot F^{22}v_{22}\cdot(v^2+(F^{23})^2) +v^4v_{11}v_{33}+v^2v_{11}v_{33}\cdot((F^{12})^2+(F^{23})^2)\bigg)\\
>&~0.
\end{align*}
Therefore, the matrix $M_2$ is positive definite modulo $\phi$. \qed 

At the end of this section, we will use Theorems \ref{Theorem-rank2} and \ref{Theorem-rank3} to establish the constant rank theorem.

\begin{corollary}[Constant rank theorem]\label{Corollary-CRT}
Let $\Omega\subset\mathbb{R}^4$ be a domain, and let $v\in C^4(\Omega)$ be a convex solution of \eqref{Equ-complexMA-v}. Then the Hessian $\nabla^2v$ has constant rank throughout $\Omega$.
\end{corollary}
\begin{proof}
Define $$l=\min_{x\in\Omega}\mathrm{rank}(\nabla^2v(x)). $$ As discussed in Section \ref{Section 1.3}, $l$ can only be $2$ or $3$. Suppose there exists some point $x_0\in\Omega$ such that $\mathrm{rank}(\nabla^2v(x_0))=l$. We will show that $\mathrm{rank}(\nabla^2v(x))=l$ for all $x\in\Omega$. By \eqref{rank D^2v=rank K+2}, namely $$\mathrm{rank}(\nabla^2v(x))=\mathrm{rank}(K(x))+2,$$ it suffices to study the constant rank property of $K(x)$. Define the set $$\Omega^\ast=\{x\in\Omega: \mathrm{rank}(K(x))=l-2\}. $$ It is clear that $\Omega^\ast$ is relatively closed in $\Omega$. By the strong maximum principle, Theorems \ref{Theorem-rank2} and \ref{Theorem-rank3} imply that $\Omega^\ast$ is also open. Hence, $\Omega^\ast=\Omega$, meaning that $K$ has constant rank $l-2$ on all of $\Omega$, which in turn implies that $\nabla^2v$ has constant rank $l$ throughout $\Omega$.

This completes the proof of Corollary \ref{Corollary-CRT}.
\end{proof}

\section{Proof of Theorem \ref{Theorem-Main}}\label{Section 5}

In this section, we prove Theorem \ref{Theorem-Main} by the continuity method. With Corollary \ref{Corollary-CRT} just established, it is well known that the proof of Theorem \ref{Theorem-Main} follows standard arguments; see, for example, Caffarelli and Friedman \cite{Caffarelli-Friedman1985}, as well as \cite{Korevaar-Lewis1987, Liu-Ma-Xu2010, Ma-Xu2008, Zhang-Zhou2023}.

Before proceeding, we require the convexity estimates near the boundary for the function $v=-\sqrt{-u/2}$, where $u$ satisfies the Dirichlet problem \eqref{Equ-CMA1}. In the following proposition, we take $\alpha=1/2$.
\begin{proposition}[Convexity estimates near the boundary, see Korevaar \cite{Korevaar1983b}, Caffarelli and Friedman \cite{Caffarelli-Friedman1985}] \label{Proposition-BCE}
Let $\Omega\subset\mathbb{R}^n$ be a bounded, smooth, strictly convex domain, and let $u\in C^2(\bar{\Omega})$ satisfy $$u<0 \quad \mathrm{in}\ \Omega, \quad u=0 \quad \mathrm{on}\ \partial\Omega, \quad |\nabla u|>0 \quad \mathrm{on}\ \partial\Omega. $$ Then there exists $\varepsilon>0$ such that the function $v=-\log(-u)$ or $v=-(-u)^\alpha$ $(0<\alpha<1)$ is strictly convex in the boundary strip $\Omega\setminus\bar\Omega_\varepsilon$, where $$\Omega_\varepsilon=\{x\in\Omega: \mathrm{d}(x, \partial\Omega)>\varepsilon\}. $$
\end{proposition}

We now moved on to the proof of Theorem \ref{Theorem-Main}.
\begin{proof}[Proof of Theorem \ref{Theorem-Main}]
If $\Omega_0=B_1$, the unit ball in $\mathbb{R}^4$, then the unique solution to the transformed Dirichlet problem \eqref{Equ-Dirichlet problem-real} is $$v_0(x)=-\sqrt{\frac{1-|z_1|^2-|z_2|^2}{4}} \quad \mathrm{for} \ x\in B_{1}, $$ where $x=(x_1, x_2, x_3, x_4)$, $z_1=x_1+\sqrt{-1}x_3$, $z_2=x_2+\sqrt{-1}x_4$, $|z_1|^2=x_1^2+x_3^2$ and $|z_2|^2=x_2^2+x_4^2$. Clearly, the function $v_0$ is strictly convex in $B_1$. For an arbitrary bounded, smooth, strictly convex domain $\Omega_1=\Omega$, define $$\Omega_t=(1-t)B_1+t\Omega \quad \mathrm{for} \ t\in(0, 1). $$ By the theory of convex bodies,  the family $\{\Omega_t\}_{t\in(0, 1)}$ provides a smooth deformation of $B_1$ into $\Omega$, preserving boundedness, smoothness, and strict convexity. For further details, see Schneider's excellent book \cite{Schneider2014}. 

For $t\in(0, 1)$, assume that the function $v_t$ solves the Dirichlet problem $$
\begin{cases}
F(\nabla^2v, \nabla v, v)=0 & \mathrm{in}\ \Omega_t, \\
\hspace{2.15cm}v=0 & \mathrm{on}\ \partial\Omega_t,
\end{cases} $$ 
where
\begin{align*}
F(\nabla^2v, \nabla v, v)=&~v^2\left[(v_{11}+v_{33})(v_{22}+v_{44})-(v_{12}+v_{34})(v_{21}+v_{43})-(v_{14}-v_{32})(v_{41}-v_{23})\right]\\
&+v\left[(v_2^2+v_4^2)(v_{11}+v_{33})+(v_1^2+v_3^2)(v_{22}+v_{44})\right.\\
&\hspace{0.8cm}\left.-(v_1v_2+v_3v_4)(v_{12}+v_{34}+v_{21}+v_{43})\right.\\
&\hspace{0.8cm}\left.-(v_1v_4-v_2v_3)(v_{14}-v_{32}+v_{41}-v_{23})\right]-1.
\end{align*}
According to the a priori estimates for the Dirichlet problem of complex Monge-Amp\`{e}re equation established by Caffarelli, Kohn, Nirenberg and Spruck \cite{CKNS1985}, we have uniform bounds for $\|v_t\|_{C^3(\bar\Omega_t)}$, depending only on the geometry of $\Omega$ and independent of $t$. On the one hand, since $v_0$ is strictly convex, it follows that $v_t$ remains strictly convex for all $t\in(0, \delta)$ for some sufficiently small $\delta>0$. On the other hand, if $v_t$ is strictly convex in $\Omega_t$ for all $t\in(0, t^\ast)$ with $t^\ast\in(0, 1)$, then $v_{t^\ast}$ must be convex in $\Omega_{t^\ast}$.

Now, suppose that $v_1$ is not strictly convex in $\Omega_1$. Then there exists some $t_0\in(0, 1)$ at which $v_{t_0}$ is convex but not strictly convex for the first time. By Corollary \ref{Corollary-CRT} (Constant rank theorem), the Hessian $\nabla^2v_{t_0}$ is degenerate at every point of $\Omega_{t_0}$. However, Proposition \ref{Proposition-BCE} (Convexity estimates near the boundary) implies that $\nabla^2v_{t_0}$ is of full rank in the boundary strip $\Omega_{t_0}\setminus\bar\Omega_{t_0, \varepsilon}$ for some $\varepsilon>0$, which depends only on the geometry of $\partial\Omega_1$. This leads to a contradiction. Therefore, $v_1$ must be strictly convex in $\Omega_1$. Thus, we complete the proof of Theorem \ref{Theorem-Main}.
\end{proof}

\begin{remark}\label{Remark}
We provide the following example to demonstrate that the convexity exponent $1/2$ is sharp. Similar examples can be found in Kennington \cite{Kennington1985}, Ma and Xu \cite{Ma-Xu2008}, and in the recent work of Zhang and Zhou \cite{Zhang-Zhou2023}, respectively.

Consider the infinite open cone $$\Gamma=\{(x_1, x_2, x_3, x_4)\in\mathbb{R}^4: \sqrt{(x_1)^2+(x_2)^2+(x_3)^2}<ax_4\}, $$where $a\in(0, 1)$. Let $\Omega\subset\Gamma$ be a bounded convex domain, with $y=(0, 0, 0, 0)\in\partial\Omega$ and $z\in(0, 0, 0, 1)\in\Omega$. Define the function $$w(x)=\frac{A}{2}(x_1^2+x_2^2+x_3^2-a^2x_4^2), $$ where $A>0$ is chosen such that $2A^2(1-a^2)=1$. Thus, the function $w$ satisfies $$
\begin{cases}
(w_{11}+w_{33})(w_{22}+w_{44})-(w_{12}+w_{34})(w_{21}+w_{43})-(w_{14}-w_{32})(w_{41}-w_{23})=1 & \mathrm{in}\ \Omega, \\
\hspace{12.4cm}w\leq 0 & \mathrm{on}\ \partial\Omega.
\end{cases}$$
Suppose that $u$ solves the Dirichlet problem $$\begin{cases}
(u_{11}+u_{33})(u_{22}+u_{44})-(u_{12}+u_{34})(u_{21}+u_{43})-(u_{14}-u_{32})(u_{41}-u_{23})=1 & \mathrm{in}\ \Omega, \\
\hspace{11.8cm}u= 0 & \mathrm{on}\ \partial\Omega.
\end{cases} $$	By the comparison principle, we have $w\leq u$ in $\Omega$. In particular, for $t\in(0, 1)$, it follows that $-Aa^2t^2/2=w(tz)\leq u(tz)\leq 0$. This implies $\lim_{t\rightarrow 0+}t^{-1/\alpha}u(tz)=0$ provided $\alpha>1/2$. Consequently, this shows that $-(-u)^\alpha$ is not convex when $\alpha>1/2$.
	
In fact, if $-(-u)^\alpha$ is convex for some $\alpha>1/2$, then we have $$-(-u)^\alpha((1-t)y+tz)\leq-(1-t)(-u)^\alpha(y) -t(-u)^\alpha(z)=-t(-u)^\alpha(z). $$ Rearranging this inequality yields $t^{-1/\alpha}u(tz)\leq u(z)<0$, which contradicts the fact that $\lim_{t\rightarrow 0+}t^{-1/\alpha}u(tz)=0$. This confirms the sharpness of the convexity exponent $1/2$.
\end{remark}

\section*{Acknowledgment}

We are grateful to Professor Jingchen Hu for many helpful discussions. We also thank Professor Chuanqiang Chen and Dr. Jin Yan for bringing to our attention an error in an earlier version of this manuscript.

\bibliographystyle{plain}

\bibliography{mybibliography}

\end{document}